\documentclass{birkjour}

\usepackage[latin9]{inputenc}
\usepackage{ifsym}
\usepackage{mathrsfs}
\usepackage{amstext}
\usepackage{amsthm}
\usepackage{amssymb}

\makeatletter
\numberwithin{equation}{section}
\numberwithin{figure}{section}
\theoremstyle{plain}
\newtheorem{thm}{\protect\theoremname}
\theoremstyle{plain}
\newtheorem{prop}[thm]{\protect\propositionname}
\theoremstyle{definition}
\newtheorem{defn}[thm]{\protect\definitionname}
\theoremstyle{plain}
\newtheorem{lem}[thm]{\protect\lemmaname}
\theoremstyle{plain}
\newtheorem{cor}[thm]{\protect\corollaryname}

\@ifundefined{date}{}{\date{}}

\usepackage{amscd}
\usepackage{thmdefs}
\usepackage[dvips]{epsfig}

\input{tcilatex}
\theoremstyle{definition}
\theoremstyle{remark}
\numberwithin{equation}{section}

\input tcilatex

\makeatother

\providecommand{\corollaryname}{Corollary}
\providecommand{\definitionname}{Definition}
\providecommand{\lemmaname}{Lemma}
\providecommand{\propositionname}{Proposition}
\providecommand{\theoremname}{Theorem}

\begin{document}
\title{Regular Functions on the Scaled Hypercomplex Numbers}

\author{Daniel Alpay}
\address{Chapman Univ., Dept. of Math., 1 University Dr., Orange, CA, 92866, USA}
\email{alpay@chapman.edu}
\author{Ilwoo Cho}
\address{St. Ambrose Univ., Dept. of Math. and Stat., 421 Ambrose
Hall, 518 W. Locust St., Davenport, Iowa, 52803, U. S. A.}
\email{choilwoo@sau.edu}
\keywords{Scaled Hypercomplex Rings, Scaled Hypercomplex $\mathbb{R}$-Spaces,
Scaled Hypercomplex $\mathbb{R}$-Algebras, Differential Operators.}

\begin{abstract}
In this paper, we study the regularity of $\mathbb{R}$-differentiable
functions on open connected subsets of the scaled hypercomplex numbers
$\left\{ \mathbb{H}_{t}\right\} _{t\in\mathbb{R}}$ by studying the
kernels of suitable differential operators $\left\{ \nabla_{t}\right\} _{t\in\mathbb{R}}$,
up to scales in the real field $\mathbb{R}$.
\end{abstract}

\maketitle

{\rm AMS Classification[2000]:20G20; 46S10; 47S10.}\\

Corresponding author: Daniel Alpay email: \email{alpay@chapman.edu}

\section{Introduction}

In this paper, we study differentiation on scaled hypercomplex numbers
scaled by an arbitrary quantity $t\in\mathbb{R}$. Roughly speaking,
scaled hypercomplex numbers are the ordered pairs of complex numbers
under an arbitrary fixed real number. We let $\mathbb{C}^{2}$ be
the usual 2-dimensional Hilbert space over the complex field $\mathbb{C}$,
and understand each vector $\left(a,b\right)\in\mathbb{C}^{2}$ as
a hypercomplex number $\left(a,b\right)\in\mathbb{H}_{t}$, where
$\mathbb{H}_{t}$ is the set of all $t$-scaled hypercomplex ring
for $t\in\mathbb{R}$. Algebraically, the triple,
\[
\mathbb{H}_{t}=\left(\mathbb{C}^{2},\;+,\;\cdot_{t}\right),
\]
with the usual vector addition ($+$) on $\mathbb{C}^{2}$, and the
$t$-scaled vector multiplication ($\cdot_{t}$),
\[
\left(a_{1},b_{1}\right)\cdot_{t}\left(a_{2},b_{2}\right)=\left(a_{1}a_{2}+tb_{1}\overline{b_{2}},\:a_{1}b_{2}+b_{1}\overline{a_{2}}\right),
\]
for all $\left(a_{l},b_{l}\right)\in\mathbb{C}^{2}$, for $l=1,2$,
forms a well-defined unital ring with its unity (or, the ($\cdot_{t}$)-identity)
$\left(1,0\right)$, where $\overline{z}$ mean the conjugates of
$z\in\mathbb{C}$ (e.g., see {[}2{]}).

From the Hilbert-space representation $\left(\mathbb{C}^{2},\pi_{t}\right)$
of the $t$-scaled hypercomplex ring $\mathbb{H}_{t}$, introduced
in {[}2{]}, a hypercomplex number $h=\left(a,b\right)\in\mathbb{H}_{t}$
is realized to be a (2$\times$2)-matrix, or a Hilbert-space operator
acting on $\mathbb{C}^{2}$,
\[
\pi_{t}\left(h\right)\overset{\textrm{denote}}{=}[h]_{t}\overset{\textrm{def}}{=}\left(\begin{array}{cc}
a & tb\\
\overline{b} & \overline{a}
\end{array}\right)\;\mathrm{in\;}M_{2}\left(\mathbb{C}\right),
\]
where $M_{2}\left(\mathbb{C}\right)$ is the matricial algebra (which
is $*$-isomorphic to the operator $C^{*}$-algebra $B\left(\mathbb{C}^{2}\right)$
of all bounded linear operators acting on the Hilbert space $\mathbb{C}^{2}$)
over $\mathbb{C}$, for $t\in\mathbb{R}$. The construction of such
rings $\left\{ \mathbb{H}_{t}\right\} _{t\in\mathbb{R}}$ provides
the generalized structures of well-known quaternions (e.g., {[}5{]},
{[}6{]}, {[}7{]}, {[}8{]}, {[}12{]}, {[}15{]}, {[}18{]} and {[}21{]}),
and split-quaternions (e.g., {[}4{]}, {[}9{]} and {[}14{]}). Indeed,
the ring $\mathbb{H}_{-1}$ is nothing but the noncommutative field
$\mathbb{H}$ of all quaternions, and the unital ring $\mathbb{H}_{1}$
is the ring of all split-quaternions (e.g., {[}1{]}, {[}2{]} and {[}3{]}).
The algebra, spectral theory, operator theory, and free probability
on $\left\{ \mathbb{H}_{t}\right\} _{t\in\mathbb{R}}$ are studied
in {[}1{]} and {[}2{]}, under the above representation $\left(\mathbb{C}^{2},\pi_{t}\right)$.
Different from the approaches of {[}1{]} and {[}2{]}, we study those
on $\left\{ \mathbb{H}_{t}\right\} _{t\in\mathbb{R}}$ by defining
suitable bilinear forms $\left\{ \left\langle ,\right\rangle _{t}\right\} _{t\in\mathbb{R}}$
on $\left\{ \mathbb{H}_{t}\right\} _{t\in\mathbb{R}}$, in {[}3{]}.
In the approaches of {[}3{]}, the pairs $\left\{ \left(\mathbb{H}_{t},\left\langle ,\right\rangle _{t}\right)\right\} _{t<0}$
form (definite) inner product spaces over $\mathbb{R}$, meanwhile,
the pairs $\left\{ \left(\mathbb{H}_{t},\left\langle ,\right\rangle _{t}\right)\right\} _{t\geq0}$
become indefinite semi-inner product spaces over $\mathbb{R}$, inducing
the complete semi-normed spaces $\left\{ \mathbf{X}_{t}\overset{\textrm{denote}}{=}\left(\mathbb{H}_{t},\left\Vert .\right\Vert _{t}\right)\right\} _{t\in\mathbb{R}}$,
having their semi-norms,
\[
\left\Vert h\right\Vert _{t}=\sqrt{\left|\left\langle h,h\right\rangle _{t}\right|},\;\;\forall h\in\mathbb{H}_{t},\;\forall t\in\mathbb{R},
\]
where $\left|.\right|$ is the absolute value on $\mathbb{R}$. We
call each $\mathbf{X}_{t}$, the $t$-scaled hypercomplex $\mathbb{R}$(-vector)-space,
for all $t\in\mathbb{R}$. (e.g., {[}3{]}). Meanwhile, it is considered
in {[}3{]} that each $t$-scaled hypercomplex number $h\in\mathbb{H}_{t}$
is regarded as a multiplication operator $M_{h}$ acting on $\mathbf{X}_{t}=\left(\mathbb{H}_{t},\left\Vert .\right\Vert _{t}\right)$,
\[
M_{h}\left(h'\right)=h\cdot_{t}h'\in\mathbf{X}_{t},\;\;\;\forall h'\in\mathbf{X}_{t},
\]
with its adjoint,
\[
M_{\left(a,b\right)}^{*}=M_{\left(\overline{a},-b\right)},\;\;\forall\left(a,b\right)\in\mathbb{H}_{t},
\]
where $\left(\overline{a},-b\right)\overset{\textrm{denote}}{=}\left(a,b\right)^{\dagger}$
is the hypercomplex-conjugate in $\mathbb{H}_{t}$, i.e., $M_{h}^{*}=M_{h^{\dagger}}$
in the operator space $B_{\mathbb{R}}\left(\mathbf{X}_{t}\right)$
of all bounded linear operators acting on $\mathbf{X}_{t}$ ``over
$\mathbb{R}$,'' for all scales $t\in\mathbb{R}$, which is a Banach
space equipped with the operator semi-norm,
\[
\left\Vert T\right\Vert =sup\left\{ \left\Vert Th\right\Vert _{t}:\left\Vert h\right\Vert _{t}=1\right\} ,\;\forall T\in B_{\mathbb{R}}\left(\mathbf{X}_{t}\right).
\]
Furthermore, the subset,
\[
\mathcal{M}_{t}\overset{\textrm{def}}{=}\left\{ M_{h}\in B_{\mathbb{R}}\left(\mathbf{X}_{t}\right):h\in\mathbb{H}_{t}\right\} ,
\]
of $B_{\mathbb{R}}\left(\mathbf{X}_{t}\right)$ forms a complete semi-normed
$*$-algebra over $\mathbb{R}$ of the adjointable operators $\left\{ M_{h}\right\} _{h\in\mathbb{R}}$
in $B_{\mathbb{R}}\left(\mathbf{X}_{t}\right)$ (e.g., see {[}3{]}),
for $t\in\mathbb{R}$. i.e., case-by-case, we understand the set $\mathbb{H}_{t}$
of $t$-scaled hypercomplex numbers as a unital ring $\mathbb{H}_{t}$
with its unity $\left(1,0\right)$ algebraically; or, as a complete
semi-normed $\mathbb{R}$-vector space $\left(\mathbb{H}_{t},\left\Vert .\right\Vert _{t}\right)$,
which is either a $\mathbb{R}$-Hilbert space if $t<0$, or an indefinite
semi-inner product $\mathbb{R}$-space if $t\geq0$, analytically;
or as a complete semi-normed $*$-algebra $\mathcal{M}_{t}$ over
$\mathbb{R}$, operator-algebra-theoretically, for all $t\in\mathbb{R}$.

In particular, in this paper, we regard our family $\mathbb{H}_{t}$
of all $t$-scaled hypercomplex numbers as the $t$-scaled hypercomplex
$\mathbb{R}$-space $\mathbf{X}_{t}$, and act a differential operators
$\nabla_{t}$ on the $\mathbb{R}$-differentiable functions on $\mathbf{X}_{t}$,
where
\[
\nabla_{t}=\frac{\partial}{\partial x_{1}}+i\frac{\partial}{\partial x_{2}}-j_{t}\frac{sgn\left(t\right)\partial}{\sqrt{\left|t\right|}\partial x_{3}}-k_{t}\frac{sgn\left(t\right)\partial}{\sqrt{\left|t\right|}\partial x_{4}},\;\mathrm{if\;}t\neq0,
\]
and\hfill{}(1.1)
\[
\nabla_{0}=\frac{\partial}{\partial x_{1}}+i\frac{\partial}{\partial x_{2}}+j_{0}\frac{\partial}{\partial x_{3}}+k_{0}\frac{\partial}{\partial x_{4}},\;\mathrm{if\;}t=0,
\]

\medskip{}

\noindent for an arbitrarily fixed scale $t\in\mathbb{R}$, where
\[
sgn\left(t\right)\overset{\textrm{def}}{=}\left\{ \begin{array}{ccc}
1 &  & \mathrm{if\;}t>0\\
-1 &  & \mathrm{if\;}t<0,
\end{array}\right.
\]
for all $t\in\mathbb{R}\setminus\left\{ 0\right\} $, and

\medskip{}

\hfill{}$i^{2}=-1,\;\mathrm{and\;}j_{t}^{2}=k_{t}^{2}=t,$\hfill{}(1.2)

\medskip{}

\noindent satisfying the commuting diagrams,
\[
\begin{array}{ccccc}
 &  & i\\
 & ^{1}\swarrow &  & \nwarrow^{-t}\\
 & j_{t} & \underset{1}{\longrightarrow} & k_{t} & ,
\end{array}
\]
and\hfill{}(1.3)
\[
\begin{array}{ccccc}
 &  & i\\
 & ^{t}\nearrow &  & \searrow^{-1}\\
 & j_{t} & \underset{-1}{\longleftarrow} & k_{t} & .
\end{array}
\]
Here, the first diagram of (1.3) means that
\[
ij_{t}=k_{t},\;j_{t}k_{t}=-ti,\;\mathrm{and\;}k_{t}i=j_{t},
\]
while the second diagram of (1.3) means that
\[
ik_{t}=-j_{t},\;k_{t}j_{t}=ti,\;\mathrm{and\;}j_{t}i=-k_{t},
\]
for $t\in\mathbb{R}$.

We study the left, or right $t$-regular functions contained in the
kernel, $ker\nabla_{t}$ of $\nabla_{t}$, and the $t$-harmonic functions
determined by the $t$-regular functions, by defining $t$-scaled
Laplacians,
\[
\varDelta_{t}=\frac{\partial^{2}}{\partial x_{1}^{2}}+\frac{\partial^{2}}{\partial x_{2}^{2}}-\frac{sgn\left(t\right)\partial^{2}}{\partial x_{3}^{2}}-\frac{sgn\left(t\right)\partial^{2}}{\partial x_{4}^{2}},\;\mathrm{if\;}t\neq0,
\]
and
\[
\varDelta_{0}=\frac{\partial^{2}}{\partial x_{1}^{2}}+\frac{\partial^{2}}{\partial x_{2}^{2}}+0\left(\frac{\partial^{2}}{\partial x_{3}^{2}}+\frac{\partial^{2}}{\partial x_{4}^{2}}\right)=\frac{\partial^{2}}{\partial x_{1}^{2}}+\frac{\partial^{2}}{\partial x_{2}^{2}},\;\mathrm{if\;}t=0,
\]
for all $t\in\mathbb{R}$.

\section{Scaled Hypercomplex Numbers}

In this section, we review scaled hypercomplex numbers (for details,
e.g., see {[}1{]}, {[}2{]} and {[}3{]}).

\subsection{Scaled Hypercomplex Rings}

Fix an arbitrarily scale $t\in\mathbb{R}$. Define an operation ($\cdot_{t}$)
on $\mathbb{C}^{2}$ by

\medskip{}

\hfill{}$\left(a_{1},b_{1}\right)\cdot_{t}\left(a_{2},b_{2}\right)\overset{\textrm{def}}{=}\left(a_{1}a_{2}+tb_{1}\overline{b_{2}},\:a_{1}b_{2}+b_{1}\overline{a_{2}}\right),$\hfill{}(2.1.1)

\medskip{}

\noindent for $\left(a_{l},b_{l}\right)\in\mathbb{C}^{2}$, for all
$l=1,2$.
\begin{prop}
The algebraic triple $\left(\mathbb{C}^{2},+,\cdot_{t}\right)$ forms
a unital ring with its unity $\left(1,0\right)$, where ($+$) is
the usual vector addition on $\mathbb{C}^{2}$, and ($\cdot_{t}$)
is the operation $\textrm{(2.1.1)}$.
\end{prop}

\begin{proof}
See {[}1{]} for details.
\end{proof}
One can understand these unital rings $\left\{ \left(\mathbb{C}^{2},+,\cdot_{t}\right)\right\} _{t\in\mathbb{R}}$
as topological rings, since the operations ($+$) and $\left\{ \left(\cdot_{t}\right)\right\} _{t\in\mathbb{R}}$
are continuous on $\mathbb{C}^{2}$.
\begin{defn}
For $t\in\mathbb{R}$, the ring $\mathbb{H}_{t}\overset{\textrm{denote}}{=}\left(\mathbb{C}^{2},+,\cdot_{t}\right)$
is called the $t$-scaled hypercomplex ring.
\end{defn}

For $t\in\mathbb{R}$, define an injective map,
\[
\pi_{t}:\mathbb{H}_{t}\rightarrow M_{2}\left(\mathbb{C}\right),
\]
by\hfill{}(2.1.2)
\[
\pi_{t}\left(\left(a,b\right)\right)=\left(\begin{array}{cc}
a & tb\\
\overline{b} & \overline{a}
\end{array}\right),\;\forall\left(a,b\right)\in\mathbb{H}_{t}.
\]
Such an injective map $\pi_{t}$ of (2.1.2) satisfies that
\[
\pi_{t}\left(h_{1}+h_{2}\right)=\pi_{t}\left(h_{1}\right)+\pi_{t}\left(h_{2}\right),
\]
and\hfill{}(2.1.3)
\[
\pi_{t}\left(h_{1}\cdot_{t}h_{2}\right)=\pi_{t}\left(h_{1}\right)\pi_{t}\left(h_{2}\right),
\]
in $M_{2}\left(\mathbb{C}\right)$, for all $h_{1},h_{2}\in\mathbb{H}_{t}$,
where $\pi_{t}\left(h_{1}\right)\pi_{t}\left(h_{2}\right)$ is the
usual matricial multiplication (e.g., see {[}1{]} for details).
\begin{prop}
The pair $\left(\mathbb{C}^{2},\:\pi_{t}\right)$ forms an injective
Hilbert-space representation of our $t$-scaled hypercomplex ring
$\mathbb{H}_{t}$, where $\pi_{t}$ is an action (2.1.2).
\end{prop}

\begin{proof}
It is shown by (2.1.3), and by the continuity of $\pi_{t}$ (e.g.,
{[}1{]} and {[}2{]}).
\end{proof}
By the injectivity of $\pi_{t}$, one can understand $\mathbb{H}_{t}$
as its realization $\pi_{t}\left(\mathbb{H}_{t}\right)$ as matrices
of $M_{2}\left(\mathbb{C}\right)$.
\begin{defn}
The subset,

\medskip{}

\hfill{}$\pi_{t}\left(\mathbb{H}_{t}\right)=\left\{ \left(\begin{array}{cc}
a & tb\\
\overline{b} & \overline{a}
\end{array}\right)\in M_{2}\left(\mathbb{C}\right):\left(a,b\right)\in\mathbb{H}_{t}\right\} ,$\hfill{}(2.1.4)

\medskip{}

\noindent of $M_{2}\left(\mathbb{C}\right)$, denoted by $\mathcal{H}_{2}^{t}$,
is called the $t$-scaled (hypercomplex-)realization of $\mathbb{H}_{t}$
(in $M_{2}\left(\mathbb{C}\right)$) for $t\in\mathbb{R}$. For convenience,
we denote the realization $\pi_{t}\left(h\right)$ of $h\in\mathbb{H}_{t}$
by $\left[h\right]_{t}$ in $\mathcal{H}_{2}^{t}$.
\end{defn}

If $\mathbb{H}_{t}^{\times}\overset{\textrm{denote}}{=}\mathbb{H}_{t}\setminus\left\{ \left(0,0\right)\right\} ,$
where $\left(0,0\right)$ is the ($+$)-identity of $\mathbb{H}_{t}$,
then, it is the maximal monoid,
\[
\mathbb{H}_{t}^{\times}\overset{\textrm{denote}}{=}\left(\mathbb{H}_{t}^{\times},\;\cdot_{t}\right),
\]
in $\mathbb{H}_{t}$, with its ($\cdot_{t}$)-identity $\left(1,0\right)$,
the unity of $\mathbb{H}_{t}$ (e.g., {[}2{]}).
\begin{defn}
The monoid $\mathbb{H}_{t}^{\times}=\left(\mathbb{H}_{t}^{\times},\:\cdot_{t}\right)$
of $\mathbb{H}_{t}$ is called the $t$-scaled hypercomplex monoid.
\end{defn}

\subsection{Invertible Hypercomplex Numbers of $\mathbb{H}_{t}$}

For an arbitrarily fixed $t\in\mathbb{R}$, let $\mathbb{H}_{t}$
be the corresponding $t$-scaled hypercomplex ring, isomorphic to
its $t$-scaled realization $\mathcal{H}_{2}^{t}$ of (2.1.4). Observe
that, for any $\left(a,b\right)\in\mathbb{H}_{t}$, one has

\medskip{}

\hfill{}$det\left(\left[\left(a,b\right)\right]_{t}\right)=det\left(\begin{array}{cc}
a & tb\\
\overline{b} & \overline{a}
\end{array}\right)=\left|a\right|^{2}-t\left|b\right|^{2}.$\hfill{}(2.2.1)

\medskip{}

\noindent where $det$ is the determinant, and $\left|.\right|$ is
the modulus on $\mathbb{C}$.
\begin{lem}
If $\left(a,b\right)\in\mathbb{H}_{t}$, then $\left|a\right|^{2}\neq t\left|b\right|^{2}$
in $\mathbb{C}$, if and only if $\left(a,b\right)$ is invertible
in $\mathbb{H}_{t}$ with its inverse,

\noindent 
\[
\left(a,b\right)^{-1}=\left(\frac{\overline{a}}{\left|a\right|^{2}-t\left|b\right|^{2}},\:\frac{-b}{\left|a\right|^{2}-t\left|b\right|^{2}}\right)\;\;\mathrm{in\;\;}\mathbb{H}_{t},
\]
satisfying\hfill{}$\textrm{(2.2.2)}$
\[
\left[\left(a,b\right)^{-1}\right]_{t}=\left[\left(a,b\right)\right]_{t}^{-1}\;\;\mathrm{in\;\;\mathcal{H}_{2}^{t}.}
\]
\end{lem}

\begin{proof}
By (2.2.1), we have $det\left(\left[\left(a,b\right)\right]_{t}\right)\neq0$,
if and only if $\left[\left(a,b\right)\right]_{t}$ is invertible
in $M_{2}\left(\mathbb{C}\right)$, with its inverse,
\[
\left[\left(a,b\right)\right]_{t}^{-1}=\left[\left(\frac{\overline{a}}{\left|a\right|^{2}-t\left|b\right|^{2}},\:\frac{-b}{\left|a\right|^{2}-t\left|b\right|^{2}}\right)\right]_{t},
\]
and this inverse is contained in $\mathcal{H}_{2}^{t}$, satisfying
the relation (2.2.2).
\end{proof}
Recall that an algebraic triple $\left(X,+,\cdot\right)$ is a noncommutative
field, if it is a unital ring, containing $\left(X^{\times},\cdot\right)$
as a non-abelian group (e.g., {[}2{]} and {[}3{]}).
\begin{prop}
We have the algebraic characterization,

\medskip{}

\hfill{}$t<0\;\mathrm{in\;}\mathbb{R}\Longleftrightarrow\mathbb{H}_{t}\textrm{ is a noncommutative field.}$\hfill{}$\textrm{(2.2.3)}$
\end{prop}

\begin{proof}
($\Rightarrow$) If $t<0$ in $\mathbb{R}$, then every $\left(a,b\right)\in\mathbb{H}_{t}^{\times}$
satisfies (2.2.2), since
\[
\left|a\right|^{2}>t\left|b\right|^{2},\;\mathrm{implying}\;\left|a\right|^{2}\neq t\left|b\right|^{2}.
\]
Indeed, $\left|a\right|^{2}$ is nonnegative, but $t\left|b\right|^{2}$
is not positive, and $\left(a,b\right)\neq\left(0,0\right)$. Therefore,
if $t<0$, then the monoid $\mathbb{H}_{t}^{\times}$ becomes a noncommutative
group.

\noindent ($\Leftarrow$) Suppose $t\geq0$. First, if $t=0$, and
$\left(0,b\right)\in\mathbb{H}_{0}^{\times}$ ($b\neq0$), then
\[
det\left(\left[\left(0,b\right)\right]_{0}\right)=det\left(\left(\begin{array}{cc}
0 & 0\\
\overline{b} & 0
\end{array}\right)\right)=0,
\]
implying that $\left[\left(0,b\right)\right]_{0}\in\mathcal{H}_{2}^{t}$
is not invertible. If, now, $t>0$, and $\left(a,b\right)\in\mathbb{H}_{t}^{\times}$,
with $\left|b\right|^{2}=\frac{\left|a\right|^{2}}{t}$ in $\mathbb{C}$,
then
\[
det\left(\left[\left(a,b\right)\right]_{t}\right)=\left|a\right|^{2}-t\left|b\right|^{2}=0,
\]
implying that $\left(a,b\right)$ is not invertible in $\mathbb{H}_{t}$.
Thus, if $t\geq0$, then $\mathbb{H}_{t}^{\times}$ cannot be a group.
\end{proof}
By the proof of the above proposition, $\left\{ \mathbb{H}_{s}\right\} _{s<0}$
are noncommutative fields, but, $\left\{ \mathbb{H}_{t}\right\} _{t\geq0}$
cannot be noncommutative fields by (2.2.3). For any scale $t\in\mathbb{R}$,
the $t$-scaled hypercomplex ring $\mathbb{H}_{t}$ is decomposed
by
\[
\mathbb{H}_{t}=\mathbb{H}_{t}^{inv}\sqcup\mathbb{H}_{t}^{sing}
\]
with\hfill{}(2.2.4)
\[
\mathbb{H}_{t}^{inv}=\left\{ \left(a,b\right):\left|a\right|^{2}\neq t\left|b\right|^{2}\right\} ,
\]
and
\[
\mathbb{H}_{t}^{sing}=\left\{ \left(a,b\right):\left|a\right|^{2}=t\left|b\right|^{2}\right\} ,
\]
where $\sqcup$ is the disjoint union.
\begin{prop}
The subset $\mathbb{H}_{t}^{inv}$ of (2.2.4) is a non-abelian group
in the monoid $\mathbb{H}_{t}^{\times}$. Meanwhile, the subset $\mathbb{H}_{t}^{\times sing}\overset{\textrm{denote}}{=}\mathbb{H}_{t}^{sing}\setminus\left\{ \left(0,0\right)\right\} $
forms a semigroup without identity in $\mathbb{H}_{t}^{\times}$.
\end{prop}

\begin{proof}
See {[}2{]} for details.
\end{proof}
\begin{defn}
The block $\mathbb{H}_{t}^{inv}$ of (2.2.4) is called the group-part
of $\mathbb{H}_{t}^{\times}$ (or, of $\mathbb{H}_{t}$), and the
other algebraic block $\mathbb{H}_{t}^{\times sing}$ of the above
proposition is called the semigroup-part of $\mathbb{H}_{t}^{\times}$
(or, of $\mathbb{H}_{t}$).
\end{defn}

By (2.2.3), if $t<0$ in $\mathbb{R}$, then the semigroup-part $\mathbb{H}_{t}^{\times sing}$
is empty in $\mathbb{H}_{t}^{\times}$, and hence, 
\[
\mathbb{H}_{t}^{\times}=\mathbb{H}_{t}^{inv}\Longleftrightarrow\mathbb{H}_{t}=\mathbb{H}_{t}^{inv}\cup\left\{ \left(0,0\right)\right\} ,
\]
Meanwhile, if $t\geq0$ in $\mathbb{R}$, then the semigroup-part
$\mathbb{H}_{t}^{\times sing}$ is non-empty, and is properly contained
in the $t$-scaled monoir $\mathbb{H}_{t}^{\times}$, satisfying (2.2.4).

\subsection{Scaled-Hypercomplex Conjugation}

In this section, we define a suitable conjugation on $\left\{ \mathbb{H}_{t}\right\} _{t\in\mathbb{R}}$,
as in {[}3{]}. For a scale $t\in\mathbb{R}$, define a unary operation
($\dagger$) on the $t$-scaled hypercomplex ring $\mathbb{H}_{t}$
by

\medskip{}

\hfill{}$\left(a,b\right)^{\dagger}\overset{\textrm{def}}{=}\left(\overline{a},\;-b\right),\;\;\forall\left(a,b\right)\in\mathbb{H}_{t}.$\hfill{}(2.3.1)

\medskip{}

\noindent This operation ($\dagger$) of (2.3.1) is indeed a well-defined
unary operation on $\mathbb{H}_{t}$, inducing the equivalent operation,
also denoted by ($\dagger$), on the $t$-scaled realization $\mathcal{H}_{2}^{t}$
of $\mathbb{H}_{t}$,

\medskip{}

\hfill{}$\left[\left(a,b\right)\right]_{t}^{\dagger}\overset{\textrm{def}}{=}\left[\left(a,b\right)^{\dagger}\right]_{t}=\left[\left(\overline{a},-b\right)\right]_{t},$\hfill{}(2.3.2)

\medskip{}

\noindent for all $\left(a,b\right)\in\mathbb{H}_{t}$. Since the
action $\pi_{t}:\mathbb{H}_{t}\rightarrow\mathcal{H}_{2}^{t}$ and
the operation ($\dagger$) of (2.3.1) are bijective, the function
(2.3.2) is also a well-defined bijection on $\mathcal{H}_{2}^{t}$.
\begin{prop}
The bijection ($\dagger$) of (2.3.1) is an adjoint on $\mathbb{H}_{t}$.
\end{prop}

\begin{proof}
It is sufficient to show that the operation ($\dagger$) of (2.3.2)
is an adjoint on $\mathcal{H}_{2}^{t}$. Observe that, for all $\left(a,b\right)\in\mathbb{H}_{t}$,
we have
\[
\left[\left(a,b\right)\right]_{t}^{\dagger\dagger}=\left[\left(a,b\right)^{\dagger}\right]_{t}^{\dagger}=\left[\left(\overline{a},-b\right)\right]_{t}^{\dagger}=\left[\left(a,b\right)\right]_{t},
\]
in $\mathcal{H}_{2}^{t}$. Also, for any $h_{1},h_{2}\in\mathbb{H}_{t}$,
\[
\left[h_{1}+h_{2}\right]_{t}^{\dagger}=\left[\left(h_{1}+h_{2}\right)^{\dagger}\right]_{t}=\left[h_{1}^{\dagger}+h_{2}^{\dagger}\right]_{t}=\left[h_{1}^{\dagger}\right]_{t}+\left[h_{2}^{\dagger}\right]_{t},
\]
implying that
\[
\left[h_{1}+h_{2}\right]_{t}^{\dagger}=\left[h_{1}\right]_{t}^{\dagger}+\left[h_{2}\right]_{t}^{\dagger},\;\mathrm{in\;}\mathcal{H}_{2}^{t};
\]
and, for any $\left(a_{l},b_{l}\right)\in\mathbb{H}_{t}$, for $l=1,2$,

\medskip{}

$\;\;\;$$\left(\left[\left(a_{1},b_{1}\right)\right]\left[\left(a_{2},b_{2}\right)\right]_{t}\right)^{\dagger}=\left(\begin{array}{ccc}
a_{1}a_{2}+tb_{1}\overline{b_{2}} &  & t\left(a_{1}b_{2}+b_{1}\overline{a_{2}}\right)\\
\\
\overline{a_{1}b_{2}+b_{1}\overline{a_{2}}} &  & \overline{a_{1}a_{2}+tb_{1}\overline{b_{2}}}
\end{array}\right)^{\dagger}$

\medskip{}

$\;\;\;\;\;\;\;\;\;\;\;\;\;\;$$=\left(\begin{array}{ccc}
\overline{a_{1}a_{2}+tb_{1}\overline{b_{2}}} &  & t\left(-a_{1}b_{2}-b_{1}\overline{a_{2}}\right)\\
\\
\overline{-a_{1}b_{2}-b_{1}\overline{a_{2}}} &  & a_{1}a_{2}+tb_{1}\overline{b_{2}}
\end{array}\right)$

\medskip{}

$\;\;\;\;\;\;\;\;\;\;\;\;\;\;$$=\left(\begin{array}{ccc}
\overline{a_{2}} &  & t\left(-b_{2}\right)\\
\\
\overline{-b_{2}} &  & a_{2}
\end{array}\right)\left(\begin{array}{ccc}
\overline{a_{1}} &  & t\left(-b_{1}\right)\\
\\
\overline{-b_{1}} &  & a_{1}
\end{array}\right)$

\medskip{}

$\;\;\;\;\;\;\;\;\;\;\;\;\;\;$$=\left[\left(a_{2},b_{2}\right)\right]_{t}^{\dagger}\left[\left(a_{1},b_{1}\right)\right]_{t}^{\dagger}$.

\medskip{}

\noindent Moreover, if $z\in\mathbb{C}$ inducing $\left(z,0\right)\in\mathbb{H}_{t}$
and $\left(a,b\right)\in\mathbb{H}_{t}$, then

\medskip{}

$\;\;\;\;$$\left(\left[\left(z,0\right)\right]_{t}\left[\left(a,b\right)\right]_{t}\right)^{\dagger}=\left(\begin{array}{ccc}
za &  & tzb\\
\\
\overline{zb} &  & \overline{za}
\end{array}\right)^{\dagger}=\left(\begin{array}{ccc}
\overline{za} &  & t\left(-zb\right)\\
\\
\overline{-zb} &  & za
\end{array}\right)$

\medskip{}

$\;\;\;\;\;\;\;\;\;\;\;\;$$=\left[\left(\overline{z},0\right)\right]_{t}\left[\left(\overline{a},-b\right)\right]_{t}=\left[\left(\overline{z},0\right)\right]_{t}\left[\left(a,b\right)\right]_{t}^{\dagger}.$

\medskip{}

\noindent Therefore, the bijection ($\dagger$) of (2.3.2) is an adjoint
on $\mathcal{H}_{2}^{t}$.
\end{proof}
Note that the adjoint ($\dagger$) of (2.3.1) (or, of (2.3.2)) is
free from the choice of scales $t\in\mathbb{R}$. So, we call ($\dagger$)
the hypercomplex-conjugate.
\begin{defn}
We call the adjoint ($\dagger$) on $\mathbb{H}_{t}$, the hypercomplex-conjugate
for all $t\in\mathbb{R}$.
\end{defn}

If $h=\left(a,b\right)\in\mathbb{H}_{t}$, then one obtains that

\medskip{}

\hfill{}$\left[h\right]_{t}^{\dagger}\left[h\right]_{t}=\left[\left(\left|a\right|^{2}-t\left|b\right|^{2},\;0\right)\right]_{t}=\left[h\right]_{t}\left[h\right]_{t}^{\dagger},$\hfill{}(2.3.3)

\medskip{}

\noindent for all $h=\left(a,b\right)\in\mathbb{H}_{t}$, for ``all''
$t\in\mathbb{R}$.

\subsection{The Normalized Trace $\tau$ on $\left\{ \mathbb{H}_{t}\right\} _{t\in\mathbb{R}}$
Over $\mathbb{R}$}

Recall that the $t$-scaled realization $\mathcal{H}_{2}^{t}$ of
the $t$-scaled hypercomplex ring $\mathbb{H}_{t}$ is an embedded
ring of $M_{2}\left(\mathbb{C}\right)$. So, the normalized trace
$\tau=\frac{1}{2}tr$ on $M_{2}\left(\mathbb{C}\right)$ is naturally
restricted to $\tau\mid_{\mathcal{H}_{2}^{t}}$, also denoted by $\tau$,
on $\mathcal{H}_{2}^{t}$, where $tr$ is the usual trace on $M_{2}\left(\mathbb{C}\right)$.
Observe that, for any $\left[\left(a,b\right)\right]_{t}\in\mathcal{H}_{2}^{t}$,
one has
\[
\tau\left(\left[\left(a,b\right)\right]_{t}\right)=\frac{1}{2}tr\left(\left(\begin{array}{cc}
a & tb\\
\overline{b} & \overline{a}
\end{array}\right)\right)=\frac{1}{2}\left(a+\overline{a}\right),
\]
i.e.,\hfill{}(2.4.1)
\[
\tau\left(\left[\left(a,b\right)\right]_{t}\right)=Re\left(a\right),\;\;\forall\left(a,b\right)\in\mathbb{H}_{t}.
\]

\medskip{}

\noindent $\mathbf{Remark\;and\;Discussion.}$ Since $\mathcal{H}_{2}^{t}$
is a sub-structure of $M_{2}\left(\mathbb{C}\right)$, the linear
functional $\tau=\frac{1}{2}tr$ on $M_{2}\left(\mathbb{C}\right)$
is well-restricted to that on $\mathcal{H}_{2}^{t}$. However, note
that, since $\tau\left(\mathcal{H}_{2}^{t}\right)\subseteq\mathbb{R}$
by (2.4.1), it becomes a linear functional ``over $\mathbb{R}$.''
So, in the following text, if we mention ``$\tau$ is a linear functional,''
then it actually means that ``this restriction is linear over $\mathbb{R}$.''
The construction of such a $\mathbb{R}$-linear functional on $\mathcal{H}_{2}^{t}$
is motivated by free probability (e.g., {[}17{]} and {[}20{]}), and
such a free-probabilistic model is considered in detail in {[}1{]},
{[}2{]} and {[}3{]}.\hfill{} \textifsymbol[ifgeo]{64} 

\medskip{}

By (2.4.1), one can define a ($\mathbb{R}$-)linear functional, also
denoted by $\tau$, on $\mathbb{H}_{t}$, by

\medskip{}

\hfill{}$\tau\left(\left(a,b\right)\right)\overset{\textrm{def}}{=}Re\left(a\right),\;\;\;\forall\left(a,b\right)\in\mathbb{H}_{t},$\hfill{}(2.4.2)

\medskip{}

\noindent by (2.4.1). By using this linear functional $\tau$ of (2.4.2)
on $\mathbb{H}_{t}$, we define a form, 
\[
\left\langle ,\right\rangle _{t}:\mathbb{H}_{t}\times\mathbb{H}_{t}\rightarrow\mathbb{R},
\]
by\hfill{}(2.4.3)
\[
\left\langle h_{1},h_{2}\right\rangle _{t}\overset{\textrm{def}}{=}\tau\left(h_{1}\cdot_{t}h_{2}^{\dagger}\right),\;\;\;\forall h_{1},h_{2}\in\mathbb{H}_{t},
\]
where the linear functional $\tau$ in (2.4.3) is in the sense of
(2.4.2), whose range is contained in $\mathbb{R}$ in $\mathbb{C}$.
Since the hypercomplex-conjugation ($\dagger$) is bijective, this
form (2.4.3) is a well-defined function.
\begin{lem}
The function $\left\langle ,\right\rangle _{t}$ of (2.4.3) is a bilinear
form on $\mathbb{H}_{t}$ ``over $\mathbb{R}$.''
\end{lem}

\begin{proof}
By the straightforward computations, one has
\[
\left\langle h_{1}+h_{2},h_{3}\right\rangle _{t}=\left\langle h_{1},h_{2}\right\rangle _{t}+\left\langle h_{2},h_{3}\right\rangle _{t},
\]
\[
\left\langle h_{1},h_{2}+h_{3}\right\rangle _{t}=\left\langle h_{1},h_{2}\right\rangle _{t}+\left\langle h_{1},h_{3}\right\rangle _{t},
\]
and
\[
\left\langle rh_{1},h_{2}\right\rangle _{t}=r\left\langle h_{1},h_{2}\right\rangle _{t},\;\left\langle h_{1},rh_{2}\right\rangle _{t}=r\left\langle h_{1},h_{2}\right\rangle _{t},
\]
for all $h_{1},h_{2},h_{3}\in\mathbb{H}_{t}$, and $r\in\mathbb{R}$,
where
\[
rh=\left(r,0\right)\cdot_{t}h,\;\;\forall h\in\mathbb{H}_{t},\;r\in\mathbb{R}.
\]
See {[}3{]} for details.
\end{proof}
The above lemma shows that the forms $\left\{ \left\langle ,\right\rangle _{t}\right\} _{t\in\mathbb{R}}$
of (2.4.3) are well-determined bilinear forms on $\left\{ \mathbb{H}_{t}\right\} _{t\in\mathbb{R}}$,
induced by the linear functionals $\tau$ of (2.4.2) and the hypercomplex-conjugate
($\dagger$) of (2.3.1).
\begin{lem}
For all $h_{1},h_{2}\in\mathbb{H}_{t}$, we have

\medskip{}

\hfill{}$\left\langle h_{1},h_{2}\right\rangle _{t}=\left\langle h_{2},h_{1}\right\rangle _{t}$
in $\mathbb{R}$.\hfill{}$\textrm{(2.4.4)}$
\end{lem}

\begin{proof}
Observe that, for any $h_{1},h_{2}\in\mathbb{H}_{t}$,

\medskip{}

$\;\;\;\;$$\left\langle h_{1},h_{2}\right\rangle _{t}=\tau\left(h_{1}\cdot_{t}h_{2}^{\dagger}\right)=\tau\left(\left[h_{1}\right]_{t}\left[h_{2}\right]_{t}^{\dagger}\right)$

\medskip{}

\noindent since $\tau$ of (2.4.1) is a trace

\medskip{}

$\;\;\;\;\;\;\;\;$$=\overline{\tau\left(\left[h_{2}\right]_{t}\left[h_{1}\right]_{t}^{\dagger}\right)}=\tau\left(\left[h_{2}\right]_{t}\left[h_{1}\right]_{t}^{\dagger}\right)=\tau\left(h_{2}\cdot_{t}h_{1}^{\dagger}\right)=\left\langle h_{2},h_{1}\right\rangle _{t},$

\medskip{}

\noindent since the normalized trace $\tau$ of (2.4.1) is a linear
functional over $\mathbb{C}$, and because the range of the linear
functional $\tau$ of (2.4.2) is contained in $\mathbb{R}$, and hence,
$\tau$ is linear on $\mathcal{H}_{2}^{t}$ over $\mathbb{R}$. See
{[}3{]} for details.
\end{proof}
The above lemma shows our bilinear form $\left\langle ,\right\rangle _{t}$
of (2.4.3) is symmetric by (2.4.4).
\begin{lem}
If $h_{1},h_{2}\in\mathbb{H}_{t}$, then

\medskip{}

\hfill{}$\left|\left\langle h_{1},h_{2}\right\rangle _{t}\right|^{2}\leq\left|\left\langle h_{1},h_{1}\right\rangle _{t}\right|^{2}\left|\left\langle h_{2},h_{2}\right\rangle _{t}\right|^{2}$,\hfill{}$\textrm{(2.4.5)}$

\medskip{}

\noindent where $\left|.\right|$ is the absolute value on $\mathbb{R}$.
\end{lem}

\begin{proof}
See {[}3{]} for details.
\end{proof}
Observe now that if $h=\left(a,b\right)\in\mathbb{H}_{t}$, then
\[
\left\langle h,h\right\rangle _{t}=\tau\left(\left(a,b\right)\cdot_{t}\left(a,b\right)^{\dagger}\right)=Re\left(\left|a\right|^{2}-t\left|b\right|^{2}\right),
\]
by (2.4.1) and (2.4.2), implying that
\[
\left\langle h,h\right\rangle _{t}=\left|a\right|^{2}-t\left|b\right|^{2}=det\left(\left[h\right]_{t}\right),
\]

\noindent and hence,
\[
\left\langle h,h\right\rangle _{t}=0\Longleftrightarrow det\left(\left[h\right]_{t}\right)=0\;\mathrm{in\;}\mathbb{R},
\]
if and only if\hfill{}(2.4.6)
\[
h=\left(a,b\right)\in\mathbb{H}_{t},\;\mathrm{with\;}\left|a\right|^{2}=t\left|b\right|^{2}.
\]

\begin{prop}
Let $h=\left(a,b\right)\in\mathbb{H}_{t}$. Then $\left\langle h,h\right\rangle _{t}=0$,
if and only if $\left|a\right|^{2}=t\left|b\right|^{2}$, if and only
if $det\left(\left[h\right]_{t}\right)=0$, if and only if $h$ is
not invertible in $\mathbb{H}_{t}$, if and only if $h\in\mathbb{H}_{t}^{sing}$.
i.e.,

\medskip{}

\hfill{}$\left\langle h,h\right\rangle _{t}=0,\Longleftrightarrow h\in\mathbb{H}_{t}^{sing},\;\;\mathrm{in\;\;}\mathbb{H}_{t}.$\hfill{}$\textrm{(2.4.7)}$
\end{prop}

\begin{proof}
The relation (2.4.7) holds by (2.2.4) and (2.4.6).
\end{proof}
If either $h_{1}$ or $h_{2}$ is contained in $\mathbb{H}_{t}^{sing}$,
then $h_{1}\cdot_{t}h_{2}^{\dagger}\in\mathbb{H}_{t}^{sing}$ in $\mathbb{H}_{t}$,
since
\[
det\left(\left[h_{1}\cdot_{t}h_{2}^{\dagger}\right]_{t}\right)=det\left(\left[h_{1}\right]_{t}\left[h_{2}\right]_{t}^{\dagger}\right),
\]
implying that
\[
det\left(\left[h_{1}\cdot_{t}h_{2}^{\dagger}\right]_{t}\right)=\left(det\left(\left[h_{1}\right]_{t}\right)\right)\left(det\left(\left[h_{2}\right]_{t}^{\dagger}\right)\right)=0.
\]

\begin{defn}
For a vector space $X$ over $\mathbb{R}$, a bilinear form $\left\langle ,\right\rangle :X\times X\rightarrow\mathbb{R}$
is called a (definite) semi-inner product on $X$ over $\mathbb{R}$,
if (i)
\[
\left\langle x_{1},x_{2}\right\rangle =\left\langle x_{2},x_{1}\right\rangle ,\;\forall x_{1},x_{2}\in X,
\]
and (ii) $\left\langle x,x\right\rangle \geq0$, for all $x\in X$.
In such a case, the pair $\left(X,\left\langle ,\right\rangle \right)$
is called a semi-inner product space over $\mathbb{R}$ (in short,
a $\mathbb{R}$-SIPS). If a semi-inner product $\left\langle ,\right\rangle $
on a $\mathbb{R}$-vector space $X$ satisfies an additional condition
(iii) $\left\langle x,x\right\rangle =0$, if and only if $x=0_{X}$,
the zero vector of $X$, then it is called an (definite) inner product
on $X$ over $\mathbb{R}$, and the pair $\left(X,\left\langle ,\right\rangle \right)$
is said to be an inner-product space over $\mathbb{R}$ (in short,
a $\mathbb{R}$-IPS).
\end{defn}

\begin{defn}
For a vector space $X$ over $\mathbb{R}$, a bilinear form $\left\langle ,\right\rangle :X\times X\rightarrow\mathbb{R}$
is called an indefinite semi-inner product on $X$ over $\mathbb{R}$,
if (i)
\[
\left\langle x_{1},x_{2}\right\rangle =\left\langle x_{2},x_{1}\right\rangle ,\;\;\forall x_{1},x_{2}\in X,
\]
and (ii) $\left\langle x,x\right\rangle \in\mathbb{R}$, for all $x\in X$.
In such a case, the pair $\left(X,\left\langle ,\right\rangle \right)$
is called an indefinite-semi-inner product space over $\mathbb{R}$
(in short, a $\mathbb{R}$-ISIPS).
\end{defn}

If an indefinite semi-inner product $\left\langle ,\right\rangle $
on a $\mathbb{R}$-vector space $X$ satisfies an additional condition;
$\left\langle x,x\right\rangle =0$, if and only if $x=0_{X}$, then
it is called an indefinite inner product on $X$ over $\mathbb{R}$,
and the pair $\left(X,\left\langle ,\right\rangle \right)$ is said
to be an indefinite-inner product space over $\mathbb{R}$ (in short,
a $\mathbb{R}$-IIPS).

By the above lemmas and proposition, one obtains the following result.
\begin{thm}
If $t<0$ in $\mathbb{R}$, then the bilinear form $\left\langle ,\right\rangle _{t}$
of (2.4.3) is a continuous inner product on $\mathbb{H}_{t}$ over
$\mathbb{R}$, i.e.,

\medskip{}

\hfill{}$t<0\Longrightarrow\left(\mathbb{H}_{t},\left\langle ,\right\rangle _{t}\right)$
is a $\mathbb{R}$-IPS.\hfill{}$\textrm{(2.4.8)}$

\medskip{}

\noindent Meanwhile, if $t\geq0$, then $\left\langle ,\right\rangle _{t}$
is a continuous indefinite semi-inner product on $\mathbb{H}_{t}$
over $\mathbb{R}$, i.e.,

\medskip{}

\hfill{}$t\geq0\Longrightarrow\left(\mathbb{H}_{t},\left\langle ,\right\rangle _{t}\right)$
is a $\mathbb{R}$-ISIPS, but not a $\mathbb{R}$-IIPS.\hfill{}$\textrm{(2.4.9)}$
\end{thm}

\begin{proof}
Assume first that a given scale $t$ is negative in $\mathbb{R}$.
Then the bilinear form $\left\langle ,\right\rangle _{t}$ forms a
semi-inner product by the symmetry (2.4.4), and
\[
\left\langle \left(a,b\right),\left(a,b\right)\right\rangle _{t}=\left|a\right|^{2}-t\left|b\right|^{2}\geq0,
\]
since $t<0$, for all $\left(a,b\right)\in\mathbb{H}_{t}$. Also,
recall that, if $t<0$, then $\mathbb{H}_{t}=\mathbb{H}_{t}^{inv}\cup\left\{ \left(0,0\right)\right\} $,
equivalently, $\mathbb{H}_{t}^{sing}=\left\{ \left(0,0\right)\right\} $.
It implies that
\[
\left\langle h,h\right\rangle _{t}=0\Longleftrightarrow h\in\mathbb{H}_{t}^{sing}\Longleftrightarrow h=\left(0,0\right)\in\mathbb{H}_{t}.
\]
i.e., $\left\langle ,\right\rangle _{t}$ is an inner product on $\mathbb{H}_{t}$
over $\mathbb{R}$, whenever $t<0$ in $\mathbb{R}$. The continuity
of $\left\langle ,\right\rangle _{t}$ is guaranteed by (2.4.5). So,
if $t<0$, then the pair $\left(\mathbb{H}_{t},\left\langle ,\right\rangle _{t}\right)$
is not only a $\mathbb{R}$-SIPS, but also a $\mathbb{R}$-IPS. Therefore,
the relation (2.4.8) holds.

Assume now that a scale $t$ is nonnegative in $\mathbb{R}$. Then
the bilinear form $\left\langle ,\right\rangle _{t}$ forms an indefinite
semi-inner product by the symmetry (2.4.4), and
\[
\left\langle h,h\right\rangle _{t}=det\left(\left[h\right]_{t}\right)\in\mathbb{R},
\]
since $t\geq0$, for all $h\in\mathbb{H}_{t}$. Note that if $t\geq0$,
then the semigroup-part $\mathbb{H}_{t}^{\times sing}$ is not empty
in $\mathbb{H}_{t}$, satisfying
\[
\left\langle h,h\right\rangle _{t}=0\Longleftrightarrow h\in\mathbb{H}_{t}^{sing},
\]
by (2.4.7). So, there are lots of ``non-zero'' hypercomplex numbers
$h$ in $\mathbb{H}_{t}^{\times}$ satisfying $\left\langle h,h\right\rangle _{t}=0$,
whenever $t\geq0$ in $\mathbb{R}$. So, if $t\geq0$, then $\left\langle ,\right\rangle _{t}$
forms an indefinite semi-inner product on $\mathbb{H}_{t}$ over $\mathbb{R}$,
which cannot be an indefinite inner product on $\mathbb{H}_{t}$.
The continuity of $\left\langle ,\right\rangle _{t}$ is shown by
(2.4.5). Therefore, if $t\geq0$, then the pair $\left(\mathbb{H}_{t},\left\langle ,\right\rangle _{t}\right)$
forms a $\mathbb{R}$-ISIPS, which is not a $\mathbb{R}$-IIPS, implying
the statement (2.4.9).
\end{proof}
The above theorem shows that $\left\{ \mathbb{H}_{t}\right\} _{t<0}$
are $\mathbb{R}$-IPSs, meanwhile, $\left\{ \mathbb{H}_{t}\right\} _{t\geq0}$
form $\mathbb{R}$-ISIPS, by (2.4.8) and (2.4.9), respectively.
\begin{defn}
Let $X$ be a vector space over $\mathbb{R}$, and $\left\Vert .\right\Vert _{X}:X\rightarrow\mathbb{R}$,
a function satisfying (i) $\left\Vert x\right\Vert _{X}\geq0$, for
all $x\in X$, and (ii) $\left\Vert rx\right\Vert _{X}=\left|r\right|\left\Vert x\right\Vert _{X}$,
for all $r\in\mathbb{R}$ and $x\in X$, and (iii) 
\[
\left\Vert x_{1}+x_{2}\right\Vert _{X}\leq\left\Vert x_{1}\right\Vert _{X}+\left\Vert x_{2}\right\Vert _{X},\;\forall x_{1},x_{2}\in X.
\]
Then the function $\left\Vert .\right\Vert _{X}$ is called a semi-norm
on $X$ over $\mathbb{R}$, and the pair $\left(X,\left\Vert .\right\Vert _{X}\right)$
is said to be a semi-normed space over $\mathbb{R}$ (in short, a
$\mathbb{R}$-SNS). If a semi-norm $\left\Vert .\right\Vert _{X}$
satisfies an additional condition (iv) $\left\Vert x\right\Vert _{X}=0$,
if and only if $x=0_{X}$ in $X$, then it is called a norm on $X$
over $\mathbb{R}$, and the pair $\left(X,\left\Vert .\right\Vert _{X}\right)$
is said to be a normed space over $\mathbb{R}$ (in short, a $\mathbb{R}$-NS).
\end{defn}

Let $\left(\mathbb{H}_{t},\left\langle ,\right\rangle _{t}\right)$
be either a $\mathbb{R}$-IPS (if $t<0$), or a $\mathbb{R}$-ISIPS
(if $t\geq0$), for an arbitrary scale $t\in\mathbb{R}$. Define a
function,
\[
\left\Vert .\right\Vert _{t}:\mathbb{H}_{t}\rightarrow\mathbb{R},
\]
by\hfill{}(2.4.10)
\[
\left\Vert h\right\Vert _{t}\overset{\textrm{def}}{=}\sqrt{\left|\left\langle h,h\right\rangle _{t}\right|}=\sqrt{\left|det\left(\left[h\right]_{t}\right)\right|}=\sqrt{\left|\tau\left(\left[h\cdot_{t}h^{\dagger}\right]\right)\right|},
\]
for all $h\in\mathbb{H}_{t}$, where $\left|.\right|$ is the absolute
value on $\mathbb{R}$. Then it is not hard to check that
\[
\left\Vert h\right\Vert _{t}\geq0,\;\mathrm{since\;}\left|\left\langle h,h\right\rangle _{t}\right|=\left|\left|a\right|^{2}-t\left|b\right|^{2}\right|\geq0,
\]
and
\[
\left\Vert rh\right\Vert _{t}=\sqrt{\left|\left\langle rh,rh\right\rangle _{t}\right|}=\sqrt{r^{2}}\sqrt{\left|\left\langle h,h\right\rangle _{t}\right|}=\left|r\right|\left\Vert h\right\Vert _{t},
\]
for all $r\in\mathbb{R}$ and $h\in\mathbb{H}_{t}$, and
\[
\left\Vert h_{1}+h_{2}\right\Vert _{t}^{2}=\left|\left\langle h_{1}+h_{2},h_{1}+h_{2}\right\rangle _{t}\right|\leq\left(\left\Vert h_{1}\right\Vert _{t}+\left\Vert h_{2}\right\Vert _{t}\right)^{2},
\]
by (2.4.5) and the bilinearity of $\left\langle ,\right\rangle _{t}$,
implying that
\[
\left\Vert h_{1}+h_{2}\right\Vert _{t}\leq\left\Vert h_{1}\right\Vert _{t}+\left\Vert h_{2}\right\Vert _{t},\;\mathrm{in\;}\mathbb{R},
\]
for all $h_{1},h_{2}\in\mathbb{H}_{t}$ (See {[}3{]} for details).
\begin{thm}
The pair $\left(\mathbb{H}_{t},\left\Vert .\right\Vert _{t}\right)$
of the $t$-scaled hypercomplex ring $\mathbb{H}_{t}$ and the function
$\left\Vert .\right\Vert _{t}$ of (2.4.10) forms a complete $\mathbb{R}$-SNS.
In particular,

\medskip{}

\hfill{}$t<0$$\Longrightarrow$$\left(\mathbb{H}_{t},\left\Vert .\right\Vert _{t}\right)$
is a complete $\mathbb{R}$-NS,\hfill{}$\textrm{(2.4.11)}$

\noindent and

\hfill{}$t\geq0\Longrightarrow\left(\mathbb{H}_{t},\left\Vert .\right\Vert _{t}\right)$
is a complete $\mathbb{R}$-SNS.\hfill{}$\textrm{(2.4.12)}$
\end{thm}

\begin{proof}
By the very above paragraph, the map $\left\Vert .\right\Vert _{t}$
of (2.4.10) is a semi-norm on $\mathbb{H}_{t}$ over $\mathbb{R}$.
In particular, the $\mathbb{R}$-SNS $\left(\mathbb{H}_{t},\left\Vert .\right\Vert _{t}\right)$
is complete by (2.4.5) and (2.4.10), since $\mathbb{H}_{t}$ is a
subspace of the finite-dimensional $\mathbb{R}$-vector space $\mathbb{R}^{4}$.
More precisely, if $t<0$ in $\mathbb{R}$, then this semi-norm $\left\Vert .\right\Vert _{t}$
satisfies the additional condition,
\[
\left\Vert h\right\Vert _{t}=0\Longleftrightarrow\left\langle h,h\right\rangle _{t}=0\Longleftrightarrow h\in\mathbb{H}_{t}^{sing}\Longleftrightarrow h=\left(0,0\right),
\]
implying that $\left\Vert .\right\Vert _{t}$ forms a norm on $\mathbb{H}_{t}$
over $\mathbb{R}$. So, the statement (2.4.11) holds true. Meanwhile
if $t\geq0$, then this semi-norm $\left\Vert .\right\Vert _{t}$
cannot be a norm, because,
\[
\left\Vert h\right\Vert _{t}=0\Longleftrightarrow h\in\mathbb{H}_{t}^{sing},
\]
by (2.4.7). Therefore, the statement (2.4.12) holds.
\end{proof}
Recall that complete IPSs (over $\mathbb{R}$, or over $\mathbb{C}$)
are said to be Hilbert spaces (over $\mathbb{R}$, respectively, over
$\mathbb{C}$). So, the statement (2.4.11) can be re-stated by that:
if $t<0$, then $\left(\mathbb{H}_{t},\left\langle ,\right\rangle _{t}\right)$
is a Hilbert space over $\mathbb{R}$ (or, a $\mathbb{R}$-Hilbert
space).
\begin{defn}
Let $\left(\mathbb{H}_{t},\left\langle ,\right\rangle _{t}\right)$
be either a $\mathbb{R}$-Hilbert space (if $t<0$), or a $\mathbb{R}$-ISIPS
(if $t\geq0$), for an arbitrary scale $t\in\mathbb{R}$, which is
a complete $\mathbb{R}$-SNS $\left(\mathbb{H}_{t},\left\Vert .\right\Vert _{t}\right)$,
where $\left\Vert .\right\Vert _{t}$ is the semi-norm (2.4.10). Then
we call $\left(\mathbb{H}_{t},\left\langle ,\right\rangle _{t}\right)$,
the $t$(-scaled)-hypercomplex $\mathbb{R}$(-vector)-space. And,
to distinguish the $t$-scaled hypercomplex ring $\mathbb{H}_{t}$
and the $t$-hypercomplex $\mathbb{R}$-space $\left(\mathbb{H}_{t},\left\langle ,\right\rangle _{t}\right)$,
we denote $\left(\mathbb{H}_{t},\left\langle ,\right\rangle _{t}\right)$
by $\mathbf{X}_{t}$, for all $t\in\mathbb{R}$.
\end{defn}

The hypercomplex $\mathbb{R}$-spaces $\left\{ \mathbf{X}_{t}\right\} _{t<0}$
are $\mathbb{R}$-Hilbert spaces by (2.4.11), while $\left\{ \mathbf{X}_{t}\right\} _{t\geq0}$
are $\mathbb{R}$-ISIPSs by (2.4.12).

\section{Scaled Hypercomplex $\mathbb{R}$-Spaces $\left\{ \mathbf{X}_{t}\right\} _{t\in\mathbb{R}}$}

In the rest of this paper, we understand the $t$-scaled hypercomplex
ring $\mathbb{H}_{t}$ as the $t$-hypercomplex $\mathbb{R}$-space
$\mathbf{X}_{t}=\left(\mathbb{H}_{t},\left\langle ,\right\rangle _{t}\right)$,
which is either a $\mathbb{R}$-Hilbert space if $t<0$, or a complete
$\mathbb{R}$-ISIPS (where the completeness is up to that of the corresponding
complete $\mathbb{R}$-SNS, $\left(\mathbb{H}_{t},\left\Vert .\right\Vert _{t}\right)$
induced by the semi-norm $\left\Vert .\right\Vert _{t}$ of (2.4.10))
if $t\geq0$. In this section, we consider a suitable setting of $\mathbf{X}_{t}$
to study the differentiation of the functions acting on $\mathbf{X}_{t}$.
More precisely, to study how an operator,
\[
D_{t}=\frac{\partial}{\partial x_{1}}+i\frac{\partial}{\partial x_{2}}+j_{t}\frac{\partial}{\partial x_{3}}+k_{t}\frac{\partial}{\partial x_{4}}
\]
under the conditions (1.2) and (1.3), acts on the set of $\mathbb{R}$-differentiable
functions on $\mathbf{X}_{t}$, we investigate an equivalent setting
for $\mathbf{X}_{t}$, which is different notationally from that of
Section 2.

Fix an arbitrary scale $t\in\mathbb{R}$, and the $t$-hypercomplex
$\mathbb{R}$-space $\mathbf{X}_{t}$. Let
\[
a=x+yi,\;\mathrm{and\;}b=u+vi
\]
be complex numbers with $i=\sqrt{-1}$ and $x,y,u,v\in\mathbb{R}$,
and let $\left(a,b\right)\in\mathbf{X}_{t}$, whose realization is

\medskip{}

\hfill{}$\left[\left(a,b\right)\right]_{t}=\left(\begin{array}{cc}
a & tb\\
\overline{b} & \overline{a}
\end{array}\right)=\left(\begin{array}{ccc}
x+yi &  & tu+tvi\\
\\
u-vi &  & x-yi
\end{array}\right),$\hfill{}(3.1)

\medskip{}

\noindent in the $t$-scaled realization $\mathcal{H}_{2}^{t}$ of
$\mathbb{H}_{t}\overset{\textrm{set}}{=}\mathbf{X}_{t}$. Then it
is identified with

\medskip{}

\hfill{}$x\left(\begin{array}{cc}
1 & 0\\
0 & 1
\end{array}\right)+y\left(\begin{array}{cc}
i & 0\\
0 & -i
\end{array}\right)+u\left(\begin{array}{cc}
0 & t\\
1 & 0
\end{array}\right)+v\left(\begin{array}{cc}
0 & ti\\
-i & 0
\end{array}\right),$\hfill{}(3.2)

\medskip{}

\noindent by (3.1), where
\[
\left(\begin{array}{cc}
1 & 0\\
0 & 1
\end{array}\right)=\left[\left(1,0\right)\right]_{t},\;\;\left(\begin{array}{cc}
i & 0\\
0 & -i
\end{array}\right)=\left[\left(i,0\right)\right]_{t},
\]
and
\[
\left(\begin{array}{cc}
0 & t\\
1 & 0
\end{array}\right)=\left[\left(0,1\right)\right]_{t},\;\;\left(\begin{array}{cc}
0 & ti\\
-i & 0
\end{array}\right)=\left[\left(0,i\right)\right]_{t},
\]
in $\mathcal{H}_{2}^{t}$. It means that every realization $T$ of
$\mathcal{H}_{2}^{t}$ is expressed by
\[
T=x\left[\left(1,0\right)\right]_{t}+y\left[\left(i,0\right)\right]_{t}+u\left[\left(0,1\right)\right]_{t}+v\left[\left(0,i\right)\right]_{t},
\]
for some $x,y,u,v\in\mathbb{R}$. i.e., the $t$-scaled realization
is spanned by
\[
\left\{ \left[\left(1,0\right)\right]_{t},\;\left[\left(i,0\right)\right]_{t},\:\left[\left(0,1\right)\right],\:\left[\left(0,i\right)\right]_{t}\right\} ,
\]
over $\mathbb{R}$.
\begin{lem}
For $t\in\mathbb{R}$, the $t$-hypercomplex $\mathbb{R}$-space $\mathbf{X}_{t}$
is spanned by the subset,

\medskip{}

\hfill{}$\mathbf{B}_{t}=\left\{ \left(1,0\right),\:\left(i,0\right),\:\left(0,1\right),\:\left(0,i\right)\right\} ,$\hfill{}$\textrm{(3.3)}$

\medskip{}

\noindent over $\mathbb{R}$, in the sense that: for every $h\in\mathbf{X}_{t}$,
there exist $x,y,u,v\in\mathbb{R}$, such that
\[
h=x\left(1,0\right)+y\left(i,0\right)+u\left(0,1\right)+v\left(0,i\right),
\]
where $r\left(a,b\right)=\left(r,0\right)\cdot_{t}\left(a,b\right)=\left(ra,rb\right)$,
for all $r\in\mathbb{R}$ and $\left(a,b\right)\in\mathbf{X}_{t}$.
i.e.,
\[
\mathbf{X}_{t}=span_{\mathbb{R}}\left(\mathbf{B}_{t}\right).
\]
\end{lem}

\begin{proof}
By the injective Hilbert-space representation $\left(\mathbb{C}^{2},\pi_{t}\right)$,
the $t$-scaled hypercomplex ring $\mathbb{H}_{t}$ is isomorphic
to the $t$-scaled realization $\mathcal{H}_{2}^{t}$. And since the
$t$-hypercomplex $\mathbb{R}$-space $\mathbf{X}_{t}$ is equipotent
to $\mathbb{H}_{t}$, $\mathcal{H}_{2}^{t}$ and $\mathbf{X}_{t}$
are isomorphic as $\mathbb{R}$-vector spaces. Therefore, the subset
$\mathbf{B}_{t}$ of (3.3) spans $\mathbf{X}_{t}$ over $\mathbb{R}$,
by (3.2).
\end{proof}
The spanning subset $\mathbf{B}_{t}$ of (3.3) satisfies the following
properties in $\mathbf{X}_{t}$.
\begin{prop}
By $\mathit{i}$, $j(t)$ and $k(t)$, we denote the spanning vectors
$\left(i,0\right)$, $\left(0,1\right)$ and $\left(0,i\right)$ of
$\mathbf{B}_{t}$, where $\mathbf{X}_{t}=span_{\mathbb{R}}\mathbf{B}_{t}$.
Then

\medskip{}

\hfill{}$\left\{ \begin{array}{c}
i^{2}=\left(-1,0\right),\;j(t)^{2}=\left(t,0\right)=k(t)^{2},\\
\\
ij(t)=k(t),\;j(t)k(t)=-ti,\;\mathrm{and\;}k(t)i=j(t),\\
\\
k(t)j(t)=ti,\;j(t)i=-k(t),\;\mathrm{and\;}ik(t)=-j(t).
\end{array}\right\} $\hfill{}$\textrm{(3.4)}$
\end{prop}

\begin{proof}
By the above lemma, equivalently, one has that
\[
\mathcal{H}_{2}^{t}=span_{\mathbb{R}}\left(\left\{ \mathbf{1},\:\mathbf{i},\:\mathbf{j}_{t},\:\mathbf{k}_{t}\right\} \right),
\]
where
\[
\mathbf{1}=\left[\left(1,0\right)\right]_{t},\;\;\mathbf{i}=\left[\left(i,0\right)\right]_{t},
\]
and
\[
\mathbf{j}_{t}=\left[\left(0,1\right)\right]_{t},\;\mathrm{and\;}\mathbf{k}_{t}=\left[\left(0,i\right)\right]_{t}.
\]
So, to prove the formulas of (3.4), it suffices to show that
\[
\left\{ \begin{array}{c}
\mathbf{i}^{2}=\left[\left(-1,0\right)\right]_{t},\;\mathbf{j}_{t}^{2}=\left[\left(t,0\right)\right]_{t}=\mathbf{k}_{t}^{2},\\
\\
\mathbf{i}\mathbf{j}_{t}=\mathbf{k}_{t},\;\mathbf{j}_{t}\mathbf{k}_{t}=-t\mathbf{i},\;\mathrm{and\;}\mathbf{k}_{t}\mathbf{i}=\mathbf{j}_{t},\\
\\
\mathbf{k}_{t}\mathbf{j}_{t}=t\mathbf{i},\;\mathbf{j}_{t}\mathbf{i}=-\mathbf{k}_{t},\;\mathrm{and\;}\mathbf{i}\mathbf{k}_{t}=-\mathbf{j}_{t},
\end{array}\right\} 
\]
in $\mathcal{H}_{2}^{t}$, respectively. And they are proven by the
straightforward computations. For instance,
\[
\mathbf{k}_{t}^{2}=\left(\begin{array}{cc}
0 & ti\\
-i & 0
\end{array}\right)\left(\begin{array}{cc}
0 & ti\\
-i & 0
\end{array}\right)=\left(\begin{array}{cc}
t & 0\\
0 & t
\end{array}\right)=\left[\left(t,0\right)\right]_{t}=t\mathbf{1};
\]
\[
\mathbf{i}\mathbf{j}_{t}=\left(\begin{array}{cc}
i & 0\\
0 & -i
\end{array}\right)\left(\begin{array}{cc}
0 & t\\
1 & 0
\end{array}\right)=\left(\begin{array}{cc}
0 & ti\\
-i & 0
\end{array}\right)=\left[\left(0,i\right)\right]_{t}=\mathbf{k}_{t},
\]
\[
\mathbf{i}\mathbf{k}_{t}=\left(\begin{array}{cc}
i & 0\\
0 & -i
\end{array}\right)\left(\begin{array}{cc}
0 & ti\\
-i & 0
\end{array}\right)=\left(\begin{array}{cc}
0 & -t\\
-1 & 0
\end{array}\right)=-\left[\left(0,1\right)\right]_{t}=-\mathbf{j}_{t},
\]
\[
\mathbf{k}_{t}\mathbf{j}_{t}=\left(\begin{array}{cc}
0 & ti\\
-i & 0
\end{array}\right)\left(\begin{array}{cc}
0 & t\\
1 & 0
\end{array}\right)=\left(\begin{array}{cc}
ti & 0\\
0 & -ti
\end{array}\right)=t\left[\left(i,0\right)\right]_{t}=t\mathbf{i},
\]
and
\[
\mathbf{j}_{t}\mathbf{k}_{t}=\left(\begin{array}{cc}
0 & t\\
1 & 0
\end{array}\right)\left(\begin{array}{cc}
0 & ti\\
-i & 0
\end{array}\right)=\left(\begin{array}{cc}
-ti & 0\\
0 & ti
\end{array}\right)=-t\left[\left(i,0\right)\right]_{t}=-t\mathbf{i}.
\]
etc.. Therefore, the formulas of (3.4) are shown.
\end{proof}
By the above proposition, one can conclude that

\medskip{}

\hfill{}$\left(i,0\right)^{2}=\left(-1,0\right),\;\mathrm{and\;}j(t)^{2}=\left(t,0\right)=k(t)^{2},$\hfill{}(3.5)

\medskip{}

\noindent and the following commuting diagrams hold;
\[
\begin{array}{ccccc}
 &  & \left(i,0\right)\\
 & ^{1}\swarrow &  & \nwarrow^{-t}\\
 & j(t) & \underset{1}{\longrightarrow} & k(t) & ,
\end{array}
\]
and\hfill{}(3.6)
\[
\begin{array}{ccccc}
 &  & \left(i,0\right)\\
 & ^{t}\nearrow &  & \searrow^{-1}\\
 & j(t) & \underset{-1}{\longleftarrow} & k(t) & ,
\end{array}
\]
where, the first diagram of (3.6) illustrates that
\[
\left(i,0\right)j(t)=k(t),\;j(t)k(t)=-t\left(i,0\right),\;\mathrm{and\;}k(t)\left(i,0\right)=j(t),
\]
meanwhile the second diagram of (3.6) illustrates that
\[
\left(i,0\right)k(t)=-j(t),\;k(t)j(t)=t\left(i,0\right),\;\mathrm{and\;}j(t)\left(i,0\right)=-k(t),
\]
for $t\in\mathbb{R}$. If we define a vector space $\mathscr{H}_{t}$
over $\mathbb{R}$, by
\[
\mathscr{H}_{t}=\left\{ x+yi+uj_{t}+vk_{t}:x,y,u,v\in\mathbb{R}\right\} ,
\]
i.e.,
\[
\mathscr{H}_{t}=span_{\mathbb{R}}\left\{ 1,i,j_{t},k_{t}\right\} ,
\]
where $i$, $j_{t}$ and $k_{t}$ are certain imaginary numbers (as
spanning vectors) satisfying (1.2) and (1.3), then the morphism,
\[
x\left(1,0\right)+y\left(i,0\right)+uj(t)+vk(t)\in\mathbf{X}_{t}\mapsto x+yi+uj_{t}+vk_{t}\in\mathscr{H}_{t},
\]
is a well-defined bijection preserving the formulas of (3.4) to the
relations of (1.2) and (1.3). Since $\mathbf{X}_{t}=span_{\mathbb{R}}\mathbf{B}_{t}$,
and 
\[
\mathscr{H}_{t}=span_{\mathbb{R}}\left(\left\{ 1,i,j_{t},k_{t}\right\} \right),
\]
this bijection becomes a $\mathbb{R}$-vector-space-isomorphism. More
precisely, we have the following structure theorem of the set $\mathscr{H}_{t}$
of (3.7).
\begin{thm}
Let $\mathscr{H}_{t}$ be the vector space (3.7) spanned by $\left\{ 1,i,j_{t},k_{t}\right\} $,
satisfying the conditions (1.2) and (1.3), over $\mathbb{R}$. Then

\medskip{}

\hfill{}$\mathscr{H}_{t}$ is a unital ring with its unity $1=1+0i+0j_{t}+0k_{t}$,\hfill{}$\textrm{(3.8)}$

\noindent algebraically, and

\hfill{}$\mathscr{H}_{t}$ is a complete $\mathbb{R}$-SNS, isomorphic
to $\mathbf{X}_{t}$,\hfill{}$\textrm{(3.9)}$

\medskip{}

\noindent analytically, especially, it is a $\mathbb{R}$-Hilbert
space if $t<0$, while it is a complete $\mathbb{R}$-ISIPS if $t\geq0$,
and

\medskip{}

\hfill{}$\mathscr{H}_{t}$ is a complete semi-normed $*$-algebra
over $\mathbb{R}$,\hfill{}$\textrm{(3.10)}$

\medskip{}

\noindent operator-algebraically, especially, it becomes a $C^{*}$-algebra
over $\mathbb{R}$ if $t<0$, while it is a complete semi-normed $*$-algebra
over $\mathbb{R}$ if $t\geq0$.
\end{thm}

\begin{proof}
By the bijection,
\[
x\left(1,0\right)+y\left(i,0\right)+uj(t)+vk(t)\in\mathbf{X}_{t}\longmapsto x+yi+uj_{t}+vk_{t}\in\mathscr{H}_{t},
\]
two mathematical structures $\mathbf{X}_{t}$ and $\mathscr{H}_{t}$
are isomorphic algebraically, and analytically, because
\[
\mathbf{X}_{t}=span_{\mathbb{R}}\left\{ \left(1,0\right),\left(i,0\right),j_{t},k_{t}\right\} ,\;\mathrm{and\;}\mathscr{H}_{t}=span_{\mathbb{R}}\left\{ 1,i,j_{t},k_{t}\right\} ,
\]
over $\mathbb{R}$. So, the statements (3.8) and (3.9) hold, since
$\mathbb{H}_{t}\overset{\textrm{set}}{=}\mathbf{X}_{t}$ is a unital
ring algebraically, and a complete $\mathbb{R}$-SNS analytically.

In {[}3{]}, we showed that, on the $t$-hypercomplex $\mathbb{R}$-space
$\mathbf{X}_{t}$, one can define the operator space $B_{\mathbb{R}}\left(\mathbf{X}_{t}\right)$
of all bounded linear transformation (simply, operators) acting on
this $\mathbb{R}$-vector space $\mathbf{X}_{t}$ ``over $\mathbb{R}$,''
equipped with the operator semi-norm,
\[
\left\Vert T\right\Vert \overset{\textrm{def}}{=}sup\left\{ \left\Vert T\left(h\right)\right\Vert _{t}:\left\Vert h\right\Vert _{t}=1\right\} ,
\]
for all $T\in B_{\mathbb{R}}\left(\mathbf{X}_{t}\right)$. Then one
can canonically define multiplication operators $M_{h}\in B_{\mathbb{R}}\left(\mathbf{X}_{t}\right)$
with their symbols $h\in\mathbf{X}_{t}$, by
\[
M_{h}\left(\eta\right)\overset{\textrm{def}}{=}h\cdot_{t}\eta\in\mathbf{X}_{t}\;\;\;\forall\eta\in\mathbf{X}_{t}.
\]
And the subset,
\[
\mathcal{M}_{t}=\left\{ M_{h}\in B_{\mathbb{R}}\left(\mathbf{X}_{t}\right):h\in\mathbf{X}_{t}\right\} ,
\]
of $B_{\mathbb{R}}\left(\mathbf{X}_{t}\right)$ forms a well-defined
complete semi-normed $*$-algebra over $\mathbb{R}$. In particular,
the operator semi-norm $\left\Vert .\right\Vert $ becomes an operator
``norm,'' if $t<0$ (because if $t<0$, then $\left\Vert .\right\Vert _{t}$
is a norm induced by the inner product $\left\langle ,\right\rangle _{t}$),
meanwhile, it is a operator ``semi-norm,'' if $t\geq0$ (because
if $t\geq0$, then $\left\Vert .\right\Vert _{t}$ forms a semi-norm
induced by the indefinite semi-inner product $\left\langle ,\right\rangle _{t}$)
. See {[}3{]} for details. Since $\mathbf{X}_{t}\overset{\textrm{set}}{=}\mathbb{H}_{t}$
is isomorphic to $\mathscr{H}_{t}$, this complete semi-normed $*$-algebra
$\mathcal{M}_{t}$ is isomorphic to $\mathscr{H}_{t}$, operator-algebraically,
too. Therefore, the statement (3.10) holds, too.
\end{proof}
The above theorem characterizes the algebraic, analytic, and operator-algebraic
structures of the set $\mathscr{H}_{t}$ of (3.7), for $t\in\mathbb{R}$.
So, one can re-define the set $\mathbb{H}_{t}$ of $t$-scaled hypercomplex
numbers by the equivalent set $\mathscr{H}_{t}$.
\begin{defn}
Let $t\in\mathbb{R}$ be fixed. Then the set, also denoted by $\mathbb{H}_{t}$,
\[
\mathbb{H}_{t}\overset{\textrm{def}}{=}\left\{ x+yi+uj_{t}+vk_{t}\left|\begin{array}{c}
x,\;y,\;u,\;v\in\mathbb{R},\\
\mathrm{and}\\
i^{2}=-1,\;\mathit{j_{t}^{2}=t=k_{t}^{2}},\\
ij_{t}=k_{t},\:j_{t}k_{t}=-ti,\:k_{t}i=j_{t},\\
j_{t}i=-k_{t},\:ik_{t}=-j_{t},\:k_{t}j_{t}=ti
\end{array}\right.\right\} ,
\]
is called the $t$-scaled hypercomplexes. All elements of this set
$\mathbb{H}_{t}$ are said to be $t$-scaled hypercomplex numbers.
\end{defn}

By the above theorem, the $t$-scaled hypercomplexes $\mathbb{H}_{t}$
is a unital ring with its unity $1$ algebraically, and a complete
$\mathbb{R}$-SNS analytically (in particular, it is a $\mathbb{R}$-Hilbert
space if $t<0$, while it is a complete $\mathbb{R}$-ISIPS if $t\geq0$),
and a complete semi-normed $*$-algebra over $\mathbb{R}$ operator-algebraically
(in particular, it forms a $C^{*}$-algebra over $\mathbb{R}$ if
$t<0$). It is not difficult to check that each hypercomplex number
$h=x+yi+uj_{t}+vk_{t}$ of $\mathbb{H}_{t}$ has its hypercomplex-conjugate,
\[
h^{\dagger}=x-yi-uj_{t}-vk_{t},\;\mathrm{in\;}\mathbb{H}_{t},
\]
and one can have a well-defined $\mathbb{R}$-linear functional $\tau$
on $\mathbb{H}_{t}$,
\[
\tau\left(h\right)=x.
\]
So, under this new setting, one can define the real part,
\[
Re\left(x+yi+uj_{t}+vk_{t}\right)=x,
\]
and the imaginary part,
\[
Im\left(x+yi+uj_{t}+vk_{t}\right)=yi+uj_{t}+vk_{t},
\]
for all $\mathit{x+yi+uj_{t}+vk_{t}\in\mathbb{H}_{t}}$, which are
naturally understood to be
\[
\left(x+yi\right)+\left(u+vi\right)j_{t},\;\mathrm{since\;}ij_{t}=k_{t}.
\]
Recall that, in {[}3{]}, we generalized hyperbolic numbers in a similar
way.
\begin{defn}
Let $t\in\mathbb{R}$, and $\mathbb{H}_{t}$, the corresponding $t$-scaled
hypercomplex ring. The subring,
\[
\mathbb{D}_{t}=\left\{ \left(x,u\right)\in\mathbb{H}_{t}:x,u\in\mathbb{R}\right\} ,
\]
is called the $t$-scaled hyperbolic (sub)ring, and all elements of
$\mathbb{D}_{t}$ are called $t$-scaled hyperbolic numbers.
\end{defn}

Note that, by definition, this $t$-scaled hyperbolic ring $\mathbb{D}_{t}$
forms a closed subspace of the $t$-hypercomplex $\mathbb{R}$-space
$\mathbf{X}_{t}=\left(\mathbb{H}_{t},\left\langle ,\right\rangle _{t}\right)$,
and it becomes a closed $*$-subalgebra of the complete semi-normed
$*$-algebra $\mathcal{M}_{t}$ over $\mathbb{R}$. By our new notational
setting, the $t$-scaled hyperbolic ring $\mathbb{D}_{t}$ can be
re-defined by

\medskip{}

\hfill{}$\mathbb{D}_{t}=\left\{ x+0i+uj_{t}+0k_{t}\in\mathbb{H}_{t}:x,u\in\mathbb{R}\right\} ,$\hfill{}(3.11)

\medskip{}

\noindent in the $t$-scaled hypercomplexes $\mathbb{H}_{t}$. It
is not hard to check that $\mathbb{D}_{-1}$ is isomorphic to the
complex field,
\[
\mathbb{C}=\left\{ x+ui:x,u\in\mathbb{C},\;i^{2}=-1\right\} ,
\]
and $\mathbb{D}_{1}$ is isomorphic to the classical hyperbolic numbers,
\[
\mathcal{D}=\left\{ x+uj:x,u\in\mathbb{R},\;j^{2}=1\right\} ,
\]
algebraically, and analytically (e.g., see {[}3{]} for details).
\begin{thm}
Let $h=x+uj_{t}\in\mathbb{D}_{t}$. Then there exists $e^{j_{t}\theta}\in\mathbb{D}_{t}$,
such that
\[
h=\left\Vert h\right\Vert _{t}e^{j_{t}\theta},\;\;\;\mathrm{in\;\;\;}\mathbb{D}_{t},
\]
where\hfill{}$\textrm{(3.12)}$
\[
\left\Vert e^{j_{t}\theta}\right\Vert _{t}=1,\;\;\mathrm{with\;\;}\theta=Arg\left(\left(x,u\right)\right)\in\left[0,2\pi\right],
\]
where $\left[0,2\pi\right]$ is the closed interval in $\mathbb{R}$,
and $Arg\left(\left(x,u\right)\right)$ is the argument of the vector
$\left(x,u\right)\in\mathbb{R}^{2}$, and
\[
e^{j_{t}\theta}=\left\{ \begin{array}{ccc}
cos\left(\sqrt{\left|t\right|}\theta\right)+j_{t}\left(\frac{sin\left(\sqrt{\left|t\right|}\theta\right)}{\sqrt{\left|t\right|}}\right) &  & \mathrm{if\;}t<0\\
\\
\pm1+uj_{t}, &  & \mathrm{if\;}t=0\\
\\
cosh\left(\sqrt{t}\theta\right)+j_{t}\left(\frac{sinh\left(\sqrt{t}\theta\right)}{\sqrt{t}}\right) &  & \mathrm{if\;}t>0.
\end{array}\right.
\]
\end{thm}

\begin{proof}
In {[}3{]}, we showed that if $\left(x,u\right)\in\mathbb{D}_{t}$
in the $t$-hypercomplex $\mathbb{R}$-space $\mathbf{X}_{t}$, then
it is expressed by
\[
\left(x,u\right)=\left\Vert \left(x,u\right)\right\Vert _{t}exp^{j_{t}\theta},\;\;\;\mathrm{in\;\;\;}\mathbb{D}_{t},
\]
where\hfill{}(3.13)
\[
exp^{j_{t}\theta}=\left\{ \begin{array}{ccc}
\left(cos\left(\sqrt{\left|t\right|}\theta\right),\:\frac{sin\sqrt{\left|t\right|}\theta}{\sqrt{\left|t\right|}}\right) &  & \mathrm{if\;}t<0\\
\\
\left(\pm1,u\right),\;\mathrm{for\;all\;}u\in\mathbb{R} &  & \mathrm{if}\;t=0\\
\\
\left(cosh\left(\sqrt{t}\theta\right),\;\frac{sinh\left(\sqrt{t}\theta\right)}{\sqrt{t}}\right) &  & \mathrm{if\;}t\geq0,
\end{array}\right.
\]
where $\theta=Arg\left(\left(x,u\right)\right)\in\left[0,2\pi\right]$
is the argument of the point $\left(x,u\right)$ in $\mathbb{R}^{2}$.
So, the $t$-scaled polar decomposition (3.12) holds by (3.11) and
(3.13).
\end{proof}
As in {[}3{]}, if we construct the unit set $\mathbb{T}_{t}=\left\{ h\in\mathbb{D}_{t}:\left\Vert h\right\Vert _{t}=1\right\} $
of all units in the $t$-scaled hyperbolic ring $\mathbb{D}_{t}$,
then
\[
\mathbb{T}_{t}=\left\{ e^{j_{t}\theta}:\theta\in\mathbb{R}\right\} ,
\]
where $e^{j_{t}\theta}$ are in the sense of (3.12) in $\mathbb{D}_{t}$. 

\section{Nonzero-Scaled Regular Functions}

Let $t\in\mathbb{R}$ be a fixed scale, and 
\[
\mathbb{H}_{t}=span_{\mathbb{R}}\left(\left\{ 1,i,j_{t},k_{t}\right\} \right),
\]
the $t$-scaled hypercomplexes. We consider functions,
\[
f:\mathbb{H}_{t}\rightarrow\mathbb{H}_{t},
\]
in the $t$-scaled hypercomplex variable,
\[
w=x_{1}+x_{2}i+x_{3}j_{t}+x_{4}k_{t},\;\mathrm{with\;}x_{l}\in\mathbb{R},
\]
for all $l=1,2,3,4$. In particular, we are interested in the case
where such functions $f$ are $\mathbb{R}$-differentiable in an open
connected set $\Omega\subset\mathbb{H}_{t}$, by regarding $\Omega$
as an open subset of $\mathbb{R}^{4}$. Let 
\[
D_{t}=\frac{\partial}{\partial x_{1}}+i\frac{\partial}{\partial x_{2}}+j_{t}\frac{\partial}{\partial x_{3}}+k_{t}\frac{\partial}{\partial x_{4}}
\]
be an operator acting on the $\mathbb{R}$-differentiable functions
on $\mathbb{H}_{t}$.

\subsection{Motivation: Hyperholomorphic Functions on $\mathbb{H}_{-1}$}

For more about hyperholomorphic theory on the $\left(-1\right)$-scaled
hypercomplexes $\mathbb{H}_{-1}$, see e.g., {[}5{]}, {[}6{]} and
{[}8{]}. Recall that the $\left(-1\right)$-scaled hypercomplexes
$\mathbb{H}_{-1}$ is the noncommutative field of the quaternions
(e.g., see Sections 2 and 3 above, or {[}1{]}, {[}2{]} and {[}3{]}).
Our study is motivated by the main results of {[}4{]}, {[}5{]}, {[}6{]},
{[}8{]}, and {[}9{]}. Recall that the quaternions $\mathbb{H}_{-1}$
is noncommutative for the ($-1$)-scaled multiplication ($\cdot_{-1}$),
and hence, the hyperholomorphic property is considered in two ways
from the left, and from the right. Also, in this quaternionic case,
the differential operator $D_{-1}$ is called the Cauchy-Fueter (differential)
operator (e.g., see {[}5{]}, {[}6{]}, and {[}8{]}).

Define two operators $D_{-1}$ and $D_{-1}^{\dagger}$ by
\[
D_{-1}=\frac{\partial}{\partial x_{1}}+i\frac{\partial}{\partial x_{2}}+j_{-1}\frac{\partial}{\partial x_{3}}+k_{-1}\frac{\partial}{\partial x_{4}},
\]
and
\[
D_{-1}^{\dagger}=\frac{\partial}{\partial x_{1}}-\frac{\partial}{\partial x_{2}}i-\frac{\partial}{\partial x_{3}}j_{-1}-\frac{\partial}{\partial x_{4}}k_{-1},
\]
on the set of all $\mathbb{R}$-differentiable quaternionic functions.
\begin{defn}
A $\mathbb{R}$-differentiable quaternionic function $f:\mathbb{H}_{-1}\rightarrow\mathbb{H}_{-1}$
is left-hyperholomorphic on an open subset $\Omega$ of $\mathbb{H}_{-1}$,
if
\[
D_{-1}f\overset{\textrm{denote}}{=}\frac{\partial f}{\partial x_{1}}+i\frac{\partial f}{\partial x_{2}}+j_{-1}\frac{\partial f}{\partial x_{3}}+k_{-1}\frac{\partial f}{\partial x_{4}}=0;
\]
and $f$ is said to be right-hyperholomorphic on $\Omega$, if
\[
fD_{-1}\overset{\textrm{denote}}{=}\frac{\partial f}{\partial x_{1}}+\frac{\partial f}{\partial x_{2}}i+\frac{\partial f}{\partial x_{3}}j_{-1}+\frac{\partial f}{\partial x_{4}}k_{-1}=0,
\]
in a $\mathbb{H}_{-1}$-variable $w=x_{1}+x_{2}i+x_{3}j_{-1}+x_{4}k_{-1}$
with $x_{1},x_{2},x_{3},x_{4}\in\mathbb{R}$. We simply say that $f$
is hyperholomorphic on $\Omega$, if it is both left and right hyperholomorphic
on $\Omega$.
\end{defn}

Readers can verify that we later define the (left, and right) scaled-hyperholomorphic
property of $\mathbb{R}$-differentiable functions on $\left\{ \mathbb{H}_{t}\right\} _{t\in\mathbb{R}}$
similarly as above. In this section, we concentrate on the left-hyperholomorphic
functions on the quaternions $\mathbb{H}_{-1}$. The study and properties
of right-hyperholomorphic functions would be similar. 

It is well-known that
\[
D_{-1}^{\dagger}D_{-1}=\varDelta_{-1},
\]
with\hfill{}(4.1.1)
\[
\varDelta_{-1}\overset{\textrm{def}}{=}\frac{\partial^{2}}{\partial x_{1}^{2}}+\frac{\partial^{2}}{\partial x_{2}^{2}}+\frac{\partial^{2}}{\partial x_{3}^{2}}+\frac{\partial^{2}}{\partial x_{4}^{2}}.
\]
Thus, hyperholomorphic functions are harmonic on $\mathbb{H}_{-1}$.
\begin{prop}
A left hyperholomorphic function $f:\mathbb{H}_{-1}\rightarrow\mathbb{H}_{-1}$
on an open connected subset $U$ of $\mathbb{H}_{-1}$ is harmonic
on $U$, i.e.,

\medskip{}

\hfill{}$\varDelta_{-1}f=f\varDelta_{-1}=\frac{\partial^{2}f}{\partial x_{1}^{2}}+\frac{\partial^{2}f}{\partial x_{2}^{2}}+\frac{\partial^{2}f}{\partial x_{3}^{2}}+\frac{\partial^{2}f}{\partial x_{4}^{2}}=0,\;\mathrm{on\;}U$.\hfill{}$\textrm{(4.1.2)}$
\end{prop}

\begin{proof}
If $f$ is left-hyperholomorphic on $\mathbb{H}_{-1}$, then
\[
\varDelta_{-1}f=D_{-1}^{\dagger}D_{-1}f=D_{-1}^{\dagger}\left(D_{-1}f\right)=0,
\]
by (4.1.1). So, the harmonicity (4.1.2) holds for a left-hyperholomorphic
function $f$.
\end{proof}
It is shown that every left-hyperholomorphic function $f$ is expressed
by the form,
\[
f\left(w\right)=f\left(0\right)+\overset{3}{\underset{n=1}{\sum}}\zeta_{n+1}\left(w\right)\mathcal{R}_{n+1}f\left(w\right),
\]
with\hfill{}(4.1.3)
\[
\zeta_{n+1}\left(w\right)=x_{n+1}-x_{1}e_{n+1},
\]
and
\[
\mathcal{R}_{n+1}f\left(w\right)=\int_{0}^{1}\frac{\partial f}{\partial x_{n+1}}\left(tw\right)dt,
\]
where
\[
e_{2}=i,\;\;e_{3}=j_{-1},\;\;\mathrm{and\;\;}e_{4}=k_{-1}
\]
(e.g., see {[}5{]}, {[}6{]} and {[}8{]}). Note that the functions
$\zeta_{2},\:\zeta_{3},$ and $\zeta_{4}$ of (4.1.3) are both left
and right entire hyperholomorphic (on $\mathbb{H}_{-1}$), called
the Cauchy-Fueter polynomials on the quaternions $\mathbb{H}_{-1}$.
Are the above harmonicity (4.1.2) and the expansion (4.1.3) generalized
for arbitrary scales $t\in\mathbb{R}$, under a certain ``scaled''
hyperholomorphic property?

Independently, in {[}4{]} and {[}9{]}, the hyperholomorphic property
(the left, or right regularity) and the harmonicity of $\mathbb{R}$-differentiable
functions on open connected subsets of the split-quaternions $\mathbb{H}_{1}$
is studied and characterized, under slightly different differential-operator,
and Laplacian settings based on {[}9{]} and {[}14{]}.

\subsection{Discussion and Assumption}

Let $\mathbb{H}_{t}$ be the $t$-scaled hypercomplexes, for a scale
$t\in\mathbb{R}$.
\begin{defn}
Suppose $\mathcal{T}_{t}$ is the semi-norm topology on $\mathbb{H}_{t}=\mathbf{X}_{t}$,
the $t$-hypercomplex $\mathbb{R}$-space induced by the semi-norm
$\left\Vert .\right\Vert _{t}$ of (2.4.10) under (3.7) (especially,
it becomes a Hilbert-norm topology if $t<0$, meanwhile, it is a semi-norm
topology if $t\geq0$). Define a set,
\[
\mathcal{F}_{t.U}\overset{\textrm{def}}{=}\left\{ f:\mathbb{H}_{t}\rightarrow\mathbb{H}_{t}\left|f\textrm{ is }\mathbb{R}\textrm{-differentiable on }U\in\mathcal{T}_{t}\right.\right\} ,
\]
and the set $\mathcal{F}_{t}\overset{\textrm{def}}{=}\underset{U\in\mathcal{T}_{t}}{\cup}\mathcal{F}_{t,U}$.
\end{defn}

Note that each $\mathbb{R}$-differentiable-function family $\mathcal{F}_{t,U}$
is understood to be

\medskip{}

\hfill{}$\mathcal{F}_{t,U}=\underset{k\in\Lambda}{\sqcup}\mathcal{F}_{t,U_{k}},\;\mathrm{if\;}U=\underset{k\in\Lambda}{\sqcup}U_{k}\in\mathcal{T}_{t},$\hfill{}(4.2.1)

\medskip{}

\noindent where $\left\{ U_{k}\right\} _{k\in\Lambda}\subset\mathcal{T}_{t}$,
for an index set $\Lambda$, where $\sqcup$ is the disjoint union,
where $\left\{ U_{k}\right\} _{k\in\Lambda}$ are connected in $\mathcal{T}_{t}$.
i.e., each $\mathbb{R}$-differentiable-function family $\mathcal{F}_{t,U}$
for $U\in\mathcal{T}_{t}$ is the disjoint union of the collection
$\left\{ \mathcal{F}_{t,U_{k}}\right\} _{k\in\Lambda}$ of $\mathbb{R}$-differentiable-function
families on open connected subsets $\left\{ U_{k}\right\} _{k\in\Lambda}$
of $\mathbb{H}_{t}$, as in (4.2.1). Of course, if $U$ is both open
and connected in $\mathcal{T}_{t}$, then the index set $\Lambda$
in (4.2.1) becomes a singleton set, and hence, $\mathcal{F}_{t,U}$
becomes itself in the sense of (4.2.1).

\medskip{}

\noindent $\mathbf{Assumption.}$ From below, if we take a $\mathbb{R}$-differentiable-function
family $\mathcal{F}_{t,U}$, then the open subset $U\in\mathcal{T}_{t}$
of $\mathbb{H}_{t}$ is automatically regarded as a connected subset,
for convenience.\hfill{}\textifsymbol[ifgeo]{64}

\medskip{}

It is not difficult to check that the $\mathbb{R}$-differentiable-function
family $\mathcal{F}_{t,U}$ of $U\in\mathcal{T}_{t}$ forms a well-defined
$\mathbb{R}$-vector space. Also, the $\mathbb{R}$-differentiable-function
family $\mathcal{F}_{t}$ forms a $\mathbb{R}$-vector space in the
sense that: if $r_{1},r_{2}\in\mathbb{R}$ and $f_{1},f_{2}\in\mathcal{F}_{t}$,
and hence, if $r_{1}f_{1}\in\mathcal{F}_{t,U_{1}}$ and $r_{2}f_{2}\in\mathcal{F}_{t,U_{2}}$,
for $U_{1},U_{2}\in\mathcal{T}_{t}$, then
\[
r_{1}f_{1}+r_{2}f_{2}\in\mathcal{F}_{t,U_{1}\cap U_{t}},\;\mathrm{in\;}\mathcal{F}_{t}.
\]

Before proceeding our works, we discuss the difficulties to maintain
similar settings of Section 4.1 for the cases where we replace the
$\left(-1\right)$-scaled hypercomplexes (which is the quaternions)
$\mathbb{H}_{-1}$ to a general $t$-scaled hypercomplexes $\mathbb{H}_{t}$,
for $t\in\mathbb{R}\setminus\left\{ -1\right\} $. Naturally, we want
to extend the results of Section 4.1 to those on $\mathbb{H}_{t}$,
for any scales $t\in\mathbb{R}$. However, if we use the operators,
\[
D_{t}=\frac{\partial f}{\partial x_{1}}+i\frac{\partial f}{\partial x_{2}}+j_{t}\frac{\partial f}{\partial x_{3}}+k_{t}\frac{\partial f}{\partial x_{4}},
\]
and\hfill{}(4.2.2)
\[
D_{t}^{\dagger}=\frac{\partial f}{\partial x_{1}}-\frac{\partial f}{\partial x_{2}}i-\frac{\partial f}{\partial x_{3}}j_{t}-\frac{\partial f}{\partial x_{4}}k_{t},
\]
for a fixed $t\in\mathbb{R}$, then one ``cannot'' have
\[
D_{t}^{\dagger}D_{t}=\frac{\partial^{2}}{\partial x_{1}^{2}}+\frac{\partial^{2}}{\partial x_{2}^{2}}+\frac{\partial^{2}}{\partial x_{3}^{2}}+\frac{\partial^{2}}{\partial x_{4}^{2}}
\]
in general, as in Section 4.1, especially, where $t\neq-1$ in $\mathbb{R}$
(e.g., see Section 4.3 below). In other words, one cannot have a natural
Laplacian-like operator from the operators (4.2.2). So, we need new
settings for studying the hyperholomorphic property, and harmonicity
of $\mathcal{F}_{t,U}$ on $\mathbb{H}_{t}$.

Also, if a given scale $t$ is zero, i.e., $t=0$ in $\mathbb{R}$,
then the analysis of $\mathcal{F}_{0,U}$ on $\mathbb{H}_{0}$, for
$U\in\mathcal{T}_{0}$, seems totally different from those on $\left\{ \mathbb{H}_{t}\right\} _{t\in\mathbb{R}\setminus\left\{ 0\right\} }$.
So, in the rest of Section 4 below, we first concentrate on the cases
where $t\neq0$ in $\mathbb{R}$. The case where $t=0$ will be considered
in Section 5, independently.

\medskip{}

\noindent $\mathbf{Assumption.}$ In the following Sections 4.3, 4.4
and 4.5, we assume $t\in\mathbb{R}\setminus\left\{ 0\right\} $.\hfill{}\textifsymbol[ifgeo]{64}

\subsection{The Scaled Regularity}

As we have seen in Section 4.1, the $\left(-1\right)$-scaled differential
operator $D_{-1}$ satisfies that
\[
D_{-1}^{\dagger}D_{-1}=\varDelta_{-1}=\frac{\partial^{2}}{\partial x_{1}^{2}}+\frac{\partial^{2}}{\partial x_{2}^{2}}+\frac{\partial^{2}}{\partial x_{3}^{2}}+\frac{\partial^{2}}{\partial x_{4}^{2}},
\]
by (4.1.1), allowing us to study harmonicity independent from the
imaginary $\mathbb{R}$-basis elements $\left\{ i,j_{-1},k_{-1}\right\} $
of $\mathbb{H}_{-1}$. In this section, let's construct ideas and
approaches to extend the results of Section 4.1 to the general cases
where $t\in\mathbb{R}\setminus\left\{ 0\right\} $. Throughout this
section, we let $t\in\mathbb{R}\setminus\left\{ 0\right\} $ as we
assumed in Section 4.2, and $\mathbb{H}_{t}$, the corresponding $t$-scaled
hypercomplexes.

First, consider the case where we have operators $D_{t}$ and $D_{t}^{\dagger}$
on $\mathcal{F}_{t}$ by
\[
D_{t}=\frac{\partial}{\partial x_{1}}+i\frac{\partial}{\partial x_{2}}+j_{t}\frac{\partial}{\partial x_{3}}+k_{t}\frac{\partial}{\partial x_{4}},
\]
and\hfill{}$\textrm{(4.3.1)}$
\[
D_{t}^{\dagger}\overset{\textrm{def}}{=}\frac{\partial}{\partial x_{1}}-\frac{\partial}{\partial x_{2}}i-\frac{\partial}{\partial x_{3}}j_{t}-\frac{\partial}{\partial x_{4}}k_{t},
\]
as in Section 4.1. Observe that

\medskip{}

$\;\;$$D_{t}^{\dagger}D_{t}=\left(\frac{\partial}{\partial x_{1}}-\frac{\partial}{\partial x_{2}}i-\frac{\partial}{\partial x_{3}}j_{t}-\frac{\partial}{\partial x_{4}}k_{t}\right)\left(\frac{\partial}{\partial x_{1}}+i\frac{\partial}{\partial x_{2}}+j_{t}\frac{\partial}{\partial x_{3}}+k_{t}\frac{\partial}{\partial x_{4}}\right)$

\medskip{}

$\;\;\;\;\;\;\;\;\;\;$$=\frac{\partial}{\partial x_{1}}\frac{\partial}{\partial x_{1}}+\frac{\partial}{\partial x_{1}}i\frac{\partial}{\partial x_{2}}+\frac{\partial}{\partial x_{1}}j_{t}\frac{\partial}{\partial x_{3}}+\frac{\partial}{\partial x_{1}}k_{t}\frac{\partial}{\partial x_{4}}$

\medskip{}

$\;\;\;\;\;\;\;\;\;\;\;\;\;\;\;\;$$-\frac{\partial}{\partial x_{2}}i\frac{\partial}{\partial x_{1}}-\frac{\partial}{\partial x_{2}}i^{2}\frac{\partial}{\partial x_{2}}-\frac{\partial}{\partial x_{2}}ij_{t}\frac{\partial}{\partial x_{3}}-\frac{\partial}{\partial x_{2}}ik_{t}\frac{\partial}{\partial x_{4}}$

\medskip{}

$\;\;\;\;\;\;\;\;\;\;\;\;\;\;\;\;$$-\frac{\partial}{\partial x_{3}}j_{t}\frac{\partial}{\partial x_{1}}-\frac{\partial}{\partial x_{3}}j_{t}i\frac{\partial}{\partial x_{2}}-\frac{\partial}{\partial x_{3}}j_{t}^{2}\frac{\partial}{\partial x_{3}}-\frac{\partial}{\partial x_{3}}j_{t}k_{t}\frac{\partial}{\partial x_{4}}$

\medskip{}

$\;\;\;\;\;\;\;\;\;\;\;\;\;\;\;\;$$-\frac{\partial}{\partial x_{4}}k_{t}\frac{\partial}{\partial x_{1}}-\frac{\partial}{\partial x_{4}}k_{t}i\frac{\partial}{\partial x_{2}}-\frac{\partial}{\partial x_{4}}k_{t}j_{t}\frac{\partial}{\partial x_{3}}-\frac{\partial}{\partial x_{4}}k_{t}^{2}\frac{\partial}{\partial x_{4}}$

\medskip{}

$\;\;\;\;\;\;\;\;\;\;$$=\frac{\partial^{2}}{\partial x_{1}^{2}}-\frac{\partial^{2}}{\partial x_{2}^{2}}i^{2}-\frac{\partial}{\partial x_{2}}ij_{t}\frac{\partial}{\partial x_{3}}-\frac{\partial}{\partial x_{2}}ik_{t}\frac{\partial}{\partial x_{4}}$

\medskip{}

$\;\;\;\;\;\;\;\;\;\;\;\;\;\;\;\;\;\;$$-\frac{\partial}{\partial x_{3}}j_{t}i\frac{\partial}{\partial x_{2}}-\frac{\partial}{\partial x_{3}}j_{t}^{2}\frac{\partial}{\partial x_{3}}-\frac{\partial}{\partial x_{3}}j_{t}k_{t}\frac{\partial}{\partial x_{4}}$

\medskip{}

$\;\;\;\;\;\;\;\;\;\;\;\;\;\;\;\;\;\;$$-\frac{\partial}{\partial x_{4}}k_{t}i\frac{\partial}{\partial x_{2}}-\frac{\partial}{\partial x_{4}}k_{t}j_{t}\frac{\partial}{\partial x_{3}}-\frac{\partial}{\partial x_{4}}k_{t}^{2}\frac{\partial}{\partial x_{4}}$

\medskip{}

$\;\;\;\;\;\;\;\;\;\;$$=\frac{\partial^{2}}{\partial x_{1}^{2}}+\frac{\partial^{2}}{\partial x_{2}^{2}}-\frac{\partial}{\partial x_{3}}\left(-ti\right)\frac{\partial}{\partial x_{4}}-\frac{\partial}{\partial x_{4}}\left(ti\right)\frac{\partial}{\partial x_{3}}-t\frac{\partial^{2}}{\partial x_{3}^{2}}-t\frac{\partial^{2}}{\partial x_{4}^{2}},$\hfill{}(4.3.2)

\medskip{}

\noindent by (3.5) and (3.6), i.e.,

\medskip{}

\hfill{}$D_{t}^{\dagger}D_{t}=\frac{\partial^{2}}{\partial x_{1}^{2}}+\frac{\partial^{2}}{\partial x_{2}^{2}}+t\frac{\partial}{\partial x_{3}}i\frac{\partial}{\partial x_{4}}-t\frac{\partial}{\partial x_{4}}i\frac{\partial}{\partial x_{3}}-t\frac{\partial^{2}}{\partial x_{3}^{2}}-t\frac{\partial^{2}}{\partial x_{4}^{2}},$\hfill{}(4.3.3)

\medskip{}

\noindent by (4.3.2).
\begin{prop}
If $D_{t}=\frac{\partial}{\partial x_{1}}+i\frac{\partial}{\partial x_{2}}+j_{t}\frac{\partial}{\partial x_{3}}+k_{t}\frac{\partial}{\partial x_{4}}$
on $\mathcal{F}_{t,U}$, for $U\in\mathcal{T}_{t}$, then

\medskip{}

\hfill{}$D_{t}^{\dagger}D_{t}=\frac{\partial^{2}}{\partial x_{1}^{2}}+\frac{\partial^{2}}{\partial x_{2}^{2}}+t\frac{\partial}{\partial x_{3}}i\frac{\partial}{\partial x_{4}}-t\frac{\partial}{\partial x_{4}}i\frac{\partial}{\partial x_{3}}-t\frac{\partial^{2}}{\partial x_{3}^{2}}-t\frac{\partial^{2}}{\partial x_{4}^{2}},\;\mathrm{on\;}U.$\hfill{}$\textrm{(4.3.4)}$

\end{prop}

\begin{proof}
The proof is done by the straightforward computations (4.3.3).
\end{proof}
Motivated by (4.3.4), we define new differential operators $\nabla_{t}$
and $\nabla_{t}^{\dagger}$ by
\[
\nabla_{t}=\frac{\partial}{\partial x_{1}}+i\frac{\partial}{\partial x_{2}}-j_{t}\frac{sgn\left(t\right)\partial}{\sqrt{\left|t\right|}\partial x_{3}}-k_{t}\frac{sgn\left(t\right)\partial}{\sqrt{\left|t\right|}\partial x_{4}},
\]
and\hfill{}(4.3.5)
\[
\nabla_{t}^{\dagger}=\frac{\partial}{\partial x_{1}}-\frac{\partial}{\partial x_{2}}i+\frac{sgn\left(t\right)\partial}{\sqrt{\left|t\right|}\partial x_{3}}j_{t}+\frac{sgn\left(t\right)\partial}{\sqrt{\left|t\right|}\partial x_{4}}k_{t},
\]
on $\mathcal{F}_{t,U}$, for any $U\in\mathcal{T}_{t}$, where
\[
\frac{sgn\left(t\right)\partial}{\sqrt{\left|t\right|}\partial x_{l}}=\left(\frac{sgn\left(t\right)}{\sqrt{\left|t\right|}}\right)\frac{\partial}{\partial x_{l}},\;\;\mathrm{for\;\;}l=3,4,
\]
where\hfill{}(4.3.6)
\[
sgn\left(t\right)=\left\{ \begin{array}{ccc}
1 &  & \mathrm{if\;}t>0\\
-1 &  & \mathrm{if\;}t<0,
\end{array}\right.
\]
for all $t\in\mathbb{R}\setminus\left\{ 0\right\} $. i.e., by (4.3.5)
and (4.3.6),
\[
\nabla_{t}\overset{\textrm{def}}{=}\left\{ \begin{array}{ccc}
\frac{\partial}{\partial x_{1}}+i\frac{\partial}{\partial x_{2}}-j_{t}\frac{\partial}{\sqrt{t}\partial x_{3}}-k_{t}\frac{\partial}{\sqrt{t}\partial x_{4}} &  & \mathrm{if\;}t>0\\
\\
\frac{\partial}{\partial x_{1}}+\frac{\partial}{\partial x_{2}}i+\frac{\partial}{\sqrt{\left|t\right|}\partial x_{3}}j_{t}+\frac{\partial}{\sqrt{\left|t\right|}\partial x_{4}}k_{t} &  & \mathrm{if\;}t<0.
\end{array}\right.
\]

\begin{lem}
Let $\nabla_{t}$ and $\nabla_{t}^{\dagger}$ be the operators (4.3.6)
on $\mathcal{F}_{t,U}$. Then

\medskip{}

\hfill{}$\nabla_{t}^{\dagger}\nabla_{t}=\frac{\partial^{2}}{\partial x_{1}^{2}}+\frac{\partial^{2}}{\partial x_{2}^{2}}-\frac{sgn\left(t\right)\partial^{2}}{\partial x_{3}^{2}}-\frac{sgn\left(t\right)\partial^{2}}{\partial x_{4}^{2}},$\hfill{}$\textrm{(4.3.7)}$

\medskip{}

\noindent on $\mathcal{F}_{t,U}$.
\end{lem}

\begin{proof}
Consider that, 

\medskip{}

$\nabla_{t}^{\dagger}\nabla_{t}=\left(\frac{\partial}{\partial x_{1}}-\frac{\partial}{\partial x_{2}}i+\frac{sgn\left(t\right)\partial}{\sqrt{\left|t\right|}\partial x_{3}}j_{t}+\frac{sgn\left(t\right)\partial}{\sqrt{\left|t\right|}\partial x_{4}}k_{t}\right)\left(\frac{\partial}{\partial x_{1}}+i\frac{\partial}{\partial x_{2}}-j_{t}\frac{sgn\left(t\right)\partial}{\sqrt{\left|t\right|}\partial x_{3}}-k_{t}\frac{sgn\left(t\right)}{\sqrt{\left|t\right|}\partial x_{4}}\right)$

\medskip{}

$\;\;\;\;$$=\frac{\partial}{\partial x_{1}}\frac{\partial}{\partial x_{1}}+\frac{\partial}{\partial x_{1}}i\frac{\partial}{\partial x_{2}}-\frac{\partial}{\partial x_{1}}j_{t}\frac{sgn\left(t\right)\partial}{\sqrt{\left|t\right|}\partial x_{3}}-\frac{\partial}{\partial x_{1}}k_{t}\frac{sgn\left(t\right)\partial}{\sqrt{\left|t\right|}\partial x_{4}}$

\medskip{}

$\;\;\;\;\;\;\;\;$$-\frac{\partial}{\partial x_{2}}i\frac{\partial}{\partial x_{1}}-\frac{\partial}{\partial x_{2}}i^{2}\frac{\partial}{\partial x_{2}}+\frac{\partial}{\partial x_{2}}ij_{t}\frac{sgn\left(t\right)\partial}{\sqrt{\left|t\right|}\partial x_{3}}+\frac{\partial}{\partial x_{2}}ik_{t}\frac{sgn\left(t\right)\partial}{\sqrt{\left|t\right|}\partial x_{4}}$

\medskip{}

$\;\;\;\;\;\;\;\;$$+\frac{sgn\left(t\right)\partial}{\sqrt{\left|t\right|}\partial x_{3}}j_{t}\frac{\partial}{\partial x_{1}}+\frac{sgn\left(t\right)\partial}{\sqrt{\left|t\right|}\partial x_{3}}j_{t}i\frac{\partial}{\partial x_{2}}-\frac{sgn\left(t\right)\partial}{\sqrt{\left|t\right|}\partial x_{3}}j_{t}^{2}\frac{sgn\left(t\right)\partial}{\sqrt{\left|t\right|}\partial x_{3}}-\frac{sgn\left(t\right)\partial}{\sqrt{\left|t\right|}\partial x_{3}}j_{t}k_{t}\frac{sgn\left(t\right)\partial}{\sqrt{\left|t\right|}\partial x_{4}}$

\medskip{}

$\;\;\;\;\;\;\;\;$$+\frac{sgn\left(t\right)\partial}{\sqrt{\left|t\right|}\partial x_{4}}k_{t}\frac{\partial}{\partial x_{1}}+\frac{sgn\left(t\right)\partial}{\sqrt{\left|t\right|}\partial x_{4}}k_{t}i\frac{\partial}{\partial x_{2}}-\frac{sgn\left(t\right)\partial}{\sqrt{\left|t\right|}\partial x_{4}}k_{t}j_{t}\frac{sgn\left(t\right)\partial}{\sqrt{\left|t\right|}\partial x_{3}}-\frac{sgn\left(t\right)\partial}{\sqrt{\left|t\right|}\partial x_{4}}k_{t}^{2}\frac{sgn\left(t\right)\partial}{\sqrt{\left|t\right|}\partial x_{4}}$

\medskip{}

$\;\;\;\;$$=\frac{\partial^{2}}{\partial x_{1}^{2}}-\frac{\partial}{\partial x_{2}}i^{2}\frac{\partial}{\partial x_{2}}+\frac{\partial}{\partial x_{2}}ij_{t}\frac{sgn\left(t\right)\partial}{\sqrt{\left|t\right|}\partial x_{3}}+\frac{\partial}{\partial x_{2}}ik_{t}\frac{sgn\left(t\right)\partial}{\sqrt{\left|t\right|}\partial x_{4}}$

\medskip{}

$\;\;\;\;\;\;\;\;\;\;\;\;\;$$+\frac{sgn\left(t\right)\partial}{\sqrt{\left|t\right|}\partial x_{3}}j_{t}i\frac{\partial}{\partial x_{2}}-\frac{1}{\left|t\right|}\frac{\partial}{\partial x_{3}}j_{t}^{2}\frac{\partial}{\partial x_{3}}-\frac{1}{\left|t\right|}\frac{\partial}{\partial x_{3}}j_{t}k_{t}\frac{\partial}{\partial x_{4}}$

\medskip{}

$\;\;\;\;\;\;\;\;\;\;\;\;\;$$+\frac{sgn\left(t\right)\partial}{\sqrt{\left|t\right|}\partial x_{4}}k_{t}i\frac{\partial}{\partial x_{2}}-\frac{1}{\left|t\right|}\frac{\partial}{\partial x_{4}}k_{t}j_{t}\frac{\partial}{\partial x_{3}}-\frac{1}{\left|t\right|}\frac{\partial}{\partial x_{4}}k_{t}^{2}\frac{\partial}{\partial x_{4}}$

\medskip{}

\noindent since $\left(sgn\left(t\right)\right)\left(sgn\left(t\right)\right)=1$
in $\left\{ \pm1\right\} $

\medskip{}

$\;\;\;\;$$=\frac{\partial^{2}}{\partial x_{1}^{2}}-\frac{\partial}{\partial x_{2}}i^{2}\frac{\partial}{\partial x_{2}}-\frac{1}{\left|t\right|}\frac{\partial}{\partial x_{3}}j_{t}^{2}\frac{\partial}{\partial x_{3}}-\frac{1}{\left|t\right|}\frac{\partial}{\partial x_{3}}j_{t}k_{t}\frac{\partial}{\partial x_{4}},$

\medskip{}

$\;\;\;\;\;\;\;\;\;\;\;\;\;\;\;\;\;\;\;\;$$-\frac{1}{\left|t\right|}\frac{\partial}{\partial x_{4}}k_{t}j_{t}\frac{\partial}{\partial x_{3}}-\frac{1}{\left|t\right|}\frac{\partial}{\partial x_{4}}k_{t}^{2}\frac{\partial}{\partial x_{4}}$

\medskip{}

$\;\;\;\;$$=\frac{\partial^{2}}{\partial x_{1}^{2}}+\frac{\partial^{2}}{\partial x_{2}^{2}}-\frac{t}{\left|t\right|}\frac{\partial^{2}}{\partial x_{3}^{2}}-\frac{1}{\left|t\right|}\frac{\partial}{\partial x_{3}}\left(-ti\right)\frac{\partial}{\partial x_{4}}-\frac{1}{\left|t\right|}\frac{\partial}{\partial x_{4}}\left(ti\right)\frac{\partial}{\partial x_{3}}-\frac{t}{\left|t\right|}\frac{\partial^{2}}{\partial x_{4}^{2}}$

\medskip{}

$\;\;\;\;$$=\frac{\partial^{2}}{\partial x_{1}^{2}}+\frac{\partial^{2}}{\partial x_{2}^{2}}-\frac{sgn\left(t\right)\partial^{2}}{\partial x_{3}^{2}}+sgn\left(t\right)\frac{\partial}{\partial x_{3}}\left(i\right)\frac{\partial}{\partial x_{4}}-sgn\left(t\right)\frac{\partial}{\partial x_{4}}\left(i\right)\frac{\partial}{\partial x_{3}}-\frac{sgn\left(t\right)\partial^{2}}{\partial x_{4}^{2}}$

\medskip{}

\noindent since $\frac{t}{\left|t\right|}=sgn\left(t\right)$

\medskip{}

$\;\;\;\;$$=\frac{\partial^{2}}{\partial x_{1}^{2}}+\frac{\partial^{2}}{\partial x_{2}^{2}}-\frac{sgn\left(t\right)\partial^{2}}{\partial x_{3}^{2}}-\frac{sgn\left(t\right)\partial^{2}}{\partial x_{4}^{2}}$,

\medskip{}

\noindent by (3.5) and (3.6), implying that
\[
\nabla_{t}^{\dagger}\nabla_{t}=\left(\frac{\partial^{2}}{\partial x_{1}^{2}}+\frac{\partial^{2}}{\partial x_{2}^{2}}-\frac{sgn\left(t\right)\partial^{2}}{\partial x_{3}^{2}}-\frac{sgn\left(t\right)\partial^{2}}{\partial x_{4}^{2}}\right).
\]
\end{proof}
The above lemma shows that if we define a $t$-scaled Laplacian operator
$\varDelta_{t}$ on $\mathcal{F}_{t}$ by
\[
\varDelta_{t}\overset{\textrm{def}}{=}\frac{\partial^{2}}{\partial x_{1}^{2}}+\frac{\partial^{2}}{\partial x_{2}^{2}}-\frac{sgn\left(t\right)\partial}{\partial x_{3}^{2}}-\frac{sgn\left(t\right)\partial}{\partial x_{4}^{2}},
\]
then
\[
\varDelta_{t}=\nabla_{t}^{\dagger}\nabla_{t},
\]
by (4.3.7).
\begin{thm}
If $\varDelta_{t}\overset{\textrm{def}}{=}\frac{\partial^{2}}{\partial x_{1}^{2}}+\frac{\partial^{2}}{\partial x_{2}^{2}}-\frac{sgn\left(t\right)\partial^{2}}{\partial x_{3}^{2}}-\frac{sgn\left(t\right)\partial^{2}}{\partial x_{4}^{2}}$
on $\mathcal{F}_{t}$, then

\medskip{}

\hfill{}$\varDelta_{t}=\nabla_{t}^{\dagger}\nabla_{t}.$\hfill{}$\textrm{(4.3.8)}$
\end{thm}

\begin{proof}
The formula (4.3.8) is shown by (4.3.7).
\end{proof}
For example, consider the formula (4.3.8) for $t=-1$ (i.e., the case
where we have the quaternions $\mathbb{H}_{-1}$). If $t=-1$, then
this formula (4.3.8) becomes
\[
\varDelta_{-1}=\frac{\partial^{2}}{\partial x_{1}^{2}}+\frac{\partial^{2}}{\partial x_{2}^{2}}+\frac{\partial^{2}}{\partial x_{3}^{2}}+\frac{\partial^{2}}{\partial x_{4}^{2}}=\nabla_{-1}^{\dagger}\nabla_{-1},
\]
as in Section 4.1. Observe now the formula (4.3.8) for $t=1$ (i.e.,
the case where we have the split-quaternions $\mathbb{H}_{1}$). If
$t=1$, then the formula (4.3.8) goes to
\[
\varDelta_{1}=\frac{\partial^{2}}{\partial x_{1}^{2}}+\frac{\partial^{2}}{\partial x_{2}^{2}}-\frac{\partial^{2}}{\partial x_{3}^{2}}-\frac{\partial^{2}}{\partial x_{4}^{2}}=\nabla_{1}^{\dagger}\nabla_{1},
\]
as in {[}4{]}.

We now consider how the new differential operator $\nabla_{t}$ of
(4.3.6) acts on the functions,
\[
\zeta_{2}=x_{2}-x_{1}i,\;\zeta_{3}=x_{3}-x_{1}j_{t},\;\mathrm{and\;}\zeta_{4}=x_{4}-x_{1}k_{t}.
\]
Observe that

\medskip{}

$\nabla_{t}\zeta_{2}=\zeta_{2}\nabla_{t}=\frac{\partial\zeta_{2}}{\partial x_{1}}+i\frac{\partial\zeta_{2}}{\partial x_{2}}-j_{t}\frac{sng\left(t\right)\partial\zeta_{2}}{\sqrt{\left|t\right|}\partial x_{3}}-k_{t}\frac{sgn\left(t\right)\partial\zeta_{2}}{\sqrt{\left|t\right|}\partial x_{4}}$

\medskip{}

$\;\;\;\;\;\;\;\;\;\;$$=\frac{\partial(x_{2}-x_{1}i)}{\partial x_{1}}+i\frac{\partial\left(x_{2}-x_{1}i\right)}{\partial x_{2}}-j_{t}\frac{sgn\left(t\right)\partial\left(x_{2}-x_{1}i\right)}{\sqrt{\left|t\right|}\partial x_{3}}-k_{t}\frac{sgn\left(t\right)\partial\left(x_{2}-x_{1}i\right)}{\sqrt{\left|t\right|}\partial x_{4}}$

\medskip{}

$\;\;\;\;\;\;\;\;\;\;$$=\left(-i\right)+\left(1\right)i-\left(0\right)j_{t}-\left(0\right)k_{t}=-i+i=0;$\hfill{}(4.3.9)

\medskip{}

$\nabla_{t}\zeta_{3}=\zeta_{3}\nabla_{t}=\frac{\partial\left(x_{3}-x_{1}j_{t}\right)}{\partial x_{1}}+i\frac{\partial\left(x_{3}-x_{1}j_{t}\right)}{\partial x_{2}}-j_{t}\frac{sgn\left(t\right)\partial\left(x_{3}-x_{1}j_{t}\right)}{\sqrt{\left|t\right|}\partial x_{3}}-k_{t}\frac{sgn\left(t\right)\partial\left(x_{3}-x_{1}j_{t}\right)}{\sqrt{\left|t\right|}\partial x_{4}}$

\medskip{}

$\;\;\;\;\;\;\;\;\;\;\;\;$$=\left(-j_{t}\right)+\left(0\right)i-\left(\frac{sgn\left(t\right)}{\sqrt{\left|t\right|}}\right)j_{t}-\left(0\right)k_{t}=\left(-1-\frac{sgn\left(t\right)}{\sqrt{\left|t\right|}}\right)j_{t};$\hfill{}(4.3.10)

\noindent and

$\nabla_{t}\zeta_{4}=\zeta_{4}\nabla_{t}=\frac{\partial\left(x_{4}-x_{1}k_{t}\right)}{\partial x_{1}}+i\frac{\partial\left(x_{4}-x_{1}k_{t}\right)}{\partial x_{2}}-j_{t}\frac{sgn\left(t\right)\partial\left(x_{4}-x_{1}k_{t}\right)}{\sqrt{\left|t\right|}\partial x_{3}}-k_{t}\frac{sgn\left(t\right)\partial\left(x_{4}-x_{1}k_{t}\right)}{\sqrt{\left|t\right|}\partial x_{4}}$

\medskip{}

$\;\;\;\;\;\;\;\;\;\;\;\;$$=\left(-k_{t}\right)+\left(0\right)i-\left(0\right)j_{t}-\left(\frac{sgn\left(t\right)}{\sqrt{\left|t\right|}}\right)k_{t}=\left(-1-\frac{sgn\left(t\right)}{\sqrt{\left|t\right|}}\right)k_{t};$\hfill{}(4.3.11)
\begin{prop}
Let $\nabla_{t}$ be the differential operator $\textrm{(4.3.6)}$,
and let $\left\{ \zeta_{l}\right\} _{l=2}^{4}\subset\mathcal{F}_{t,\mathbb{H}_{t}}$
be the functions introduced as in the very above paragraph. Then
\[
\nabla_{t}\zeta_{2}=\zeta_{2}\nabla_{t}=0,
\]
and\hfill{}$\textrm{(4.3.12)}$
\[
\nabla_{t}\zeta_{3}=\zeta_{3}\nabla_{t}=\rho j_{t},\;\;\nabla\zeta_{4}=\zeta_{4}\nabla_{t}=\rho k_{t},
\]
where
\[
\rho=-1-\frac{sgn\left(t\right)}{\sqrt{\left|t\right|}},\;\;\;\mathrm{in\;\;\;}\mathbb{R}.
\]
\end{prop}

\begin{proof}
The formula (4.3.12) is proven by (4.3.9), (4.3.10) and (4.3.11).
\end{proof}
The formulas in (4.3.12) show how the operator $\nabla_{t}$ acts
on $span_{\mathbb{R}}\left(\left\{ \zeta_{l}\right\} _{l=2}^{4}\right)$
in $\mathcal{F}_{t,\mathbb{H}_{t}}$.
\begin{defn}
Let $\nabla_{t}$ be the operator (4.3.6) on $\mathcal{F}_{t,U}$,
for $U\in\mathcal{T}_{t}$ in $\mathbb{H}_{t}$, and $f\in\mathcal{F}_{t,U}$.
If 
\[
\nabla_{t}f=\frac{\partial f}{\partial x_{1}}+i\frac{\partial f}{\partial x_{2}}-j_{t}\frac{sgn\left(t\right)\partial f}{\sqrt{\left|t\right|}\partial x_{3}}-k_{t}\frac{sgn\left(t\right)\partial f}{\sqrt{\left|t\right|}\partial x_{4}}=0,
\]
then $f$ is said to be left $t$(-scaled)-regular on $U$. If
\[
f\nabla_{t}=\frac{\partial f}{\partial x_{1}}+\frac{\partial f}{\partial x_{2}}i-\frac{sgn\left(t\right)\partial f}{\sqrt{\left|t\right|}\partial x_{3}}j_{t}-\frac{sgn\left(t\right)\partial f}{\sqrt{\left|t\right|}\partial x_{4}}k_{t}=0,
\]
then $f$ is said to be right $t$(-scaled)-regular on $U$. If $f\in\mathcal{F}_{t,U}$
is both left and right $t$-regular, then it is said to be $t$(-scaled)-regular.

A function $f\in\mathcal{F}_{t,U}$ is called a $t$(-scaled)-harmonic
function, if
\[
\varDelta_{t}f=\frac{\partial^{2}f}{\partial x_{1}^{2}}+\frac{\partial^{2}f}{\partial x_{2}^{2}}-\frac{sgn\left(t\right)\partial^{2}f}{\partial x_{3}^{2}}-\frac{sgn\left(t\right)\partial^{2}f}{\partial x_{4}^{2}}=0,
\]
where 
\[
\varDelta_{t}=\frac{\partial^{2}}{\partial x_{1}^{2}}+\frac{\partial^{2}}{\partial x_{2}^{2}}-\frac{sgn\left(t\right)\partial^{2}}{\partial x_{3}^{2}}-\frac{sgn\left(t\right)\partial^{2}}{\partial x_{4}^{2}}.
\]
is in the sense of (4.3.8). 
\end{defn}

By (4.3.12), one can realize that $\zeta_{2}=x_{2}-x_{1}i$ is $t$-regular,
but, $\zeta_{3}=x_{3}-x_{1}j_{t}$ and $\zeta_{4}=x_{4}-x_{1}k_{t}$
are neither left nor right $t$-regular, in general, especially, if
$t\neq-1$ in $\mathbb{R}\setminus\left\{ 0\right\} $. As in Section
4.1, we concentrate on the left $t$-regular functions. We finish
this section with the following theorem.
\begin{thm}
Let $f\in\mathcal{F}_{t,U}$, for $U\in\mathcal{T}_{t}$ in $\mathbb{H}_{t}$.
If
\[
f\;\mathrm{is\;left\;}t\textrm{-regular,}
\]
then\hfill{}$\textrm{(4.3.13)}$
\[
f\;\mathrm{is\;}t\textrm{-harmonic on }U.
\]
\end{thm}

\begin{proof}
Suppose $f$ is left $t$-regular in $\mathcal{F}_{t,U}$, i.e., $\nabla_{t}f=0$.
Then $\varDelta_{t}f=0$, because
\[
\varDelta_{t}f=\nabla_{t}^{\dagger}\left(\nabla_{t}f\right)=\nabla_{t}^{\dagger}\left(0\right)=0,
\]
implying that it is $t$-harmonic on $U$.
\end{proof}
The above theorem shows that the left $t$-regularity implies the
$t$-harmonicity on $U\in\mathcal{T}_{t}$ in $\mathbb{H}_{t}$, by
(4.3.13).

\subsection{Certain $t$-Regular Functions on $\mathbb{H}_{t}$}

Throughout this section, we automatically assume that a given scale
$t$ is non-zero, i.e., $t\in\mathbb{R}\setminus\left\{ 0\right\} $.
Let $\nabla_{t}$ and $\nabla_{t}^{\dagger}$ be the operators (4.3.6),
\[
\nabla_{t}=\frac{\partial}{\partial x_{1}}+i\frac{\partial}{\partial x_{2}}-j_{t}\frac{sgn\left(t\right)\partial}{\sqrt{\left|t\right|}\partial x_{3}}-k_{t}\frac{sgn\left(t\right)\partial}{\sqrt{\left|t\right|}\partial x_{4}},
\]
and
\[
\nabla_{t}^{\dagger}=\frac{\partial}{\partial x_{1}}-\frac{\partial}{\partial x_{2}}i+\frac{sgn\left(t\right)\partial}{\sqrt{\left|t\right|}\partial x_{3}}j_{t}+\frac{sgn\left(t\right)\partial}{\sqrt{\left|t\right|}\partial x_{4}}k_{t},
\]
on $\mathcal{F}_{t,U}$, for any $U\in\mathcal{T}_{t}$ in the $t$-scaled
hypercomplexes $\mathbb{H}_{t}$, and let
\[
\varDelta_{t}=\frac{\partial^{2}}{\partial x_{1}^{2}}+\frac{\partial^{2}}{\partial x_{2}^{2}}-\frac{sgn\left(t\right)\partial^{2}}{\partial x_{3}^{2}}-\frac{sgn\left(t\right)\partial^{2}}{\partial x_{4}^{2}},
\]
on $\mathcal{F}_{t,U}$, satisfying 
\[
\varDelta_{t}=\nabla_{t}^{\dagger}\nabla_{t},\;\;\;\;\mathrm{on\;\;\;\;}\mathcal{F}_{t,U},
\]
by (4.3.8). Motivated by (4.3.12) and (4.3.13), define entire $\mathbb{R}$-differentiable
functions $\left\{ \eta_{l}\right\} _{l=2}^{4}$ on $\mathbb{H}_{t}$
by
\[
\eta_{2}\left(w\right)=x_{2}-x_{1}i,
\]
and\hfill{}(4.4.1)
\[
\eta_{3}\left(w\right)=x_{3}+\frac{sgn\left(t\right)x_{1}}{\sqrt{\left|t\right|}}j_{k},\;\eta_{4}\left(w\right)=x_{4}+\frac{sgn\left(t\right)x_{1}}{\sqrt{\left|t\right|}}k_{t},
\]
in a $\mathbb{H}_{t}$-variable $w=x_{1}+x_{2}i+x_{3}j_{t}+x_{4}k_{t}$,
with $x_{1},x_{2},x_{3},x_{4}\in\mathbb{R}$.

Observe first that
\[
\eta_{2}=x_{2}-x_{1}i\;\mathrm{in\;}\mathcal{F}_{t,\mathbb{H}_{t}},
\]
is $t$-regular in $\mathcal{F}_{t,\mathbb{H}_{t}}$, i.e., \hfill{}(4.4.2)
\[
\eta_{2}\nabla_{t}=0=\nabla_{t}\eta_{2},\;\;\mathrm{on\;\;}\mathbb{H}_{t},
\]
by (4.3.12), since $\eta_{2}$ is identical to the function $\zeta_{2}$
introduced in Section 4.3.
\begin{lem}
If $\eta_{2}=x_{2}-x_{1}i\in\mathcal{F}_{1,\mathbb{H}_{t}}$ is in
the sense of (4.4.1), then

\medskip{}

\hfill{}$\eta_{2}\;\mathrm{is\;a\;}t\textrm{-harmonic }t\textrm{-regular function on }\mathbb{H}_{t}.$\hfill{}$\textrm{(4.4.3)}$
\end{lem}

\begin{proof}
The function $\eta_{2}$ is $t$-regular on $\mathbb{H}_{t}$ by (4.4.2).
Therefore, by (4.3.13), it is $t$-harmonic on $\mathbb{H}_{t}$,
too.
\end{proof}
From below, denote the quantity $\sqrt{\left|t\right|}$ by $\rho$
and $sgn\left(t\right)$ by $s_{t}\in\left\{ \pm\right\} $, and let
$\eta_{3}=x_{3}+\frac{s_{t}x_{1}}{\rho}j_{t}\in\mathcal{F}_{t,\mathbb{H}_{t}}$
be in the sense of (4.4.1). Observe that

\medskip{}

$\;\;\;$$\nabla_{t}\eta_{3}=\frac{\partial\left(x_{3}+\frac{s_{t}x_{1}}{\rho}j_{t}\right)}{\partial x_{1}}+i\frac{\partial\left(x_{3}+\frac{s_{t}x_{1}}{\rho}j_{t}\right)}{\partial x_{2}}-j_{t}\frac{s_{t}\partial\left(x_{3}+\frac{sx_{t1}}{\rho}j_{t}\right)}{\rho\partial x_{3}}-k_{t}\frac{s_{t}\partial\left(x_{3}+\frac{s_{t}x_{1}}{\rho}j_{t}\right)}{\rho\partial x_{4}}$

\medskip{}

$\;\;\;\;\;\;\;\;\;\;\;\;$$=\left(\frac{s_{t}}{\rho}j_{t}\right)+i\left(0\right)-j_{t}\left(\frac{s_{t}}{\rho}\right)-k_{t}\left(0\right)=\frac{s_{t}}{\rho}j_{t}-\frac{s_{t}}{\rho}j_{t}=0,$

\medskip{}

\noindent and\hfill{}(4.4.4)

\medskip{}

$\;\;\;$$\eta_{3}\nabla_{t}=\frac{\partial\left(x_{3}+\frac{s_{t}x_{1}}{\rho}j_{t}\right)}{\partial x_{1}}+\frac{\partial\left(x_{3}+\frac{s_{t}x_{1}}{\rho}j_{t}\right)}{\partial x_{2}}i-\frac{s_{t}\partial\left(x_{3}+\frac{s_{t}x_{1}}{\rho}j_{t}\right)}{\rho\partial x_{3}}j_{t}-\frac{s_{t}\partial\left(x_{3}+\frac{s_{t}x_{1}}{\rho}j_{t}\right)}{\rho\partial x_{4}}k_{t}$

\medskip{}

$\;\;\;\;\;\;\;\;\;\;\;\;$$=\left(\frac{s_{t}}{\rho}j_{t}\right)+\left(0\right)i-\left(\frac{s_{t}}{\rho}\right)j_{t}-\left(0\right)k_{t}=\frac{s_{t}}{\rho}j_{t}-\frac{s_{t}}{\rho}j_{t}=0$;

\medskip{}

\noindent also, if $\eta_{4}=x_{4}+s_{t}\rho x_{1}k_{t}\in\mathcal{F}_{1,\mathbb{H}_{t}}$
is in the sense of (4.4.1), then

\medskip{}

$\;\;\;$$\nabla_{t}\eta_{4}=\frac{\partial\left(x_{4}+\frac{s_{t}x_{1}}{\rho}k_{t}\right)}{\partial x_{1}}+i\frac{\partial\left(x_{4}+\frac{s_{t}x_{1}}{\rho}k_{t}\right)}{\partial x_{2}}-j_{t}\frac{s_{t}\partial\left(x_{4}+\frac{s_{t}x_{1}}{\rho}k_{t}\right)}{\rho\partial x_{3}}-k_{t}\frac{s_{t}\partial\left(x_{4}+\frac{s_{t}x_{1}}{\rho}k_{t}\right)}{\rho\partial x_{4}}$

\medskip{}

$\;\;\;\;\;\;\;\;\;\;\;\;$$=\left(\frac{s_{t}}{\rho}k_{t}\right)+i\left(0\right)-j_{t}\left(0\right)-k_{t}\left(\frac{s_{t}}{\rho}\right)=\frac{s_{t}}{\rho}k_{t}-\frac{s_{t}}{\rho}k_{t}=0,$

\medskip{}

\noindent and\hfill{}(4.4.5)

\medskip{}

$\;\;\;$$\eta_{4}\nabla_{t}=\frac{\partial\left(x_{4}+\frac{s_{t}x_{1}}{\rho}k_{t}\right)}{\partial x_{1}}+\frac{\partial\left(x_{4}+\frac{s_{t}x_{1}}{\rho}k_{t}\right)}{\partial x_{2}}i-\frac{s_{t}\partial\left(x_{4}+\frac{x_{1}}{\rho}k_{t}\right)}{\rho\partial x_{3}}j_{t}-\frac{s_{t}\partial\left(x_{4}+\frac{x_{1}}{\rho}k_{t}\right)}{\rho\partial x_{4}}k_{t}$

\medskip{}

$\;\;\;\;\;\;\;\;\;\;\;\;$$=\left(\frac{s_{t}}{\rho}k_{t}\right)+\left(0\right)i-\left(0\right)j_{t}-\left(\frac{s_{t}}{\rho}\right)k_{t}=\frac{s_{t}}{\rho}k_{t}-\frac{s_{t}}{\rho}k_{t}=0$.

\medskip{}

\begin{lem}
If $\eta_{3}=x_{3}+\frac{sng\left(t\right)x_{1}}{\sqrt{\left|t\right|}}j_{t}$,
and $\eta_{4}=x_{4}+\frac{sgn\left(t\right)x_{1}}{\sqrt{\left|t\right|}}k_{t}$
are in the sense of (4.4.1) in $\mathcal{F}_{1,\mathbb{H}_{t}}$,
then

\medskip{}

\hfill{}$\eta_{3}\;\mathrm{and\;}\eta_{4}\;\mathrm{are\;}t\textrm{-harmonic }t\textrm{-regular on }\mathbb{H}_{t}$.\hfill{}$\textrm{(4.4.6)}$
\end{lem}

\begin{proof}
The functions $\eta_{3},\eta_{4}\in\mathcal{F}_{1,\mathbb{H}_{t}}$
are both left and right $t$-regular on $\mathbb{H}_{t}$, by (4.4.4)
and (4.4.5). Therefore, the functions $\eta_{3}$ and $\eta_{4}$
are $t$-harmonic on $\mathbb{H}_{t}$, too, by (4.3.13).
\end{proof}
By the above two lemmas, one obtains the following result.
\begin{thm}
The functions $\left\{ \eta_{l}\right\} _{l=2}^{4}$ of (4.4.1) are
$t$-harmonic $t$-regular functions on $\mathbb{H}_{t}$.
\end{thm}

\begin{proof}
The proof is done by (4.4.3) and (4.4.6).
\end{proof}
The following corollary is an immediate consequence of the above theorem.
\begin{cor}
Let $f=s_{1}\eta_{2}+s_{2}\eta_{3}+s_{3}\eta_{4}\in\mathcal{F}_{t,\mathbb{H}_{t}}$
be a $\mathbb{R}$-linear combination of the $t$-regular functions
$\left\{ \eta_{l}\right\} _{l=2}^{4}$ of (4.4.1), where $s_{1},s_{2},s_{3}\in\mathbb{R}$.
Then $f$ is not only $t$-regular, but also $t$-harmonic on $\mathbb{H}_{t}$.
\end{cor}

\begin{proof}
Let $f=\overset{3}{\underset{n=1}{\sum}}s_{n}\eta_{n+1}\in\mathcal{F}_{1,\mathbb{H}_{t}}$
be a $\mathbb{R}$-linear combination of $\left\{ \eta_{l}\right\} _{l=2}^{4}$,
for $s_{1},s_{2},s_{3}\in\mathbb{R}$. Then
\[
\nabla_{t}f=\overset{3}{\underset{n=1}{\sum}}s_{n}\left(\nabla_{t}\eta_{n+1}\right)=0,
\]
and
\[
f\nabla_{t}=\overset{3}{\underset{n=1}{\sum}}s_{n}\left(\eta_{n+1}\nabla_{t}\right)=0,
\]
by (4.4.2), (4.4.4) and (4.4.5). Therefore, this function $f$ is
$t$-regular on $\mathbb{H}_{t}$. So, by (4.3.13), this function
$f$ is $t$-harmonic on $\mathbb{H}_{t}$, too.
\end{proof}

\subsection{$t$-Regular Functions of $\mathcal{F}_{t}$}

In this section, we consider $t$-regular functions in $\mathcal{F}_{t}$
more in detail. As in Sections 4.3 and 4.4, throughout this section,
we assume that $t\in\mathbb{R}\setminus\left\{ 0\right\} $. Recall
that the functions,
\[
\eta_{2}\left(w\right)=x_{2}-x_{1}i,
\]
and
\[
\eta_{3}\left(w\right)=x_{3}+\frac{sgn\left(t\right)x_{1}}{\sqrt{\left|t\right|}}j_{t},\;\;\eta\left(w\right)=x_{4}+\frac{sgn\left(t\right)x_{1}}{\sqrt{\left|t\right|}}k_{t},
\]
of (4.4.1) in a $\mathbb{H}_{t}$-variable $w=x_{1}+x_{2}i+x_{3}j_{t}+x_{4}k_{t}$,
with $x_{1},x_{2},x_{3},x_{4}\in\mathbb{R}$, are $t$-harmonic $t$-regular
functions on $\mathbb{H}_{t}$. Therefore, all $\mathbb{R}$-linear
combinations $f=\overset{3}{\underset{n=1}{\sum}}s_{n}\eta_{n+1}\in\mathcal{F}_{t,\mathbb{H}_{t}}$
of $\left\{ \eta_{l}\right\} _{l=2}^{4}\subset\mathcal{F}_{t,\mathbb{H}_{t}}$,
for $s_{1},s_{2},s_{3}\in\mathbb{R}$, are $t$-harmonic $t$-regular
functions on $\mathbb{H}_{t}$. We here consider $t$-regular functions
induced by $\left\{ \eta_{l}\right\} _{l=2}^{4}$. To do that we define
and study the following new operation ($\times$) on the $t$-scaled
hypercomplexes $\mathbb{H}_{t}$.
\begin{defn}
Let $h_{1},...,h_{N}\in\mathbb{H}_{t}$, for $N\in\mathbb{N}$. Then
the symmetrized product of $h_{1},...,h_{N}$ is defined by a new
hypercomplex number,

\medskip{}

\hfill{}$\overset{N}{\underset{n=1}{\times}}h_{n}\overset{\textrm{denote}}{=}h_{1}\times...\times h_{N}\overset{\textrm{def}}{=}\frac{1}{N!}\underset{\sigma\in S_{N}}{\sum}h_{\sigma(1)}\cdot_{t}h_{\sigma(2)}\cdot_{t}...\cdot_{t}h_{\sigma\left(N\right)},$\hfill{}(4.5.1)

\medskip{}

\noindent where $S_{N}$ is the symmetric (or, permutation) group
over $\left\{ 1,...,N\right\} $, where ($\cdot_{t}$) is the $t$-scaled
multiplication (2.1.1) on $\mathbb{H}_{t}$.
\end{defn}

\medskip{}

\noindent $\mathbf{Notation.}$ From below, if there are no confusions,
we denote the operation ($\cdot_{t}$) simply by ($\cdot$) for convenience.
i.e., the above definition (4.5.1) can be re-written to be

\medskip{}

\hfill{}$\overset{N}{\underset{n=1}{\times}}h_{n}=\frac{1}{N!}\underset{\sigma\in S_{N}}{\sum}h_{\sigma(1)}h_{\sigma(2)}...h_{\sigma(N)}$,\hfill{}(4.5.2)

\medskip{}

\noindent for $h_{1},...,h_{N}\in\mathbb{H}_{t}$, for $t\in\mathbb{R}\setminus\left\{ 0\right\} $.\hfill{}\textifsymbol[ifgeo]{64}

\medskip{}

\noindent $\mathbf{Remark.}$ The above symmetrized product ($\times$)
of (4.5.1) is well-defined for ``all'' $t\in\mathbb{R}$, including
the case where $t=0$. Of course, as we mentioned, we restrict our
interests to the cases where $t\neq0$ in $\mathbb{R}$ in this section,
but the above symmetrized product of (4.5.1) is also defined for the
case where $t=0$, too (See Section 5 below).\hfill{}\textifsymbol[ifgeo]{64}

\medskip{}

Remark that, actually, the above definition (4.5.1) of the symmetrized
product implies the following case. If we consider
\[
h^{(n)}\overset{\textrm{denote}}{=}\underset{n\textrm{-times}}{\underbrace{h\times h\times...\times h}},\;\;\mathrm{for\;\;}n\in\mathbb{N},
\]
then, all permutations $\sigma$ of $S_{n}$ induce exactly same element,
so,
\[
h^{(n)}=\frac{1}{n!}\underset{\sigma\in S_{n}}{\sum}\left(h_{\sigma(1)}...h_{\sigma(n)}\right),\;\mathrm{with\;}h_{\sigma(j)}=h,\;\forall j,
\]
in $\mathbb{H}_{t}$, implying that

\medskip{}

\hfill{}$h^{(n)}=\frac{1}{n!}\underset{\sigma\in S_{n}}{\sum}h_{\sigma\left(1\right)}...h_{\sigma\left(n\right)}=\frac{1}{n!}\left(n!h^{n}\right)=h^{n},$\hfill{}(4.5.3)

\medskip{}

\noindent where $h^{n}=\:\underset{n\textrm{-times}}{\underbrace{hh.......h}}\:\overset{\textrm{denote}}{=}\:\underset{n\textrm{-times}}{\underbrace{h\cdot_{t}h\cdot_{t}...\cdot_{t}h}}$,
for all $n\in\mathbb{N}$, by (4.5.2). So, if we consider ``mutually
distinct'' $h_{1}$ and $h_{2}$ in $\mathbb{H}_{t}$, and
\[
h_{1}^{(n)}\times h_{2}\;\;\;\mathrm{in\;\;\;}\mathbb{H}_{t},
\]
then
\[
h_{1}^{(n)}\times h_{2}=\underset{n\textrm{-times}}{\underbrace{h_{1}\times...\times h_{1}}}\times h_{2},
\]
satisfying
\[
h_{1}^{(n)}\times h_{2}=\frac{1}{(n+1)!}\underset{\sigma\in S_{n+1}}{\sum}\left(\left(h_{1}\right)_{\sigma(1)}...\left(h_{1}\right)_{\sigma(n)}\right)\left(h_{2}\right)_{\sigma(n+1)},
\]
and hence,
\[
h_{1}^{(n)}\times h_{2}=\frac{n!}{\left(n+1\right)!}\underset{\sigma\in S_{n+1}}{\sum}\left(h_{1}^{n}\right)\left(h_{2}\right)_{\sigma(n+1)},
\]
by (4.5.3). More generally, if $h_{1},...,h_{N}\in\mathbb{H}_{t}$
are ``mutually distinct,'' and $n_{1},...,n_{N}\in\mathbb{N}$,
then

\medskip{}

\hfill{}$\overset{N}{\underset{j=1}{\times}}h_{j}^{(n_{j})}=\frac{\overset{N}{\underset{j=1}{\prod}}\left(n_{j}!\right)}{\left(\overset{N}{\underset{j=1}{\sum}}n_{j}\right)!}\left(\underset{\sigma\in S_{\overset{N}{\underset{j=1}{\sum}}n_{j}}}{\sum}h_{\sigma\left(1\right)}...h_{\sigma\left(\overset{N}{\underset{j=1}{\sum}}n_{j}\right)}\right),$\hfill{}(4.5.4)

\medskip{}

\noindent by (4.5.3), where $h_{\sigma\left(j\right)}\in\left\{ h_{1},...,h_{N}\right\} $,
for all $j\in\left\{ 1,...,\overset{n}{\underset{j=1}{\sum}}n_{j}\right\} $.
\begin{prop}
Let $h_{1},...,h_{N}$ be mutually distinct elements of the $t$-scaled
hypercomplexes $\mathbb{H}_{t}$, for $N\in\mathbb{N}\setminus\left\{ 1\right\} $,
and let $n_{1},...,n_{N}\in\mathbb{N}$. If
\[
n\overset{\textrm{denote}}{=}\overset{N}{\underset{j=1}{\sum}}n_{j}=n_{1}+...+n_{N}\in\mathbb{N},
\]
then
\[
\overset{N}{\underset{j=1}{\times}}h^{(n_{j})}=\frac{\overset{N}{\underset{j=1}{\prod}}n_{j}}{n!}\left(\underset{\sigma\in S_{n}}{\sum}h_{\sigma\left(1\right)}h_{\sigma\left(2\right)}...h_{\sigma\left(n\right)}\right),
\]
with\hfill{}$\textrm{(4.5.5)}$
\[
h_{j}^{(k)}=\:\underset{k\textrm{-times}}{\underbrace{h_{j}\times h_{j}\times...\times h_{j}}}=h_{j}^{k},\;\forall j=1,...,N,
\]
for all $k\in\left\{ 1,...,n\right\} $, where $h_{\sigma\left(j\right)}\in\left\{ h_{1},...,h_{N}\right\} $,
for all $j=1,...,n$.
\end{prop}

\begin{proof}
The refined computation (4.5.5) of the definition (4.5.1) (denoted
simply by (4.5.2)) is obtained by (4.5.4).
\end{proof}
Let $f_{1},...,f_{N}:\mathbb{H}_{t}\rightarrow\mathbb{H}_{t}$ be
functions for $N\in\mathbb{N}$. Then, similar to (4.5.1) (expressed
by (4.5.2)), one can define a symmetrized-product function of them
by
\[
\overset{N}{\underset{n=1}{\times}}f_{n}=\frac{1}{N!}\underset{\sigma\in S_{N}}{\sum}f_{\sigma(1)}f_{\sigma(2)}...f_{\sigma(N)},
\]
where\hfill{}(4.5.6)
\[
\left(f_{\sigma(1)}f_{\sigma(2)}...f_{\sigma(N)}\right)\left(h\right)=f_{\sigma(1)}\left(h\right)\cdot_{t}...\cdot_{t}f_{\sigma\left(N\right)}\left(h\right),
\]
simply denoted by
\[
\left(f_{\sigma(1)}f_{\sigma(2)}...f_{\sigma(N)}\right)\left(h\right)=f_{\sigma(1)}\left(h\right)...f_{\sigma\left(N\right)}\left(h\right),
\]
as in (4.5.2), for all $\sigma\in S_{N}$ and $h\in\mathbb{H}_{t}$.
Also, similar to (4.5.3), one write 
\[
\mathit{f_{j}^{(n)}}=\underset{n\textrm{-times}}{\underbrace{f_{j}\times f_{j}\times...\times f_{j}}}=f_{j}^{n},\;\;\forall n\in\mathbb{N},
\]
in terms of (4.5.6). Then, for any $h\in\mathbb{H}_{t}$, the image
$\left(\overset{N}{\underset{n=1}{\times}}f_{n}^{(n_{j})}\right)\left(h\right)$
is expressed similar to (4.5.5), if we replace $h_{\sigma\left(j\right)}$
to $f_{\sigma\left(j\right)}\left(h\right)$ in (4.5.5). Without loss
of generality, let's axiomatize that
\[
f^{(0)}=1,\;\mathrm{the\;constant\;\textrm{1-}function\;on\;}\mathbb{H}_{t},
\]
for all functions $f:\mathbb{H}_{t}\rightarrow\mathbb{H}_{t}$.
\begin{defn}
Let $\mathbf{n}\overset{\textrm{denote}}{=}\left(n_{1},n_{2},n_{3}\right)\in\mathbb{N}_{0}^{3}$
be a triple of numbers in $\mathbb{N}_{0}=\mathbb{N}\cup\left\{ 0\right\} $,
and let $\left\{ \eta_{l}\right\} _{l=2}^{4}$ be the $t$-harmonic
$t$-regular functions of (4.4.1). Define a new function $\eta^{\mathbf{n}}$
by
\[
\eta^{\mathbf{n}}\overset{\textrm{def}}{=}\frac{1}{\mathbf{n}!}\left(\eta_{2}^{(n_{1})}\times\eta_{3}^{(n_{2})}\times\eta_{4}^{(n_{3})}\right),
\]
where\hfill{}(4.5.7)
\[
\mathbf{n}!=\left(n_{1}!\right)\left(n_{2}!\right)\left(n_{3}!\right)\in\mathbb{N},
\]
and
\[
\eta_{l+1}^{(n_{l})}=\;\underset{n_{l}\textrm{-times}}{\underbrace{\eta_{l+1}\times\eta_{l+1}\times...\times\eta_{l+1}}}=\eta_{l+1}^{n_{j}},\;\forall l=1,2,3.
\]
\end{defn}

For an arbitrary $\mathbf{n}=\left(n_{1},n_{2},n_{3}\right)\in\mathbb{N}_{0}^{3}$,
let $\eta^{\mathbf{n}}:\mathbb{H}_{t}\rightarrow\mathbb{H}_{t}$ be
a function (4.5.7). By the very construction, this function $\eta^{\mathbf{n}}$
is contained in $\mathcal{F}_{1,\mathbb{H}_{t}}$. Observe that, if
we let
\[
e_{2}=-i,\;e_{3}=\frac{sgn\left(t\right)}{\sqrt{\left|t\right|}}j_{t},\;\mathrm{and\;}e_{4}=\frac{sgn\left(t\right)}{\sqrt{\left|t\right|}}k_{t},
\]
making
\[
\eta_{l}=x_{l}+x_{1}e_{l},\;\;\forall l=2,3,4,
\]
then

\medskip{}

$\;\;\;\;$$\left(\mathbf{n}!\right)\nabla_{t}\eta^{\mathbf{n}}=\left(\mathbf{n}!\right)\nabla_{t}\left(\frac{1}{\mathbf{n}!}\left(\eta_{2}^{\times n_{1}}\times\eta_{3}^{\times n_{2}}\times\eta_{4}^{\times n_{3}}\right)\right)$

\medskip{}

$\;\;\;\;\;\;\;\;\;\;$$=\underset{\sigma\in S_{N}}{\sum}\overset{4}{\underset{l=2}{\sum}}\underset{k\in\left\{ 1,...,N\right\} ,\sigma\left(k\right)=l}{\sum}\eta_{\sigma\left(1\right)}...\eta_{\sigma\left(k-1\right)}e_{l}\eta_{\sigma\left(k+1\right)}...\times\eta_{\sigma\left(N\right)}$

\medskip{}

$\;\;\;\;\;\;\;\;\;\;\;\;\;\;$$-\underset{\sigma\in S_{N}}{\sum}\overset{4}{\underset{l=2}{\sum}}\underset{k\in\left\{ 1,...,N\right\} ,\sigma\left(k\right)=l}{\sum}e_{l}\eta_{\sigma\left(1\right)}...\times\eta_{\sigma\left(k-1\right)}\eta_{\sigma\left(k+1\right)}...\eta_{\sigma\left(N\right)},$\hfill{}(4.5.8)

\medskip{}

\noindent where
\[
N=n_{1}+n_{2}+n_{3}\;\;\mathrm{in\;\;}\mathbb{N}.
\]
Now, note that

\medskip{}

$\;\;\;\;$$0=\underset{\sigma\in S_{N}}{\sum}\overset{4}{\underset{l=2}{\sum}}\underset{k\in\left\{ 1,...,N\right\} ,\sigma\left(k\right)=l}{\sum}x_{l}\eta_{\sigma\left(1\right)}...\eta_{\sigma\left(k-1\right)}\eta_{\sigma\left(k+1\right)}...\times\eta_{\sigma\left(N\right)}$

\medskip{}

$\;\;\;\;\;\;\;\;\;\;\;\;$$-\underset{\sigma\in S_{N}}{\sum}\overset{4}{\underset{l=2}{\sum}}\underset{k\in\left\{ 1,...,N\right\} ,\sigma\left(k\right)=l}{\sum}\eta_{\sigma\left(1\right)}...\times\eta_{\sigma\left(k-1\right)}x_{l}\eta_{\sigma\left(k+1\right)}...\eta_{\sigma\left(N\right)}$,\hfill{}(4.5.9)

\medskip{}

\noindent by the formula (7) in the page 99 of {[}8{]}. Let's multiply
$x_{1}$ to the formula (4.5.8). Then

\medskip{}

$\;\;$$x_{1}\left(\mathbf{n}!\right)\nabla_{t}\eta^{\mathbf{n}}=x_{1}\left(n!\right)\nabla_{t}\eta^{\mathbf{n}}+0$

\medskip{}

$\;\;\;\;$$=x_{1}\left(\underset{\sigma\in S_{N}}{\sum}\overset{4}{\underset{l=2}{\sum}}\underset{k\in\left\{ 1,...,N\right\} ,\sigma\left(k\right)=l}{\sum}\eta_{\sigma\left(1\right)}...\eta_{\sigma\left(k-1\right)}e_{l}\eta_{\sigma\left(k+1\right)}...\times\eta_{\sigma\left(N\right)}\right.$

\medskip{}

$\;\;\;\;\;\;\;\;$$\left.-\underset{\sigma\in S_{N}}{\sum}\overset{4}{\underset{l=2}{\sum}}\underset{k\in\left\{ 1,...,N\right\} ,\sigma\left(k\right)=l}{\sum}e_{l}\eta_{\sigma\left(1\right)}...\times\eta_{\sigma\left(k-1\right)}\eta_{\sigma\left(k+1\right)}...\eta_{\sigma\left(N\right)}\right)+0$

\medskip{}

\noindent by (4.5.8)

\medskip{}

$\;\;\;\;$$=\underset{\sigma\in S_{N}}{\sum}\overset{4}{\underset{l=2}{\sum}}\underset{k\in\left\{ 1,...,N\right\} ,\sigma\left(k\right)=l}{\sum}\eta_{\sigma\left(1\right)}...\eta_{\sigma\left(k-1\right)}\left(x_{1}e_{l}\right)\eta_{\sigma\left(k+1\right)}...\times\eta_{\sigma\left(N\right)}$

\medskip{}

$\;\;\;\;\;\;\;\;$$-\underset{\sigma\in S_{N}}{\sum}\overset{4}{\underset{l=2}{\sum}}\underset{k\in\left\{ 1,...,N\right\} ,\sigma\left(k\right)=l}{\sum}\left(x_{1}e_{l}\right)\eta_{\sigma\left(1\right)}...\times\eta_{\sigma\left(k-1\right)}\eta_{\sigma\left(k+1\right)}...\eta_{\sigma\left(N\right)}$

\medskip{}

$\;\;\;\;\;\;\;\;$$+\underset{\sigma\in S_{N}}{\sum}\overset{4}{\underset{l=2}{\sum}}\underset{k\in\left\{ 1,...,N\right\} ,\sigma\left(k\right)=l}{\sum}x_{l}\eta_{\sigma\left(1\right)}...\eta_{\sigma\left(k-1\right)}\eta_{\sigma\left(k+1\right)}...\times\eta_{\sigma\left(N\right)}$

\medskip{}

$\;\;\;\;\;\;\;\;$$-\underset{\sigma\in S_{N}}{\sum}\overset{4}{\underset{l=2}{\sum}}\underset{k\in\left\{ 1,...,N\right\} ,\sigma\left(k\right)=l}{\sum}\eta_{\sigma\left(1\right)}...\times\eta_{\sigma\left(k-1\right)}x_{l}\eta_{\sigma\left(k+1\right)}...\eta_{\sigma\left(N\right)}$

\medskip{}

\noindent by (4.5.9)

\medskip{}

$\;\;\;\;$$=\left(\underset{\sigma\in S_{N}}{\sum}\overset{4}{\underset{l=2}{\sum}}\underset{k\in\left\{ 1,...,N\right\} ,\sigma\left(k\right)=l}{\sum}\eta_{\sigma\left(1\right)}...\eta_{\sigma\left(k-1\right)}\left(x_{1}e_{l}\right)\eta_{\sigma\left(k+1\right)}...\times\eta_{\sigma\left(N\right)}\right.$

\medskip{}

$\;\;\;\;\;\;\;\;$$\left.+\underset{\sigma\in S_{N}}{\sum}\overset{4}{\underset{l=2}{\sum}}\underset{k\in\left\{ 1,...,N\right\} ,\sigma\left(k\right)=l}{\sum}x_{l}\eta_{\sigma\left(1\right)}...\eta_{\sigma\left(k-1\right)}\eta_{\sigma\left(k+1\right)}...\times\eta_{\sigma\left(N\right)}\right)$

\medskip{}

$\;\;\;\;\;\;\;\;$$+\left(-\underset{\sigma\in S_{N}}{\sum}\overset{4}{\underset{l=2}{\sum}}\underset{k\in\left\{ 1,...,N\right\} ,\sigma\left(k\right)=l}{\sum}\left(x_{1}e_{l}\right)\eta_{\sigma\left(1\right)}...\times\eta_{\sigma\left(k-1\right)}\eta_{\sigma\left(k+1\right)}...\eta_{\sigma\left(N\right)}\right.$

\medskip{}

$\;\;\;\;\;\;\;\;$$\left.-\underset{\sigma\in S_{N}}{\sum}\overset{4}{\underset{l=2}{\sum}}\underset{k\in\left\{ 1,...,N\right\} ,\sigma\left(k\right)=l}{\sum}\eta_{\sigma\left(1\right)}...\times\eta_{\sigma\left(k-1\right)}x_{l}\eta_{\sigma\left(k+1\right)}...\eta_{\sigma\left(N\right)}\right)$

\medskip{}

$\;\;\;\;$$=\underset{\sigma\in S_{N}}{\sum}\overset{4}{\underset{l=2}{\sum}}\underset{k\in\left\{ 1,...,N\right\} ,\sigma\left(k\right)=l}{\sum}\eta_{\sigma\left(1\right)}...\eta_{\sigma\left(k-1\right)}\eta_{l}\eta_{\sigma\left(k+1\right)}...\times\eta_{\sigma\left(N\right)}$

\medskip{}

$\;\;\;\;\;\;\;\;$$-\underset{\sigma\in S_{N}}{\sum}\overset{4}{\underset{l=2}{\sum}}\underset{k\in\left\{ 1,...,N\right\} ,\sigma\left(k\right)=l}{\sum}\eta_{l}\eta_{\sigma\left(1\right)}...\eta_{\sigma\left(k-1\right)}\eta_{\sigma\left(k+1\right)}...\times\eta_{\sigma\left(N\right)}$

\medskip{}

$\;\;\;\;$$=0,$\hfill{}(4.5.10)

\medskip{}

\noindent under the sum over $S_{N}$.
\begin{lem}
For any $\mathbf{n}=\left(n_{1},n_{2},n_{3}\right)\in\mathbb{N}_{0}^{3}$,
the function \textup{$\eta^{\mathbf{n}}$ of $\textrm{(4.5.7)}$is
$t$-regular on $\mathbb{H}_{t}$, i.e.,}

\medskip{}

\hfill{}$\nabla_{t}\eta^{\mathbf{n}}=0=\eta^{\mathbf{n}}\nabla_{t}.$\hfill{}$\textrm{(4.5.11)}$
\end{lem}

\begin{proof}
The first equality of the formula (4.5.11) is obtained by (4.5.10).
The second equality is obtained similarly by (4.5.10), by the definition
of the symetrized product (4.5.2) and (4.5.7), with help of the $t$-regularity
(4.4.3) of $\left\{ \eta_{l}\right\} _{l=2}^{4}$.
\end{proof}
By the above lemma, we obtain the following result.
\begin{thm}
Let $\mathbf{n}=\left(n_{1},n_{2},n_{3}\right)\in\mathbb{N}_{0}^{3}$,
and let $\eta^{\mathbf{n}}\in\mathcal{F}_{t,\mathbb{H}_{t}}$ be the
function (4.5.7). Then it is a $t$-harmonic $t$-regular function,
i.e.,

\medskip{}

\hfill{}$\eta^{\mathbf{n}}\nabla_{t}=\nabla_{t}\eta^{\mathbf{n}}=0,\;\mathrm{and\;}$
$\varDelta_{t}\eta^{\mathbf{n}}=0,\;\mathrm{on\;}\mathbb{H}_{t}.$\hfill{}$\textrm{(4.5.12)}$
\end{thm}

\begin{proof}
The $t$-regularity of $\eta^{\mathbf{n}}$ on $\mathbb{H}_{t}$ is
shown by (4.5.11), for all $\mathbf{n}\in\mathbb{N}_{0}^{3}$. Thus,
by (4.3.13), the functions $\left\{ \eta^{\mathbf{n}}\right\} _{\mathbf{n}\in\mathbb{N}_{0}^{3}}$
are $t$-harmonic on $\mathbb{H}_{t}$. Therefore, the relation (4.5.12)
holds.
\end{proof}
By the above theorem, all functions $\left\{ \eta^{\mathbf{n}}:\mathbf{n}\in\mathbb{N}_{0}^{3}\right\} $
are $t$-harmonic $t$-regular functions on $\mathbb{H}_{t}$. By
(4.5.12), one obtains the following result.

Now, recall the complete semi-norm $\left\Vert .\right\Vert _{t}$
of (2.4.10) on the $t$-scaled hypercomplexes $\mathbb{H}_{t}$,
\[
\left\Vert h\right\Vert _{t}=\sqrt{\left|\left\langle h,h\right\rangle _{t}\right|}=\sqrt{\left|\tau\left(hh^{\dagger}\right)\right|}=\sqrt{\left|\left|a\right|^{2}-t\left|b\right|^{2}\right|},
\]
for all $h=a+bj_{t}\in\mathbb{H}_{t}$, regarded as $\left(a,b\right)\in\mathbb{H}_{t}$
with $a,b\in\mathbb{C}$ (in the sense of Section 2) satisfying (2.4.11)
and (2.4.12). i.e., if $t<0$, then it is a complete norm on $\mathbb{H}_{t}$,
while if $t>0$, then it is a complete semi-norm on $\mathbb{H}_{t}$.
If we understand,
\[
\eta_{2}=x_{2}-x_{1}i,\;\eta_{3}=x_{3}+\frac{sgn\left(t\right)x_{1}}{\sqrt{\left|t\right|}}j_{t},\;\eta_{4}=x_{4}+\frac{sgn\left(t\right)x_{1}}{\sqrt{\left|t\right|}}k_{t},
\]
as their images of $\mathbb{H}_{t}$, then they are regarded as
\[
\eta_{2}=\left(x_{2}-x_{1}i,\:0\right),\;\eta_{3}=\left(x_{3},\;\frac{sgn\left(t\right)x_{1}}{\sqrt{\left|t\right|}}\right),
\]
and
\[
\eta_{4}=\left(x_{4},\;\frac{sgn\left(t\right)x_{1}}{\sqrt{\left|t\right|}}\right),
\]
in ``the $t$-scaled hypercomplex ring $\mathbb{H}_{t}$ in the sense
of Section 2.'' So, one can compute their norms on $\mathbb{H}_{t}$,
\[
\left\Vert \eta_{2}\right\Vert _{t}=\sqrt{\left|\left|x_{2}-x_{1}i\right|^{2}-t\left|0\right|^{2}\right|}=\sqrt{\left|x_{2}^{2}+\left(-x_{1}\right)^{2}\right|}=\sqrt{x_{1}^{2}+x_{2}^{2}},
\]
\[
\left\Vert \eta_{3}\right\Vert _{t}=\sqrt{\left|\left|x_{3}\right|^{2}-t\left|\frac{sgn\left(t\right)x_{1}}{\sqrt{\left|t\right|}}\right|^{2}\right|}=\sqrt{\left|x_{3}^{2}-sgn\left(t\right)x_{1}^{2}\right|},
\]
and\hfill{}(4.5.13)
\[
\left\Vert \eta_{4}\right\Vert _{t}=\sqrt{\left|\left|x_{4}\right|^{2}-t\left|\frac{sgn\left(t\right)x_{1}}{\sqrt{\left|t\right|}}\right|^{2}\right|}=\sqrt{\left|x_{4}^{2}-sgn\left(t\right)x_{1}^{2}\right|}.
\]

\begin{lem}
Let $\eta^{\mathbf{n}}\in\mathcal{F}_{t,\mathbb{H}_{t}}$ be a function
$\mathrm{(4.5.7)}$ for $\mathbf{n}=\left(n_{1},n_{2},n_{3}\right)\in\mathbb{N}_{0}^{3}$.
Then

\medskip{}

\hfill{}$\left\Vert \eta^{\mathbf{n}}\right\Vert _{t}\leq\left(\sqrt{x_{2}^{2}+x_{1}^{2}}\right)^{n_{1}}\left(\sqrt{\left|x_{3}^{2}-sgn\left(t\right)x_{1}^{2}\right|}\right)^{n_{2}}\left(\sqrt{\left|x_{4}^{2}-sgn\left(t\right)x_{1}^{2}\right|}\right)^{n_{3}}$.\hfill{}$\mathrm{(4.5.14)}$
\end{lem}

\begin{proof}
By (4.5.7), for any $\mathbf{n}=\left(n_{1},n_{2},n_{3}\right)\in\mathbb{N}_{0}^{3}$,

\medskip{}

$\;\;$$\left\Vert \eta^{\mathbf{n}}\right\Vert _{t}=\left\Vert \frac{1}{\mathbf{n}!}\left(\eta_{2}^{(n_{1})}\times\eta_{3}^{(n_{2})}\times\eta_{4}^{(n_{3})}\right)\right\Vert _{t}$

\medskip{}

$\;\;\;\;$$=\left\Vert \frac{1}{\left(n_{1}!\right)\left(n_{2}!\right)\left(n_{3}!\right)}\left(\frac{\left(n_{1}!\right)\left(n_{2}!\right)\left(n_{3}!\right)}{n!}\underset{\sigma\in S_{n}}{\sum}h_{\sigma\left(1\right)}h_{\sigma\left(2\right)}...h_{\sigma\left(n\right)}\right)\right\Vert _{t}$

\medskip{}

\noindent where $h_{\sigma\left(l\right)}\in\left\{ \eta_{2},\eta_{3},\eta_{4}\right\} $,
for all $\sigma\in S_{n}$, with $n=n_{1}+n_{2}+n_{3}\in\mathbb{N}$

\medskip{}

$\;\;\;\;$$=\left\Vert \eta_{2}^{n_{1}}\eta_{3}^{n_{2}}\eta_{4}^{n_{3}}\right\Vert _{t}\leq\left\Vert \eta_{2}\right\Vert _{t}^{n_{1}}\left\Vert \eta_{3}\right\Vert _{t}^{n_{2}}\left\Vert \eta_{4}\right\Vert _{t}^{n_{3}}$

\medskip{}

\noindent by (4.5.5)

\medskip{}

$\;\;\;\;$$=\left(\sqrt{x_{2}^{2}+x_{1}^{2}}\right)^{n_{1}}\left(\sqrt{\left|x_{3}^{2}-sgn\left(t\right)x_{1}^{2}\right|}\right)^{n_{2}}\left(\sqrt{\left|x_{4}^{2}-sgn\left(t\right)x_{1}^{2}\right|}\right)^{n_{3}},$

\medskip{}

\noindent by (4.5.13). Therefore, the inequality (4.5.14) holds.
\end{proof}
The above lemma shows that, for any arbitrarily fixed $\mathbf{n}\in\mathbb{N}_{0}^{3}$,
the corresponding function $\eta^{\mathbf{n}}$ is bounded by (4.5.14).
\begin{thm}
Let $t\neq0$ in $\mathbb{R}$, and $f\in\mathcal{F}_{t,U}$, a $\mathbb{R}$-differentiable
function, where $U\in\mathcal{T}_{t}$ containing $0=0+0i+0j_{t}+0k_{t}\in\mathbb{H}_{t}$.
If $f$ is $\mathbb{R}$-analytic on $U$, then
\[
f\textrm{ is left }t\textrm{-regular on }U,
\]
if and only if\hfill{}$\mathrm{(4.5.15)}$
\[
f=f\left(0\right)+\underset{\mathbf{n}\in\mathbb{N}^{3}}{\sum}\eta^{\mathbf{n}}f_{\mathbf{n}},
\]
with
\[
f_{\mathbf{n}}=\frac{1}{\mathbf{n}!}\frac{\partial^{n_{1}+n_{2}+n_{3}}f}{\partial x_{2}^{n_{1}}\partial x_{3}^{n_{2}}\partial x_{4}^{n_{3}}}\left(0\right),\;\forall\mathbf{n}\in\mathbb{N}^{3}.
\]
\end{thm}

\begin{proof}
Clearly, if $f=f\left(0\right)+\underset{\mathbf{n}\in\mathbb{N}^{3}}{\sum}\eta^{\mathbf{n}}f_{\mathbf{n}}$,
then it is left $t$-regular, by the $t$-regularity (4.5.12) of $\left\{ \eta^{\mathbf{n}}:\mathbf{n}\in\mathbb{N}_{0}^{3}\right\} $.

Suppose $f$ is left $t$-regular in $\mathcal{F}_{t,U}$, satisfying
$\nabla_{t}f=0$. Then
\[
f\left(w\right)-f\left(0\right)=\overset{4}{\underset{n=2}{\sum}}\left(\eta_{n}\left(w\right)\right)\left(\left(R_{n}f\right)\left(w\right)\right),
\]
where\hfill{}(4.5.16)
\[
\left(R_{n}f\right)\left(w\right)=\int_{0}^{1}\frac{\partial f\left(tw\right)}{\partial x_{n}}dt,\;\forall n=2,3,4,
\]
in a $\mathbb{H}_{t}$-variable $w=x_{1}+x_{2}i+x_{3}j_{t}+x_{4}k_{t}$,
with $x_{1},x_{2},x_{3},x_{4}\in\mathbb{R}$. Indeed, 
\[
\frac{df\left(tw\right)}{dt}=x_{1}\frac{\partial f\left(tw\right)}{\partial x_{1}}+\overset{4}{\underset{l=2}{\sum}}x_{l}\frac{\partial f\left(tw\right)}{\partial x_{l}},
\]
identical to
\[
\frac{df\left(tw\right)}{dt}=x_{1}\left(-i\frac{\partial}{\partial x_{2}}+j_{t}\frac{sgn\left(t\right)\partial}{\sqrt{\left|t\right|}\partial x_{3}}+k_{t}\frac{sgn\left(t\right)\partial}{\sqrt{\left|t\right|}\partial x_{4}}\right)f\left(tw\right)+\overset{4}{\underset{l=2}{\sum}}x_{l}\frac{\partial f\left(tw\right)}{\partial x_{l}},
\]
implying that
\[
\frac{df\left(tw\right)}{dt}=\overset{4}{\underset{l=2}{\sum}}\eta_{l}\left(w\right)\frac{\partial f\left(tw\right)}{\partial x_{l}}.
\]
Thus, by iterating (4.5.16), one can get that
\[
f=f\left(0\right)+\underset{\mathbf{n}\in\mathbb{N}^{3}}{\sum}\eta^{\mathbf{n}}f_{\mathbf{n}}
\]
where
\[
f_{\mathbf{n}}=\frac{1}{\mathbf{n}!}\frac{\partial^{n_{1}+n_{2}+n_{3}}f}{\partial x_{2}^{n_{1}}\partial x_{3}^{n_{2}}\partial x_{4}^{n_{3}}}\left(0\right),\;\forall\mathbf{n}\in\mathbb{N}^{3}.
\]
Note that, such an iteration can be done by the assumption that $f$
is $\mathbb{R}$-analytic on $U$, by the boundedness condition (4.5.14).
Therefore, the characterization (4.5.15) holds true.
\end{proof}
Remark that, if a given scale $t$ is negatice, i.e., $t<0$ in $\mathbb{R}$,
then the $\mathbb{R}$-regularity automatically implies the $\mathbb{R}$-analyticity
by (4.3.5), (4.3.6) and (4.5.14). So, if $t<0$, then the above theorem
holds without the $\mathbb{R}$-analyticity assumption for $f\in\mathcal{F}_{t,U}$.

\section{$0$-Regular Functions on $\mathbb{H}_{0}$}

In Section 4, we studied the $t$-regularity and the $t$-hamonicity
of $\mathbb{R}$-differentiable functions of $\mathcal{F}_{t,U}$,
for $U\in\mathcal{T}_{t}$ in the $t$-scaled hypercomplexes $\mathbb{H}_{t}$,
where a given scale is nonzero, i.e., $t\in\mathbb{R}\setminus\left\{ 0\right\} $.
To do that, we defined the operators $\nabla_{t}$ and $\nabla_{t}^{\dagger}$
on $\mathcal{F}_{t}$ by
\[
\nabla_{t}=\frac{\partial}{\partial x_{1}}+i\frac{\partial}{\partial x_{2}}-j_{t}\frac{sgn\left(t\right)\partial}{\sqrt{\left|t\right|}\partial x_{3}}-k_{t}\frac{sgn\left(t\right)\partial}{\sqrt{\left|t\right|}\partial x_{4}},
\]
and\hfill{}(5.1)
\[
\nabla_{t}^{\dagger}=\frac{\partial}{\partial x_{1}}-\frac{\partial}{\partial x_{2}}i+\frac{sgn\left(t\right)\partial}{\sqrt{\left|t\right|}\partial x_{3}}j_{t}+\frac{sgn\left(t\right)\partial}{\sqrt{\left|t\right|}\partial x_{4}}k_{t},
\]
satisfying\hfill{}(5.2)
\[
\varDelta_{t}=\frac{\partial^{2}}{\partial x_{1}^{2}}+\frac{\partial^{2}}{\partial x_{2}^{2}}-\frac{sgn\left(t\right)\partial^{2}}{\partial x_{3}^{2}}-\frac{sgn\left(t\right)\partial^{2}}{\partial x_{4}^{2}}=\nabla_{t}^{\dagger}\nabla_{t},
\]
in a $\mathbb{H}_{t}$-variable $w=x_{1}+x_{2}i+x_{3}j_{t}+x_{4}k_{t}$,
with $x_{1},x_{2},x_{3},x_{4}\in\mathbb{R}$. We showed in Section
4 that the functions,

\medskip{}

\hfill{}$\eta_{2}\left(w\right)=x_{2}-x_{1}i,\;\eta_{3}\left(w\right)=x_{3}+\frac{sgn\left(t\right)x_{1}}{\sqrt{\left|t\right|}}j_{t},\;\eta_{4}+\frac{sgn\left(t\right)x_{1}}{\sqrt{\left|t\right|}},$\hfill{}(5.3)

\medskip{}

\noindent are $t$-harmonic $t$-regular in $\mathcal{F}_{t,\mathbb{H}_{t}}$,
furthermore, all functions of $\mathcal{F}_{t,\mathbb{H}_{t}}$ formed
by

\medskip{}

\hfill{}$\eta^{\mathbf{n}}=\frac{1}{\mathbf{n}!}\left(\eta_{2}^{\left(n_{1}\right)}\times\eta_{3}^{\left(n_{2}\right)}\times\eta_{4}^{\left(n_{3}\right)}\right),\;\mathrm{with\;}\mathbf{n}!=\overset{3}{\underset{l=1}{\prod}}\left(n_{l}!\right)$\hfill{}(5.4)

\medskip{}

\noindent are $t$-harmonic $t$-regular in $\mathcal{F}_{t,\mathbb{H}_{t}}$,
for all $\mathbf{n}=\left(n_{1},n_{2},n_{3}\right)\in\mathbb{N}_{0}^{3}$,
characterizing the left $t$-regularity (4.5.13) on $\mathcal{F}_{t}$,
whenever $t\neq0$ in $\mathbb{R}$.

In this section, we consider the case where the given scale $t$ is
zero in $\mathbb{R}$, i.e., we study $0$-regularity on $\mathcal{F}_{0}$
acting on the $0$-scaled hypercomplexes $\mathbb{H}_{0}$. Recall
that if
\[
h=a+bi+cj_{0}+dk_{0}\in\mathbb{H}_{0},\;\mathrm{with\;}a,b,c,d\in\mathbb{R},
\]
then
\[
h=\left(a+bi\right)+\left(c+di\right)j_{0}=z_{1}+z_{2}j_{0},
\]
with
\[
z_{1}=a+bi,\;\mathrm{and\;}z_{2}=c+di,\;\mathrm{in\;}\mathbb{C},
\]
realized to be
\[
\left[\left(z_{1},z_{2}\right)\right]_{0}=\left(\begin{array}{cc}
z_{1} & 0\cdot z_{2}\\
\overline{z_{2}} & \overline{z_{1}}
\end{array}\right)=\left(\begin{array}{cc}
z_{1} & 0\\
\overline{z_{2}} & \overline{z_{1}}
\end{array}\right),\;\mathrm{in\;}\mathcal{H}_{2}^{0},
\]
satisfying
\[
\mathbb{H}_{0}=span_{\mathbb{R}}\left\{ 1,i,j_{0},k_{0}\right\} .
\]
So, the $\mathbb{R}$-basis elements $\left\{ i,j_{0},k_{0}\right\} $
of $\mathbb{H}_{0}$ satisfy

\medskip{}

\hfill{}$i^{2}=-1,\;\;j_{0}^{2}=0=k_{0}^{2},$\hfill{}(5.5)

\medskip{}

\noindent and the commuting diagrams,
\[
\begin{array}{ccccc}
 &  & i\\
 & ^{1}\swarrow &  & \nwarrow^{-0=0}\\
 & j_{0} & \underset{1}{\longrightarrow} & k_{0} & ,
\end{array}\;\mathrm{and\;}\begin{array}{ccccc}
 &  & i\\
 & ^{0}\nearrow &  & \searrow^{-1}\\
 & j_{0} & \underset{-1}{\longleftarrow} & k_{0} & ,
\end{array}
\]
saying that
\[
ij_{0}=k_{0},\;j_{0}k_{0}=-0i=0,\;k_{0}i=j_{0},
\]
and\hfill{}(5.6)
\[
j_{0}i=-k_{0},\;ik_{0}=-j_{0},\;k_{0}j_{0}=0i=0,
\]
by (3.4), (3.5) and (3.6).

Define first the operators $\nabla_{0}$ and $\nabla_{0}^{\dagger}$
by
\[
\nabla_{0}=\frac{\partial}{\partial x_{1}}+i\frac{\partial}{\partial x_{2}}+j_{0}\frac{\partial}{\partial x_{3}}+k_{0}\frac{\partial}{\partial x_{4}},
\]
and\hfill{}(5.7)
\[
\nabla_{0}^{\dagger}=\frac{\partial}{\partial x_{1}}-\frac{\partial}{\partial x_{2}}i-\frac{\partial}{\partial x_{3}}j_{0}-\frac{\partial}{\partial x_{4}}k_{0},
\]
on $\mathcal{F}_{0,U}$, for any $U\in\mathcal{T}_{t}$ in $\mathbb{H}_{0}$,
similar to (5.1) in a $\mathbb{H}_{0}$-variable $x_{1}+x_{2}i+x_{3}j_{0}+x_{4}k_{0}$,
with $x_{1},x_{2},x_{3},x_{4}\in\mathbb{R}$. Observe that

\medskip{}

$\;\;\;\;$$\nabla_{0}^{\dagger}\nabla_{0}=\left(\frac{\partial}{\partial x_{1}}-\frac{\partial}{\partial x_{2}}i-\frac{\partial}{\partial x_{3}}j_{0}-\frac{\partial}{\partial x_{4}}k_{0}\right)\left(\frac{\partial}{\partial x_{1}}+i\frac{\partial}{\partial x_{2}}+j_{0}\frac{\partial}{\partial x_{3}}+k_{0}\frac{\partial}{\partial x_{4}}\right)$

\medskip{}

$\;\;\;\;\;\;\;\;\;\;\;\;$$=\frac{\partial^{2}}{\partial x_{1}^{2}}+\frac{\partial}{\partial x_{1}}i\frac{\partial}{\partial x_{2}}+\frac{\partial}{\partial x_{1}}j_{0}\frac{\partial}{\partial x_{3}}+\frac{\partial}{\partial x_{1}}k_{0}\frac{\partial}{\partial x_{4}}$

\medskip{}

$\;\;\;\;\;\;\;\;\;\;\;\;\;\;\;\;\;\;$$-\frac{\partial}{\partial x_{2}}i\frac{\partial}{\partial x_{1}}-\frac{\partial}{\partial x_{2}}i^{2}\frac{\partial}{\partial x_{2}}-\frac{\partial}{\partial x_{2}}ij_{0}\frac{\partial}{\partial x_{3}}-\frac{\partial}{\partial x_{2}}ik_{0}\frac{\partial}{\partial x_{4}}$

\medskip{}

$\;\;\;\;\;\;\;\;\;\;\;\;\;\;\;\;\;\;$$-\frac{\partial}{\partial x_{3}}j_{0}\frac{\partial}{\partial x_{1}}-\frac{\partial}{\partial x_{3}}j_{0}i\frac{\partial}{\partial x_{2}}-\frac{\partial}{\partial x_{3}}j_{0}^{2}\frac{\partial}{\partial x_{3}}-\frac{\partial}{\partial x_{3}}j_{0}k_{0}\frac{\partial}{\partial x_{4}}$

\medskip{}

$\;\;\;\;\;\;\;\;\;\;\;\;\;\;\;\;\;\;$$-\frac{\partial}{\partial x_{4}}k_{0}\frac{\partial}{\partial x_{1}}-\frac{\partial}{\partial x_{4}}k_{0}i\frac{\partial}{\partial x_{2}}-\frac{\partial}{\partial x_{4}}k_{0}j_{0}\frac{\partial}{\partial x_{3}}-\frac{\partial}{\partial x_{4}}k_{0}^{2}\frac{\partial}{\partial x_{4}}$

\medskip{}

$\;\;\;\;\;\;\;\;\;\;\;\;$$=\frac{\partial^{2}}{\partial x_{1}^{2}}+\frac{\partial^{2}}{\partial x_{2}^{2}}-\frac{\partial}{\partial x_{2}}ij_{0}\frac{\partial}{\partial x_{3}}-\frac{\partial}{\partial x_{2}}ik_{0}\frac{\partial}{\partial x_{4}}$

\medskip{}

$\;\;\;\;\;\;\;\;\;\;\;\;\;\;\;\;\;\;\;\;\;\;\;\;$$-\frac{\partial}{\partial x_{3}}j_{0}i\frac{\partial}{\partial x_{2}}-\frac{\partial}{\partial x_{3}}j_{0}^{2}\frac{\partial}{\partial x_{3}}-\frac{\partial}{\partial x_{3}}j_{0}k_{0}\frac{\partial}{\partial x_{4}}$

\medskip{}

$\;\;\;\;\;\;\;\;\;\;\;\;\;\;\;\;\;\;\;\;\;\;\;\;$$-\frac{\partial}{\partial x_{4}}k_{0}i\frac{\partial}{\partial x_{2}}-\frac{\partial}{\partial x_{4}}k_{0}j_{0}\frac{\partial}{\partial x_{3}}-\frac{\partial}{\partial x_{4}}k_{0}^{2}\frac{\partial}{\partial x_{4}}$

\medskip{}

$\;\;\;\;\;\;\;\;\;\;\;\;$$=\frac{\partial^{2}}{\partial x_{1}^{2}}+\frac{\partial^{2}}{\partial x_{2}^{2}}-\frac{\partial}{\partial x_{3}}j_{0}^{2}\frac{\partial}{\partial x_{3}}-\frac{\partial}{\partial x_{3}}j_{0}k_{0}\frac{\partial}{\partial x_{4}}$

\medskip{}

$\;\;\;\;\;\;\;\;\;\;\;\;\;\;\;\;\;\;\;\;\;\;\;\;$$-\frac{\partial}{\partial x_{4}}k_{0}j_{0}\frac{\partial}{\partial x_{3}}-\frac{\partial}{\partial x_{4}}k_{0}^{2}\frac{\partial}{\partial x_{4}}$

\medskip{}

$\;\;\;\;\;\;\;\;\;\;\;\;$$=\frac{\partial^{2}}{\partial x_{1}^{2}}+\frac{\partial^{2}}{\partial x_{2}^{2}}-\frac{\partial}{\partial x_{3}}\left(0\right)\frac{\partial}{\partial x_{3}}-\frac{\partial}{\partial x_{3}}\left(-0i\right)\frac{\partial}{\partial x_{4}}-\frac{\partial}{\partial x_{4}}\left(0i\right)\frac{\partial}{\partial x_{3}}-\frac{\partial}{\partial x_{4}}\left(0\right)\frac{\partial}{\partial x_{4}}$

\medskip{}

$\;\;\;\;\;\;\;\;\;\;\;\;$$=\frac{\partial^{2}}{\partial x_{1}^{2}}+\frac{\partial^{2}}{\partial x_{2}^{2}}$,\hfill{}(5.8)

\medskip{}

\noindent by (5.5) and (5.6).
\begin{thm}
Let $\nabla_{0}$ be the operator (5.7) on $\mathcal{F}_{0}$. Then

\medskip{}

\hfill{}$\nabla_{0}^{\dagger}\nabla_{0}=\frac{\partial^{2}}{\partial x_{1}^{2}}+\frac{\partial^{2}}{\partial x_{2}^{2}}.$\hfill{}$\textrm{(5.9)}$
\end{thm}

\begin{proof}
The formula (5.9) is obtained by the straight-forward computations
(5.8).
\end{proof}
One can recognize that the formula (5.9) is similar to the formula
(5.2). However, we do not have $x_{l}$-depending double-partial-derivatives,
for $l=3,4$, in (5.9). It shows that, in the $0$-scaled case, the
Laplacian on $\mathcal{F}_{0,U}$ for an open connected subset $U\subseteq\mathbb{H}_{0}$
is depending on only the first and the second $\mathbb{R}$-variables
$x_{1}$ and $x_{2}$. More generally, for any $u_{3},u_{4}\in\mathbb{R}\setminus\left\{ 0\right\} $,
if we define operators $D_{u_{3},u_{4}}$ and $D_{u_{3},u_{4}}^{\dagger}$
by 
\[
D_{u_{3},u_{4}}=\frac{\partial}{\partial x_{1}}+i\frac{\partial}{\partial x_{2}}+j_{0}\frac{u_{3}\partial}{\partial x_{3}}+k_{0}\frac{u_{4}\partial}{\partial x_{4}},
\]
and\hfill{}(5.10)
\[
D_{u_{3},u_{4}}^{\dagger}=\frac{\partial}{\partial x_{1}}-\frac{\partial}{\partial x_{2}}i-\frac{u_{3}\partial}{\partial x_{3}}j_{0}-\frac{u_{4}\partial}{\partial x_{4}}k_{0},
\]
as in (5.1), then, for any dilations (5.10) of $\nabla_{0}$, we have

\medskip{}

\hfill{}$D_{u_{3},u_{4}}^{\dagger}D_{u_{3},u_{4}}=\frac{\partial^{2}}{\partial x_{1}^{2}}+\frac{\partial^{2}}{\partial x_{2}^{2}},$\hfill{}(5.11)

\medskip{}

\noindent similar to (5.8), by (5.5) and (5.6). It shows that, indeed,
the Laplacian $\varDelta_{0}$ on $\mathcal{F}_{0}$ can be well-defined
to be
\[
\varDelta_{0}=\frac{\partial^{2}}{\partial x_{1}^{2}}+\frac{\partial^{2}}{\partial x_{2}^{2}},
\]
by (5.11).
\begin{defn}
Define an operator $\varDelta_{0}$ on $\mathcal{F}_{0,U}$, for any
open connected subsets $U$ of the $0$-scaled hypercomplexes $\mathbb{H}_{0}$,
by

\medskip{}

\hfill{}$\varDelta_{0}\overset{\textrm{def}}{=}\frac{\partial^{2}}{\partial x_{1}^{2}}+\frac{\partial^{2}}{\partial x_{2}^{2}}=\left(\frac{\partial^{2}}{\partial x_{1}^{2}}+\frac{\partial^{2}}{\partial x_{2}^{2}}\right)+0\cdot\left(\frac{\partial^{2}}{\partial x_{3}^{2}}+\frac{\partial^{2}}{\partial x_{4}^{2}}\right).$\hfill{}(5.12)

\medskip{}

\noindent A function $f\in\mathcal{F}_{0,U}$ is said to be left $0$(-scaled)-regular
on $U$, if 
\[
\nabla_{0}f=\frac{\partial f}{\partial x_{1}}+i\frac{\partial f}{\partial x_{2}}+j_{0}\frac{\partial f}{\partial x_{3}}+k_{0}\frac{\partial f}{\partial x_{4}}=0,
\]
and it is said to be right $0$(-scaled)-regular on $U$, if 
\[
f\nabla_{0}=\frac{\partial f}{\partial x_{1}}+\frac{\partial f}{\partial x_{2}}i+\frac{\partial f}{\partial x_{3}}j_{0}+\frac{\partial f}{\partial x_{4}}k_{0}=0.
\]
If $f\in\mathcal{F}_{0,U}$ is both left and right $0$-regular, then
$f$ is called a $0$-regular function on $U$. Also, a function $f\in\mathcal{F}_{0,U}$
is said to be $0$(-scaled)-harmonic on $U$, if
\[
\varDelta_{0}f=0,
\]
where $\varDelta_{0}$ is in the sense of (5.12).
\end{defn}

By definition, we have the following result.
\begin{thm}
Let $f\in\mathcal{F}_{0,U}$, for an open connected subset $U$ of
$\mathbb{H}_{0}$. Then

\medskip{}

\hfill{}$f$ is left $0$-regular$\;$$\Longrightarrow$$\;$$f$
is $0$-harmonic on $U$.\hfill{}$\textrm{(5.13)}$
\end{thm}

\begin{proof}
By (5.9), one has $\varDelta_{0}=\nabla_{0}^{\dagger}\nabla_{0},$
on $\mathcal{F}_{0,U}$. So, if $f$ is left $0$-regular on $U$,
then
\[
\varDelta_{0}f=\nabla_{0}^{\dagger}\left(\nabla_{0}f\right)=\nabla_{0}^{\dagger}\left(0\right)=0,
\]
and hence, it is $0$-harmonic on $U$.
\end{proof}
Similar to (5.3), we let
\[
\eta_{2}\left(w\right)=x_{2}-x_{1}i,
\]
and\hfill{}(5.14)
\[
\eta_{3}\left(w\right)=x_{3}-x_{1}j_{0},\;\mathrm{and\;}\eta_{4}\left(w\right)=x_{4}-x_{1}k_{0},
\]
in a $\mathbb{H}_{0}$-variable $w=x_{1}+x_{2}i+x_{3}j_{0}+x_{4}k_{0}$,
with $x_{1},x_{2},x_{3},x_{4}\in\mathbb{R}$. Then these functions
$\left\{ \eta_{l}\right\} _{l=2}^{4}$ are contained in $\mathcal{F}_{0,\mathbb{H}_{t}}$.
\begin{lem}
Let $\left\{ \eta_{l}\right\} _{l=2}^{4}$ be in the sense of (5.14).
Then they are not only $0$-regular, but also, $0$-harmonic on the
$0$-scaled hypercomplexes $\mathbb{H}_{0}$, i.e.,

\medskip{}

\hfill{}$\left\{ \eta_{l}\right\} _{l=2}^{4}$ are $0$-harmonic
$0$-regular functions on $\mathbb{H}_{0}$.\hfill{}$\textrm{(5.15)}$

\end{lem}

\begin{proof}
If $\eta_{2}=x_{2}-x_{1}i\in\mathcal{F}_{0,\mathbb{H}_{0}}$, then

\medskip{}

\medskip{}

$\;\;$$\nabla_{0}\eta_{2}=\eta_{2}\nabla_{0}=\frac{\partial\left(x_{2}-x_{1}i\right)}{\partial x_{1}}+i\frac{\partial\left(x_{2}-x_{1}i\right)}{\partial x_{2}}+j_{0}\frac{\partial\left(x_{2}-x_{1}i\right)}{\partial x_{3}}+k_{0}\frac{\partial\left(x_{2}-x_{1}i\right)}{\partial x_{4}}$

\medskip{}

$\;\;\;\;\;\;\;\;\;\;\;\;\;\;\;\;$$=\left(-i\right)+i\left(1\right)+j_{0}\left(0\right)+k_{0}\left(0\right)=-i+i=0$;

\medskip{}

\medskip{}

\noindent and if $\eta_{3}=x_{3}-x_{1}j_{0}\in\mathcal{F}_{0,\mathbb{H}_{0}}$,
then

\medskip{}

\medskip{}

$\;\;$$\nabla_{0}\eta_{3}=\eta_{3}\nabla_{0}=\frac{\partial\left(x_{3}-x_{1}j_{0}\right)}{\partial x_{1}}+i\frac{\partial\left(x_{3}-x_{1}j_{0}\right)}{\partial x_{2}}+j_{0}\frac{\partial\left(x_{3}-x_{1}j_{0}\right)}{\partial x_{3}}+k_{0}\frac{\partial\left(x_{3}-x_{1}j_{0}\right)}{\partial x_{4}}$

\medskip{}

$\;\;\;\;\;\;\;\;\;\;\;\;\;\;\;\;$$=\left(-j_{0}\right)+i\left(0\right)+j_{0}\left(1\right)+k_{0}\left(0\right)=-j_{0}+j_{0}=0$;

\medskip{}

\medskip{}

\noindent and if $\eta_{4}=x_{4}-x_{1}k_{0}\in\mathcal{F}_{0,\mathbb{H}_{0}}$,
then

\medskip{}

\medskip{}

$\;\;$$\nabla_{0}\eta_{3}=\eta_{3}\nabla_{0}=\frac{\partial\left(x_{4}-x_{1}k_{0}\right)}{\partial x_{1}}+i\frac{\partial\left(x_{4}-x_{1}k_{0}\right)}{\partial x_{2}}+j_{0}\frac{\partial\left(x_{4}-x_{1}k_{0}\right)}{\partial x_{3}}+k_{0}\frac{\partial\left(x_{4}-x_{1}k_{0}\right)}{\partial x_{4}}$

\medskip{}

$\;\;\;\;\;\;\;\;\;\;\;\;\;\;\;\;$$=\left(-k_{0}\right)+i\left(0\right)+j_{0}\left(0\right)+k_{0}\left(1\right)=-k_{0}+k_{0}=0$.

\medskip{}

\medskip{}

\noindent Therefore, the functions $\left\{ \eta_{l}\right\} _{l=2}^{4}$
of (5.14) are $0$-regular on $\mathbb{H}_{0}$. Therefore, these
0-regular functions $\left\{ \eta_{l}\right\} _{l=2}^{4}$ are $0$-harmonic
by (5.13).
\end{proof}
By (5.15), we consider the symmetrized product $\eta^{\mathbf{n}}$,
for $\mathbf{n}=\left(n_{1},n_{2},n_{3}\right)\in\mathbb{N}_{0}^{3}$,
defined by

\medskip{}

\hfill{}$\eta^{\mathbf{n}}=\frac{1}{\mathbf{n}!}\left(\eta_{2}^{\left(n_{1}\right)}\times\eta_{3}^{\left(n_{2}\right)}\times\eta_{4}^{\left(n_{3}\right)}\right)\in\mathcal{F}_{0,\mathbb{H}_{0}}$,\hfill{}(5.16)

\medskip{}

\noindent like (5.4), as in Section 4.5, where ($\times$) is the
symmetrized product (4.5.1), expressed to be (4.5.2).
\begin{thm}
The symmetrized products $\eta^{\mathbf{n}}\in\mathcal{F}_{0,\mathbb{H}_{0}}$
of (5.16) are $0$-harmonic $0$-regular functions on $\mathbb{H}_{0}$,
for all $\mathbf{n}\in\mathbb{N}_{0}^{3}$, i.e.,

\medskip{}

\hfill{}$\eta^{\mathbf{n}}$ are $0$-harmonic $0$-regular functions
on $\mathbb{H}_{0}$, $\forall\mathbf{n}\in\mathbb{N}_{0}^{3}$.\hfill{}$\textrm{(5.17)}$
\end{thm}

\begin{proof}
The proof of the $0$-regularity of $\eta^{\mathbf{n}}$ are similar
to that of (4.5.11). In particular, if we let
\[
e_{1}=-i,\;e_{2}=-j_{0},\;\mathrm{and\;}e_{3}=-k_{0},
\]
and if we replace $\nabla_{t}$ to $\nabla_{0}$, then the formula
(4.5.8) is obtained similarly, and hence, the computation (4.5.10)
holds by (4.5.9). i.e.,
\[
x_{1}\mathbf{n}!\nabla_{0}\eta^{\mathbf{n}}=0,\Longrightarrow\nabla_{0}\eta^{\mathbf{n}}=0.
\]
So, the symmetrized products $\eta^{\mathbf{n}}$ are $0$-regular,
for $\mathbf{n}\in\mathbb{N}_{0}^{3}$. By the $0$-regularity, $\eta^{\mathbf{n}}$
is 0-harmonic on $\mathbb{H}_{0}$, too, by (5.13). Therefore, the
relation (5.17) holds true.
\end{proof}
The above theorem shows that all symmetrized products $\eta^{\mathbf{n}}$
of (5.16) are 0-harmonic 0-regular functions on $\mathbb{H}_{0}$
by (5.17).

Now, consider that if $\left\{ \eta_{2},\eta_{3},\eta_{4}\right\} \subset\mathcal{F}_{0,\mathbb{H}_{0}}$
are the 0-harmonic 0-regular functions (5.14), then they are understood
to be their images,
\[
\eta_{2}=\left(x_{2}-x_{1}i,0\right),\;\eta_{3}=\left(x_{3},-x_{1}\right),\;\eta_{4}=\left(x_{4},-x_{1}\right),
\]
as elements of the 0-scaled hypercomplex ring $\mathbb{H}_{0}$, having
their norms,

\[
\left\Vert \eta_{2}\right\Vert _{0}=\sqrt{\left|\left|x_{2}-x_{1}i\right|^{2}-0\left|0\right|^{2}\right|}=\sqrt{x_{2}^{2}+\left(-x_{1}\right)^{2}}=\sqrt{x_{2}^{2}+x_{1}^{2}},
\]
\[
\left\Vert \eta_{3}\right\Vert _{0}=\sqrt{\left|\left|x_{3}\right|^{2}-0\left|-x_{1}\right|^{2}\right|}=\sqrt{\left|x_{3}^{2}\right|}=\left|x_{3}\right|,
\]
and\hfill{}(5.18)
\[
\left\Vert \eta_{4}\right\Vert _{0}=\sqrt{\left|\left|x_{4}\right|^{2}-0\left|-x_{1}\right|^{2}\right|}=\sqrt{\left|x_{4}^{2}\right|}=\left|x_{4}\right|.
\]

\begin{lem}
Let $\eta^{\mathbf{n}}\in\mathcal{F}_{0,\mathbb{H}_{t}}$ be a function
$\mathrm{(5.16)}$ for $\mathbf{n}=\left(n_{1},n_{2},n_{3}\right)\in\mathbb{N}_{0}^{3}$.
Then

\medskip{}

\hfill{}$\left\Vert \eta^{\mathbf{n}}\right\Vert _{0}\leq\left(\sqrt{x_{2}^{2}+x_{1}^{2}}^{n_{1}}\right)\left(\left|x_{3}\right|^{n_{2}}\right)\left(\left|x_{4}\right|^{n_{3}}\right)$.\hfill{}$\mathrm{(5.19)}$
\end{lem}

\begin{proof}
By (5.16), for any $\mathbf{n}=\left(n_{1},n_{2},n_{3}\right)\in\mathbb{N}_{0}^{3}$,

\medskip{}

$\;\;$$\left\Vert \eta^{\mathbf{n}}\right\Vert _{0}=\left\Vert \frac{1}{\mathbf{n}!}\left(\eta_{2}^{(n_{1})}\times\eta_{3}^{(n_{2})}\times\eta_{4}^{(n_{3})}\right)\right\Vert _{0}$

\medskip{}

$\;\;\;\;$$=\left\Vert \frac{1}{\left(n_{1}!\right)\left(n_{2}!\right)\left(n_{3}!\right)}\left(\frac{\left(n_{1}!\right)\left(n_{2}!\right)\left(n_{3}!\right)}{n!}\underset{\sigma\in S_{n}}{\sum}h_{\sigma\left(1\right)}h_{\sigma\left(2\right)}...h_{\sigma\left(n\right)}\right)\right\Vert _{0}$

\medskip{}

\noindent where $h_{\sigma\left(l\right)}\in\left\{ \eta_{2},\eta_{3},\eta_{4}\right\} $,
for all $\sigma\in S_{n}$, with $n=n_{1}+n_{2}+n_{3}\in\mathbb{N}$

\medskip{}

$\;\;\;\;$$=\left\Vert \eta_{2}^{n_{1}}\eta_{3}^{n_{2}}\eta_{4}^{n_{3}}\right\Vert _{0}\leq\left\Vert \eta_{2}\right\Vert _{0}^{n_{1}}\left\Vert \eta_{3}\right\Vert _{0}^{n_{2}}\left\Vert \eta_{4}\right\Vert _{0}^{n_{3}}$

\medskip{}

\noindent by (4.5.5)

\medskip{}

$\;\;\;\;$$=\left|x_{2}\right|^{n_{1}}\left|x_{3}\right|^{n_{2}}\left|x_{4}\right|^{n_{3}},$

\medskip{}

\noindent by (5.18). Therefore, the inequality (5.19) holds.
\end{proof}
The above lemma shows that, for any arbitrarily fixed $\mathbf{n}\in\mathbb{N}_{0}^{3}$,
the corresponding function $\eta^{\mathbf{n}}$ is bounded by (5.19).
\begin{thm}
Let $f\in\mathcal{F}_{0,U}$ be a $\mathbb{R}$-differentiable function,
where $U$ is an open connected subset of $\mathbb{H}_{0}$ containing
$0\in\mathbb{H}_{0}$. If $f$ is $\mathbb{R}$-analytic at $0\in U$,
then 
\[
f\textrm{ is left }0\textrm{-regular on }U,
\]
if and only if\hfill{}$\mathrm{(5.20)}$
\[
f=f\left(0\right)+\underset{\mathbf{n}\in\mathbb{N}^{3}}{\sum}\eta^{\mathbf{n}}f_{\mathbf{n}},
\]
with
\[
f_{\mathbf{n}}=\frac{1}{\mathbf{n}!}\frac{\partial^{n_{1}+n_{2}+n_{3}}f}{\partial x_{2}^{n_{1}}\partial x_{3}^{n_{2}}\partial x_{4}^{n_{3}}}\left(0\right),\;\forall\mathbf{n}\in\mathbb{N}^{3}.
\]
\end{thm}

\begin{proof}
If $f=f\left(0\right)+\underset{\mathbf{n}\in\mathbb{N}^{3}}{\sum}\eta^{\mathbf{n}}f_{\mathbf{n}}$,
then it is left $0$-regular, by the $0$-harmonic-$0$-regularity
(5.17) of $\left\{ \eta^{\mathbf{n}}:\mathbf{n}\in\mathbb{N}_{0}^{3}\right\} $.

Suppose $f$ is left $0$-regular in $\mathcal{F}_{0,U}$, satisfying
$\nabla_{0}f=0$. Then
\[
f\left(w\right)-f\left(0\right)=\overset{4}{\underset{n=2}{\sum}}\left(\eta_{n}\left(w\right)\right)\left(\left(R_{n}f\right)\left(w\right)\right),
\]
where\hfill{}(5.21)
\[
\left(R_{n}f\right)\left(w\right)=\int_{0}^{1}\frac{\partial f\left(tw\right)}{\partial x_{n}}dt,\;\forall n=2,3,4,
\]
in the $\mathbb{H}_{t}$-variable $w=x_{1}+x_{2}i+x_{3}j_{t}+x_{4}k_{t}$,
with $x_{1},x_{2},x_{3},x_{4}\in\mathbb{R}$. By iterating (5.21),
one can get the function in (5.20), by the boundedness condition (5.19).
\end{proof}

{\bf Acknowledgements and funding:}  Daniel Alpay thanks
the Foster G. and Mary McGaw Professorship in Mathematical Sciences, which
supported this research. No other funds, grants, or other support were received
during the preparation of this manuscript.\\

{\bf Data Availability Statement:} The authors confirm that no data known is used
in our manuscript.

{\bf Statements and Declarations of Competing Interests:} We declare that we have no competing interests.

\end{document}